\documentclass[11pt]{amsart}
\usepackage[mathscr]{euscript}
\usepackage{amssymb}

\usepackage{amscd}
\usepackage{amsmath, amssymb}
\usepackage{amsfonts}
\usepackage[colorlinks,linkcolor=blue,citecolor=blue, pdfstartview=FitH]{hyperref}
\usepackage[all]{xy}
\usepackage{diagbox}
\numberwithin{equation}{section}

\theoremstyle{plain}
\newtheorem{thm}{Theorem}[section]
\newtheorem{theorem}[thm]{Theorem}
\newtheorem{lemma}[thm]{Lemma}
\newtheorem{corollary}[thm]{Corollary}
\newtheorem{proposition}[thm]{Proposition}

\theoremstyle{definition}

\newtheorem{question}[thm]{Question}

\newtheorem{remark}[thm]{Remark}

\newtheorem{definition}[thm]{Definition}

\newtheorem{example}[thm]{Example}

\newtheorem{defn-thm}[thm]{Definition-Theorem}

\def\wt{\widetilde}
\newcommand{\im}{{ \mathrm{im}\,}}

\newcommand{\C}{{ \mathbb{C} }}
\newcommand{\Z}{{ \mathbb{Z} }}

\newcommand{\g}{{ \mathfrak{g} }}

\newcommand{\p}{{ \partial }}
\newcommand{\pb}{{ \bar{\partial} }}
\newcommand{\ot}{{\,\otimes\,}}

\usepackage{fancyhdr}
\makeatletter
\@namedef{subjclassname@2020}{\textup{2020} Mathematics Subject Classification}
\makeatother

\begin{document}
\title{Deformations of $(p,q)$-forms and degenerations of the Fr\"olicher spectral sequence}

\author{Xueyuan Wan}
\address{Xueyuan Wan, Mathematical Science Research Center, Chongqing University of Technology, Chongqing 400054, China} \email{xwan@cqut.edu.cn}

\author{Wei Xia}
\address{Wei Xia, Mathematical Science Research Center, Chongqing University of Technology, Chongqing 400054, China} \email{xiaweiwei3@126.com, xiawei@cqut.edu.cn}


\thanks{Xueyuan Wan is partially supported by the National Natural Science Foundation of China (Grant No. 12101093) and the Natural Science Foundation of Chongqing (Grant No. CSTB2022NSCQ-JQX0008), the Scientific Research Foundation of the Chongqing University of Technology. Wei Xia is supported by the National Natural Science Foundation of China No. 11901590, Natural Science Foundation of Chongqing (China) No. CSTB2022NSCQ-MSX0876 and Scientific Research Foundation of Chongqing University of Technology.}
\date{\today}

\begin{abstract}
It is well-known that Hodge numbers are invariant under deformations of complex structures if the Fr\"olicher spectral sequence of the central fiber degenerates at the first page (i.e. $E_1=E_\infty$). As a result, the deformations of $(p,q)$-forms are unobstructed for all $(p,q)$ if $E_1=E_\infty$. We refine this classical result by showing that for any fixed $(p,q)$ the deformations of $(p,q)$-forms are unobstructed if the differentials $d_r^{p,q}$ in the Fr\"olicher spectral sequence satisfy
  \[
  \bigoplus_{r\geq 1}d_r^{p,q}=0\quad\text{and}\quad
  \bigoplus_{ r\geq i\geq 1 } d_r^{p-i,q+i}=0.
  \]
Moreover, the deformation stability of the degeneration property for Fr\"olicher spectral sequences in the first page and higher pages is also studied. In particular, we have found suitable conditions to ensure the deformation stability of $E_r^{p,q}=E_\infty^{p,q}$ ($r\geq1$) for fixed $(p,q)$.
\end{abstract}

\maketitle

\tableofcontents
\section{Introduction}
As one of the most basic holomorphic invariants of complex manifolds, Hodge numbers are studied extensively in complex geometry, algebraic geometry, and related subjects. For a smooth family of compact K\"ahler manifolds, it is well-known that Hodge numbers do not vary as the complex structure deforms which is an easy corollary of Kodaira-Spencer's upper semi-continuous theorem (see e.g. \cite[pp.\,46]{BDIP02}). The K\"ahler assumption may be greatly weakened.

Indeed, consider a complex analytic family $\pi: \mathcal{X}\to \Delta$ over a small disc $\Delta\subset \C$, if the Fr\"olicher spectral sequence of the central fiber $X_0$ degenerate at the first page (denoted by $E_1=E_\infty$), then the Hodge numbers $h^{p,q}(X_t)$ of the fibers $X_t:=\pi^{-1}(t)$ does not depend on the parameter $t\in\Delta$. Furthermore, for any $t\in\Delta$ the Fr\"olicher spectral sequence of $X_t$ also degenerates at the first page (see \cite[Coro.\,2.15]{Bin83} or \cite[Prop.\,9.20]{Voi02I}).

Although $E_1=E_\infty$ seems a reasonable sufficient condition for the deformation invariance of Hodge numbers, there are many concrete examples of complex analytic family (see \cite{Nak75,Ang13,AK17b,AK17a} and the references therein) where $E_1\neq E_\infty$ and $h^{p,q}$ are deformation invariant only for some particular choices of $(p,q)$. Taking the Kuranishi family of the Iwasawa manifold (complex dimension equal to $3$) as an example, $h^{p,q}$ are deformation invariant if and only if $(p,q)=(0,0), (0,1),(0,2),(3,0),(0,3),(3,1),(3,2),(3,3)$.

On the other hand, it is interesting to notice that the deformation invariance of $h^{0,1}$ is often involved in studies related to deformation limits of projective manifolds or Moishezon manifolds, see \cite{Pop13,Bar15,RT21}. We are therefore led to the following natural question:
\begin{question}\label{que-inv-hpq} Given fixed $(p,q)$ and a complex analytic family $\pi: \mathcal{X}\to \Delta$, under what conditions on the central fiber $X_0$ (preferably weaker than $E_1=E_\infty$), will the Hodge number $h^{p,q}(X_t)$ be independent of $t$?
\end{question}
This is in fact part of a question posed by Rao-Zhao \cite[Question\,1.6]{RZ18} where they are seeking sufficient and necessary conditions that make $h^{p,q}$ invariant under deformations of complex structures when $(p,q)$ is fixed and arbitrary, respectively. In a recent work by the second named author, we have managed to prove the following criterion (\cite[Thm\,1.2]{Xia19dDol} in the case of $E=\Omega^p$) regarding this question:
\begin{theorem}\label{thm-invariance-unob-0} Given fixed $(p,q)$ and a complex analytic family $\pi: \mathcal{X}\to \Delta$, the Hodge number $h^{p,q}(X_t)$ is independent of $t$ if and only if the deformations of $(p,q)$ and $(p,q-1)$-forms on $X_0$ are canonically unobstructed.
\end{theorem}
For the readers' convenience, let us briefly recall the definition of canonically unobstructedness in the case of $(p,q)$-forms (see \cite[Def.\,5.1]{Xia19dDol} for more information). In the situation of Theorem \ref{thm-invariance-unob-0}, let $\phi(t)\in A^{0,1}(X,T_X^{1,0})$ be the Beltrami differential \cite{MK71,Kod86} which represents the complex structures of the fibers $X_t:=\pi^{-1}(t)$ where $t\in\Delta$.
\begin{definition}\label{def-defor-pq-forms-0}
Given
\[
\sigma_0\in \ker \pb\cap A^{p,q}(X).
\]
Suppose there exists an analytic subset $T\subseteq \Delta$ containing the origin $0\in\Delta$, and a family of $(p,q)$-forms
\[
\{\sigma(t)\in A^{p,q}(X)\}_{t\in T}
\]
such that
\begin{enumerate}
  \item $\sigma (t)$ is holomorphic in $t$ and $\sigma (0) = \sigma_0$;
  \item $\pb_{\phi(t)}\sigma (t)=0$, for any $t\in T$,
\end{enumerate}
where $\pb_{\phi(t)}=\pb-\mathcal{L}^{1,0}_{\phi(t)}$ and $\mathcal{L}^{1,0}_{\phi(t)}=i_{\phi(t)}\p-\p i_{\phi(t)}$.

In this case, we call $\sigma (t)$ a (Dolbeault) \textit{deformation} of $\sigma_0$ (with respect to $\pi$) on $T$. Assume $X$ has been equipped with a fixed Hermitian metric. A deformation $\{\sigma(t)\in A^{p,q}(X)\}_{t\in T}$ of $\sigma_0$ is said to be \textit{canonical} if it satisfies
$\sigma(t) = \sigma_0 + \bar{\partial}^*G \mathcal{L}^{1,0}_{\phi(t)}\sigma(t)$ for any $t\in T$ where $G$ is the $\bar{\partial}$-Green's operator. We say that the (Dolbeault) deformation of $\sigma_0$ is \textit{canonically unobstructed} if there is a canonical deformation of $\sigma_0$ on $\Delta$. Finally, the (Dolbeault) deformations of $(p,q)$-forms are \textit{canonically unobstructed} if for every class in the Dolbeault cohomology $H_\pb^{p,q}(X)$, it has a representative which is canonically unobstructed. In this sense, the notion of canonically unobstructedness is independent of the choices of the Hermitian metric \cite[Coro\,4.8]{HX24A}.
\end{definition}
By Theorem \ref{thm-invariance-unob-0}, Question \ref{que-inv-hpq} can be reduced to the following
\begin{question}\label{que-deformation-unob-pq-0} Given fixed $(p,q)$ and a complex analytic family $\pi: \mathcal{X}\to \Delta$, under what conditions on the central fiber $X_0$ (preferably weaker than $E_1=E_\infty$), the deformations of $(p,q)$-forms are canonically unobstructed?
\end{question}
In this paper, we will provide an answer to Question \ref{que-deformation-unob-pq-0} which shows in order that the deformations of $(p,q)$-forms are canonically unobstructed, only less than one-quarter of the classical condition $E_1=E_\infty$ will be already enough.

It is well known that $E_1=E_\infty$ is a deformation open property in the sense of Popovici \cite{Pop14}. This follows easily from the fact that Hodge numbers vary upper semi-continuously in a complex analytic family. When $r>1$, the situation is very different. In fact, Ceballos-Otal-Ugarte-Villacampa find that 1. $E_2=E_\infty$ is not a deformation open property in general \cite[Coro.\,5.12]{COUV16};
2. the varying of $\dim E_r^{\bullet,\bullet}$ ($r>1$) in a complex analytic family may be neither upper nor lower semi-continuous, see \cite[Coro.\,4.9]{COUV16}. Nevertheless, Popovici recently posed the following
\begin{question}\cite[Pro.\,1.1]{Pop24} Let $r\geq 2$ be an integer. Find a geometric property $(P)$ that certain compact complex manifolds satisfy such that, whenever $X_0$ has property $(P)$ and $E_r=E_\infty$, then $E_r=E_\infty$ also holds for $X_t$ with $t\in\Delta$ close enough to $0$.
\end{question}
Related to this question, Maschio \cite{Mas20} have shown that for any fixed $k$, assume for every $(p,q)\in\mathbb{N}\times\mathbb{N}$ with $p+q=k$ the Hodge number $h^{p,q}$ is deformation invariant and $E_2^{p,q}(X)=E_\infty^{p,q}(X)$, then for any $t\in\Delta$, $E_2^{p,q}(X_t)=E_\infty^{p,q}(X_t)$ whenever $p+q=k$. In the same spirit of studying unobstructed deformations of $(p,q)$-forms, we will be much interested in the case of fixed $(p,q)$. More precisely, we will provide suitable conditions (these are motivated by \cite{Mas20}) to ensure the deformation openness of $E_r^{p,q}=E_\infty^{p,q}$ ($r\geq1$) for fixed $(p,q)$.
\subsection{Main results}
\subsubsection{Unobstructed deformations of $(p,q)$-forms}
Recall that in the Fr\"olicher spectral sequence of any complexes manifold $X$, there is a sequence of differential cochain complexes $(E_r^{\bullet,\bullet}(X),d_r^{\bullet,\bullet}(X))$ with differentials
\[
d_r^{p,q}(X):E_r^{p,q}(X)\longrightarrow E_r^{p+r,q-r+1}(X),~r=0,1,2,\cdots,
\]
where $(p,q)\in\mathbb{N}\times\mathbb{N}$. Our first main result is the following:
\begin{thm}[=Theorem \ref{mainthm}+Corollary \ref{coro-altsum-hp,q}]\label{mainthm-0}
Let $\pi: (\mathcal{X}, X)\to (\Delta,0)$ be a complex analytic family over a small disc $\Delta\subset \C$. Assume
  \[
  \bigoplus_{r\geq 1}d_r^{p,q}(X)=0\quad\text{and}\quad
  \bigoplus_{ r\geq i\geq 1 }d_r^{p-i,q+i}(X)=0.
  \]
Then
\begin{enumerate}
  \item The Dolbeault deformations of $(p,q)$-forms on $X$ are canonically unobstructed or equivalently, the alternating sum of Hodge numbers: $\sum_{i=0}^{q}(-1)^{q-i}h^{p,i}$, is deformation invariant.
  \item For any $t\in\Delta$ we have
  \[
     \bigoplus_{r\geq 1}d_r^{p,q}(X_t)=0.
  \]
\end{enumerate}
\end{thm}

The condition $\bigoplus_{r\geq i\geq 1}d_r^{p-i,q+i}(X)=0$ is equivalent to the statement that
\[
F^{p}A^{p+q+1}(X)\cap dA^{p+q}(X)=dF^{p}A^{p+q}(X)
\]
holds, see Proposition \ref{coro-Fp+1capdA=dFp+1}. Moreover, we remark that Theorem \ref{mainthm-0} is a refinement of the classical fact:
\begin{equation}\label{eq-classical-fact-a0}
\begin{split}
&E_1^{p-i,q+i}(X)\cong E_\infty^{p-i,q+i}(X)~\text{for~any}~i\in \mathbb{Z}\\
\Longrightarrow &~\text{deformation~invariance~of}~h^{p-i,q+i}~\text{for~any}~i\in \mathbb{Z}.
\end{split}
\end{equation}
Indeed, in view of Theorem \ref{thm-invariance-unob-0}, statement \eqref{eq-classical-fact-a0} is equivalent to the following
\begin{equation}\label{eq-classical-fact-b0}
\begin{split}
&E_1^{p-i,q+i}(X)\cong E_\infty^{p-i,q+i}(X)~\text{for~any}~i\in \mathbb{Z}\\
\Longrightarrow&~\text{canonically~unobstructed~deformations~for}\\
&(p-i,q+i), (p-i,q-1+i)\text{-forms~with}~i\in \mathbb{Z}.
\end{split}
\end{equation}
On the other hand, there holds the following (see Proposition \ref{prop-Er=Er+1}),
\begin{align*}
E_1^{p-i,q+i}(X)\cong E_\infty^{p-i,q+i}(X)~\text{for~any}~i\in \mathbb{Z}
&\Longleftrightarrow \bigoplus_{\substack{i\in \mathbb{Z},\\r\geq 1 }}d_r^{p-i,q+i}(X)=0=\bigoplus_{\substack{i\in \mathbb{Z},\\r\geq 1 }}d_r^{p-r-i,q+r-1+i}(X)\\
&\Longleftrightarrow \bigoplus_{\substack{i\in \mathbb{Z},\\r\geq 1 }}d_r^{p-i,q+i}(X)=0=\bigoplus_{\substack{i\in \mathbb{Z},\\r\geq 1 }}d_r^{p-i,q-1+i}(X).
\end{align*}

In particular, we obtain the following answer to Question \ref{que-inv-hpq}.
\begin{corollary}[=Corollary \ref{coro-def-invariance-hpq}]
Assume
  \[
  \bigoplus_{\substack{j=0,1\\ r\geq 1}}d_r^{p,q-j}(X)=0=
  \bigoplus_{\substack{j=0,1\\ r\geq i\geq 1}}d_r^{p-i,q+i-j}(X).
  \]
Then the Hodge number $h^{p,q}$ is deformation invariant.
\end{corollary}

In the case that $\pi$ is the Kuranishi family of the Iwasawa manifold (see Example \ref{example Case III-(2)}), we find Theorem \ref{mainthm-0} is applicable for
\[
(p,q)=(0,1),(0,2),(0,3),(3,0),(3,1),(3,2),(2,3),(3,3).
\]
It is asked in \cite{Xia19dDol} whether $d_1^{p,q}(X)=0=d_1^{p-1,q+1}(X)$ is enough for the unobstructed deformations of $(p,q)$-forms? It seems that no counter example is known to this question. For necessity, concrete computations in Section \ref{sec-example} (including the Iwasawa manifolds) show that Theorem \ref{mainthm-0} is not optimal. Indeed, in Example \ref{example-E3} we exhibit a compact nilmanifold with nilpotent complex structure $X$ such that
\[
d_r^{p,q}\neq0\Longleftrightarrow r=2~\text{and}~(p,q)=(0,2), (1,2),
\]
while the deformations of $(p,q)$-forms are canonically unobstructed for all $(p,q)\in\mathbb{N}\times\mathbb{N}$.
\begin{corollary}[=Corollary \ref{coro-h0,q}]
The deformations of $(0,q)$-forms on $X$ are canonically unobstructed if $\bigoplus_{\substack{r\geq 1}}d_r^{0,q}(X)=0$. In particular, $h^{0,1}$ is deformation invariant if $\bigoplus_{\substack{r\geq 1}}d_r^{0,1}(X)=0$.
\end{corollary}

\subsubsection{Deformation stability of $E_r^{p,q}=E_\infty^{p,q}$ for $r\geq1$ and fixed $(p,q)$}
Our second main result (see Theorem \ref{mainthm-2nd}) studies the deformation stability of the condition $\bigoplus_{\lambda\geq r}d_\lambda^{p,q}=0$ for general $r\geq1$. To avoid technical notions, we only state a simplified version of it for $r=2$ as the following:
\begin{thm}[=Corollary \ref{coro-to-mainthm2}]\label{mainthm-2nd-restate-0}
Let $\pi: (\mathcal{X}, X)\to (\Delta,0)$ be a complex analytic family over a small disc $\Delta\subset \C$.
\begin{enumerate}
\item Assume
\begin{itemize}
  \item the deformations of $(p,q)$-forms on $X$ are canonically unobstructed or equivalently, the alternating sum of Hodge numbers $\sum_{i=0}^{q}(-1)^{q-i}h^{p,i}$ is deformation invariant;
  \item there holds
\begin{equation}
\bigoplus_{\lambda\geq i\geq 1}d_\lambda^{p+1-i,q-1+i}(X)=0.
\end{equation}
\end{itemize}
Then
\[
\bigoplus_{\lambda\geq 1}d_\lambda^{p,q}(X_t)=0,\quad \text{for~any}~t\in\Delta.
\]
\item Assume
\begin{itemize}
  \item the deformations of $(p+1,q)$-forms, $(p+1,q-1)$-forms and $2$-filtered $(p,q)$-forms on $X$ are canonically unobstructed;
  \item there holds
\begin{equation}
\bigoplus_{ \lambda\geq i\geq 1 }d_\lambda^{p+2-i,q-2+i}(X)=0.
\end{equation}
\end{itemize}
Then
\[
\bigoplus_{\lambda\geq 2}d_\lambda^{p,q}(X_t)=0,\quad \text{for~any}~t\in\Delta.
\]
\end{enumerate}
\end{thm}
For any $r\geq1$, \textit{$r$-filtered $(p,q)$-forms} on a complex manifold $X$ are defined to be elements in the vector space $\bigoplus_{i=0}^{r-1}A^{p+i,q-i}(X)$. Note that $1$-filtered $(p,q)$-forms are just $(p,q)$-forms in the usual sense, so $r$-filtered $(p,q)$-forms are natural generalizations of $(p,q)$-forms. Canonically unobstructedness for deformations of $r$-filtered $(p,q)$-forms can be defined in a similar way as in Definition \ref{def-defor-pq-forms-0}. To avoid redundance, we refer the readers to Definition \ref{def-cano-unobst}.

Finally, we can show a generalized version of $(1)$ in Theorem \ref{mainthm-0}, which deals with deformation of $r$-filtered forms, also holds:
\begin{corollary}[=Corollary \ref{coro-deformation of r-filtered-2}]\label{coro-deformation of r-filtered-2-0}
Let $\pi: (\mathcal{X}, X)\to (\Delta,0)$ be a complex analytic family over a small disc $\Delta\subset \C$. Assume for some $r\geq1$, there holds
\begin{equation}
\bigoplus_{ \substack{\lambda\geq 1\\ 1\leq i\leq r+1}} d_\lambda^{p+r-i,q-r+i}(X)=0=\bigoplus_{\lambda\geq i\geq 2}d_\lambda^{p-i,q+i}(X).
\end{equation}
Then the deformation of $r$-filtered $(p,q)$-forms on $X$ are canonically unobstructed.
\end{corollary}

\subsection{Related works in the literature and ideas of proof}
As has been mentioned, the problem of deformation stability of $E_r^{p,q}=E_\infty^{p,q}$ for $r=2$ was studied by Maschio \cite{Mas20}. Besides of this, it was shown by Stelzig \cite{Ste22} that for $q=1$ or $n-1$, $\dim E^{0,q}_1/E^{0,q}_\infty$ is upper semi-continuous in deformations of complex structures. As a result,
\[
\bigoplus_{\lambda\geq 1}d_\lambda^{0,q}(X)=0\Longrightarrow\bigoplus_{\lambda\geq 1}d_\lambda^{0,q}(X_t)=0,\quad~\text{for~any}~t\in\Delta,
\]
which can also be implied by Theorem \ref{mainthm-0}. By using Grauert’s direct image theorem, Ye \cite[Coro.\,3.7]{Ye08} showed that $h^{p,q}$ is deformation invariant if for any finite order infinitesimal deformation of $\pi$, the $d_1$ of the Fr\"olicher spectral sequence vanishes at $(p,q)$-th and $(p-1,q+1)$-th position. We want to point out that the results of Maschio and Stelzig are intimately related to the upper semi-continuous theorem for elliptic differential operators. In fact, Maschio proves his result by showing $\dim E_2^{p,q}$ is upper semi-continuous which is modeled on the classical proof of upper semi-continuous theorem given by Kodaira-Spencer, while Stelzig's proof used the ellipticity of the Schweitzer complex and the structure theorem of double complex~\cite{Ste21,KQ20}. Our approach is very different.

The proof of Theorem \ref{mainthm-0} and Theorem \ref{mainthm-2nd-restate-0} is based on the power series method which was initiated in the work of Liu-Sun-Yau \cite{LSY09}. This method and many of its subsequent developments (see \cite{LZ20} for a recent survey) has its origins in Kodaira-Spencer's classical work on deformations of complex manifolds (see e.g. \cite[Ch.\,4]{MK71}) and has been applied to various problems related to deformation invariance or stability. For example, these include the following:
\begin{enumerate}
\item deformation invariance of Hodge numbers \cite{RZ15,RZ18,Xia19dDol,WZ20,RZ22}, Bott-Chern numbers \cite{RWZ19,Xia19dBC} and plurigenera \cite{LZ18};
\item deformation stability of special metrics \cite{RWZ19,RWZ21}, Fano K\"ahler-Einstein manifolds \cite{Sun12,CSYZ22,FSZ22};
\item Taylor expansions of period mappings \cite{LRY15,Yin10,ZR13,ZR15}.
\end{enumerate}
Compared to previous works mentioned above, there are several special features in this work:
\begin{itemize}
\item Assumptions: conditions like those in Theorem \ref{mainthm-0} seem not appear elsewhere before. In fact, almost all previous works in this direction assume some form of the $\p\pb$-lemma (see e.g. \cite{LRY15,RZ18,RWZ19,Xia19dDol,WZ20}). By\footnote{Similar ideas has been adopted in the work of K. Liu, S. Rao, and the first named author \cite{LRW19} where they presented geometric and simpler proofs of Deligne's degeneracy theorem for the logarithmic Hodge to de Rham spectral sequences at $E_1$-level.} using the explicit descriptions of the Fr\"olicher spectral sequence $(E_r^{p,q}(X),d_r(X))$ provided by Cordero-Fern\'andez-Ugarte-Gray \cite{CFUG97}, our assumptions are readily translated into statements that we use in practice. For example, the condition $\bigoplus_{r\geq 1}d_r^{p,q}(X)=0$ in Theorem \ref{mainthm-0} is equivalent to the statement that for any given $\alpha^{p,q}\in \ker\pb\cap A^{p-i,q+i}(X)$ there exist $\alpha\in \ker d\cap F^{p}A^{p+q}(X)$ such that $\Pi^{p,q}\alpha$, the $(p,q)$-th component of $\alpha$, equals to $\alpha^{p,q}$. These discussions will be collected in Appendix \ref{appendix-B}.
  \item Extension equation: the extension equation for deformations of $(p,q)$-forms is by definition the following
  \begin{equation}\label{eq-extsion-pb-0}
  \pb\alpha^{p,q}(t)=\mathcal{L}_{\phi(t)}^{1,0}\alpha^{p,q}(t),
  \end{equation}
  where $\alpha^{p,q}(t)$ is a power series of $t$ with coefficients in $A^{p,q}(X)$. Instead of solving \eqref{eq-extsion-pb-0} directly which is equivalent to a linear system of $\pb$-equations, we try to solve
  \begin{equation}\label{eq-extsion-d-0}
  (\pb_{\phi(t)}+\p)\alpha(t)=0\quad\text{or}\quad d_{\phi(t)}\alpha(t)=0,
  \end{equation}
  for $\alpha(t)=\alpha^{p,q}(t)+\alpha^{p+1,q-1}(t)+\cdots\in F^pA^{p+q}(X)$ where $d_{\phi(t)}:=e^{-i_{\phi(t)}}de^{i_{\phi(t)}}=\pb_{\phi(t)}+\p$. Notice that \eqref{eq-extsion-pb-0} is exactly the $(p,q+1)$-component of \eqref{eq-extsion-d-0}. Decomposing \eqref{eq-extsion-d-0} according to degrees, we need to solve a system of $d$-equations:
  \begin{equation}\label{eq-extsion-d-expanded-0}
  d\alpha_k=\sum_{1\leq j\leq k}\mathcal{L}_{\phi_j}^{1,0}\alpha_{k-j},~\alpha_k\in F^{p}A^{p+q}(X),\quad  k\geq 1.
  \end{equation}
  \item Iteration process: when trying to solve \eqref{eq-extsion-d-expanded-0} inductively, we find that it seems not possible to show directly the obstruction forms (terms in the right hand side of \eqref{eq-extsion-d-expanded-0}) are $d$-exact. This is very different from the situation we are familiar with previously, see e.g. \cite{RZ18,Xia19dDol}. In fact, here we need a key lemma (i.e. Lemma \ref{lem-induction}) to overcome this difficulty. Lemma \ref{lem-induction} basically says if \eqref{eq-extsion-d-expanded-0} can be solved for $k\leq N$, then it can also be ``solved" when $k=N+1$. The only problem is that at this moment, the ``solution", say $\tilde{\alpha}_{N+1}$, provided by Lemma \ref{lem-induction} does not belong to $F^{p}A^{p+q}(X)$. Instead, we only know $\tilde{\alpha}_{N+1}\in F^{p-N-1}A^{p+q}(X)$. But now we can use our assumptions (e.g. $\bigoplus_{ r\geq i\geq 1 }d_r^{p-i,q+i}(X)=0$ in Theorem \ref{mainthm-0}) to replace $\tilde{\alpha}_{N+1}$ by a desired solution $\alpha_{N+1}\in F^{p}A^{p+q}(X)$.

  \item Proof of convergence and canonical solutions: in order to show the solutions $\alpha_k$ of \eqref{eq-extsion-d-expanded-0} constitute a convergent power series $\sum_{k\geq0}\alpha_k$, we must choose the solutions $\alpha_k$ carefully in each step. For this, we developed a Hodge theory which works for forms in $F^pA^{p+q}(X)$, see Appendix \ref{appendix-A}.
\end{itemize}
\subsection{Conventions and notations}
To fix our notations and avoid redundance, we make the following conventions which will be used throughout this paper:
\begin{enumerate}
  \item Unless otherwise stated, all complex manifolds are assumed to be compact and of dimension\footnote{For compact complex manifolds of dimension $\leq2$, $E_1=E_\infty$ always holds, see \cite[Thm.\,2.8]{BHPV04}.} $\geq3$. If $X$ is a complex manifold, we use $A^{p,q}(X)$ to denote the space of smooth (complex valued) differential forms of $(p,q)$-type on $X$ and $A^k(X)$ is the space of smooth (complex valued) differential forms of degree $k$ on $X$. We use $\Pi^{p,q}:A^{p+q}(X)\to A^{p,q}(X)$ to denote the natural projection operator. Furthermore, we make the convention that $A^{p,q}(X)=0$ if $p,q<0$ or $p,q>\dim X$.
  \item For a differential form $\bullet\in A^*(X):=\oplus_{k\geq0} A^k(X)$ where $A^k(X)=\oplus_{p+q=k}A^{p,q}(X)$, we will denote by $\bullet^{p,q}$ its $(p,q)$-component. Sometimes we write a power series $\alpha(t)$ as $\alpha$ for simplicity.
  \item We always write $\pi: (\mathcal{X}, X)\to (\Delta,0)$ for a complex analytic family over a small disc $\Delta\subset \C$ where $X=X_0:=\pi^{-1}(0)$ is the central fiber and $X_t:=\pi^{-1}(t)$ is the fiber corresponding to $t\in\Delta$. We assume the central fiber of any given complex analytic family has been equipped with a fixed Hermitian metric. This is necessary because we need a Hermitian metric to construct canonical and convergent deformations. We use $\pb^*, G_{\pb}, \mathcal{H}$ to denote the formal adjoint of $\pb$, the $\pb$-Green's operator on $X$ and the projection operator onto $\pb$-Harmonic spaces $\mathcal{H}^{p,q}(X)$, respectively.
\end{enumerate}
\vskip 1\baselineskip \textbf{Acknowledgements.} We would like to thank Prof. Kefeng Liu for his constant encouragement and many useful comments. Many thanks to Prof. Sheng Rao who helped us greatly improved the expositions of this paper. We would also like to thank the anonymous referees for their careful reading and helpful suggestions.
\section{Fr\"olicher spectral sequence and $r$-filtered $(p,q)$-forms}
In this section, we will recall some basic facts on the Fr\"olicher spectral sequence, introduce the definition of deformations of $r$-filtered $(p,q)$ forms, and give a characterization of the condition $d_r^{p,q}(X_t)=0$ via canonical Dolbeault deformations.
\subsection{The Fr\"olicher spectral sequence}
Let $X$ be a complex manifold of dimension $n$. The Fr\"olicher spectral sequence of $X$ is by definition the spectral sequence associated to the filtered complex
\[
(F^\bullet A^{\bullet}(X),d),~\text{with}~F^pA^{k}(X)=\oplus_{\lambda\geq p} A^{\lambda,k-\lambda}(X).
\]
In fact, there are complexes $(E_r^{p,q}(X),d_r(X))$ with differentials
\[
d_r(X):E_r^{p,q}(X)\longrightarrow E_r^{p+r,q-r+1}(X),~r=0,1,2,\cdots
\]
such that
\begin{itemize}
  \item[1.] $E_0^{p,q}(X)=\frac{F^pA^{p+q}(X)}{F^{p+1}A^{p+q}(X)}\cong A^{p,q}(X)$ and $E_1^{p,q}(X)\cong H^{p,q}_{\pb}(X)$;
  \item[2.] for sufficiently large $r$, $E_r^{p,q}(X)\cong E_\infty^{p,q}(X)\cong \frac{F^pH^{p+q}(A^{\bullet}(X),d)}{F^{p+1}H^{p+q}(A^{\bullet}(X),d)}$,
\end{itemize}
where
\begin{align*}
E_r^{p,q}(X):=&Z_r^{p,q}(X)/\left(Z_{r-1}^{p+1,q-1}(X)+ B_{r-1}^{p,q}(X)\right)\\
F^pH^{p+q}(A^{\bullet}(X),d):=&\im \left(H^{p+q}(F^pA^{\bullet}(X),d)\longrightarrow H^{p+q}(A^{\bullet}(X),d)\right),
\end{align*}
and
\begin{align*}
Z_r^{p,q}(X):=&F^pA^{p+q}(X)\cap d^{-1}\left(F^{p+r}A^{p+q+1}(X)\right)\\
B_r^{p,q}(X):=&F^pA^{p+q}(X)\cap d\left(F^{p-r}A^{p+q-1}(X)\right)\\
Z_\infty^{p,q}(X):=&\ker d\cap F^pA^{p+q}(X)\\
B_\infty^{p,q}(X):=&\im d\cap F^pA^{p+q}(X).
\end{align*}
\subsection{The deformed Fr\"olicher spectral sequence}
Let $\pi: (\mathcal{X}, X)\to (B,0)$ be a complex analytic family such that for each $t\in B$ the complex structure on $X_t$ is represented by a Beltrami differential $\phi(t)$, then\footnote{We learned of this observation from a seminar talk given by Dingchang Wei, see also \cite{FM06,FM09,WZ20}.}
\begin{equation}\label{eq-weidingchang}
e^{i_{\phi(t)}}~: F^pA^{k}(X)\longrightarrow F^pA^{k}(X_t),
\end{equation}
is an isomorphism.

Set (see e.g. \cite{LRY15,FM06,Xia19deri} and \cite[pp.\,78]{Ma05}):
\begin{equation}\label{eq-LRY15}
d_{\phi(t)}:=e^{-i_{\phi(t)}}de^{i_{\phi(t)}} = \p + \pb_{\phi(t)}=\partial + \bar{\partial}-\mathcal{L}_{\phi(t)}^{1,0}~,
\end{equation}
where $\mathcal{L}_{\phi(t)}^{1,0}:=i_{\phi(t)}\p-\p i_{\phi(t)}$ is the Lie derivative. We then have the following commutative diagram:
\begin{equation}\label{eq-d-exp-dphi}
\xymatrix{
  F^pA^{k}(X) \ar[d]^-{d_{\phi(t)} } \ar[r]^-{e^{i_{\phi(t)}} } & F^pA^{k}(X_t) \ar[d]^-{d} \\
  F^pA^{k+1}(X) \ar[r]^-{e^{i_{\phi(t)}} } & F^pA^{k+1}(X_t), }
\end{equation}
which means exactly that
\[
e^{i_{\phi(t)}}:(F^\bullet A^{\bullet}(X),d_{\phi(t)} )\to (F^\bullet A^{\bullet}(X_t),d)
\]
is an isomorphism of filtered complexes\footnote{Note that $(F^\bullet A^{\bullet}(X),d_{\phi(t)} )$ is just the filtered complex induced by the double complex $(A^{\bullet,\bullet}(X), \p, \pb_{\phi(t)} )$.}. This induces isomorphisms
\[
E_r^{\bullet,\bullet}(X,d_{\phi(t)})\cong  E_r^{\bullet,\bullet}(X_t),
\]
and the following commutative diagram:
\begin{equation}\label{eq-drt-exp-dXt}
\xymatrix{
  E_r^{\bullet,\bullet}(X,d_{\phi(t)}) \ar[d]^-{d_r(X, d_{\phi(t)}) } \ar[r]^-{e^{i_{\phi(t)}} } & E_r^{\bullet,\bullet}(X_t) \ar[d]^-{d_r(X_t) } \\
  E_r^{\bullet+r,\bullet-r+1}(X,d_{\phi(t)}) \ar[r]^-{e^{i_{\phi(t)}} } & E_r^{\bullet+r,\bullet-r+1}(X_t), }
\end{equation}
where $E_r^{\bullet,\bullet}(X,d_{\phi(t)})$ is the $r$-th page of the spectral sequence associated to $(F^\bullet A^{\bullet}(X),d_{\phi(t)} )$ and $d_r(X, d_{\phi(t)})$ the corresponding differential operator on the $r$-th page. We will call $\left(E_r^{\bullet,\bullet}(X,d_{\phi(t)}), d_r(X, d_{\phi(t)})\right)$ the \emph{deformed Fr\"olicher spectral sequence}.
\subsection{An explicit description of the Fr\"olicher spectral sequence}
In \cite{CFUG97}, the authors gave an explicit description of the Fr\"olicher spectral sequence $(E_r^{p,q}(X),d_r(X))$ which can be described as follows. First, set
\begin{equation*}
   \wt{Z}^{p,q}_1(X)=A^{p,q}(X)\cap \ker\pb,\quad \wt{B}^{p,q}_1(X)=\pb(A^{p,q-1}(X)),
\end{equation*}
for $r\geq 2$
\begin{align*}
\begin{split}
  \wt{Z}^{p,q}_r(X)=&\{\alpha^{p,q}\in A^{p,q}(X)|\pb\alpha_{p,q}=0 \text{ and there exist }\\& \alpha^{p+i,q-i}\in A^{p+i,q-i}(X)
 \text{ such that }\\
 & \p\alpha^{p+i-1,q-i+1}+\pb\alpha^{p+i,q-i}=0,\,1\leq i\leq r-1\}\\
 =&\{\alpha^{p,q}\in A^{p,q}(X)| \exists~\alpha=\alpha^{p,q}+\cdots+\alpha^{p+r-1,q-r+1}\\
   & \text{s.t.}~ d\alpha-\p\alpha^{p+r-1,q-r+1}=0         \}
 \end{split}
\end{align*}
and
\begin{align*}
\begin{split}
  \wt{B}^{p,q}_r(X)=&\{\p\beta^{p-1,q}+\pb\beta^{p,q-1}\in A^{p,q}(X)|\text{ there exist}\\
  & \beta^{p-i,q+i-1}\in A^{p-i,q+i-1}(X),\,2\leq i\leq r-1,\\
  &\text{ satisfying } \p\beta^{p-i,q+i-1}+\pb\beta^{p-i+1,q+i-2}=0,\\
  &\pb\beta^{p-r+1,q+r-2}=0 \}\\
  =&\{ \alpha^{p,q}\in A^{p,q}(X)| \exists~\beta=\beta^{p,q-1}+\beta^{p-1,q}+\cdots+\beta^{p-r+1,q+r-2}\\
  &\quad\text{s.t.}~ \alpha^{p,q}=d\beta-\p\beta^{p-r,q+r-1} \}.
\end{split}
\end{align*}
For any $r\geq 1$ there is a canonical isomorphism
\begin{equation}\label{eq-ErtoA/B}
E_r^{p,q}(X)\longrightarrow \wt{Z}_r^{p,q}/\wt{B}_r^{p,q}:[\sum_{0\leq i\leq \min\{q,n-p\} }\alpha^{p+i,q-i}]\longmapsto [\alpha^{p,q}],
\end{equation}
where $\sum_{0\leq i\leq \min\{q,n-p\}}\alpha^{p+i,q-i}\in Z_r^{p,q}$.
In particular,
\begin{equation*}
  E_1^{p,q}(X)\cong \frac{A^{p,q}(X)\cap \ker\pb }{\pb(A^{p,q-1}(X))}=H^{p,q}_{\pb }(X)
\end{equation*}
and
\begin{equation*}
  E^{p,q}_2(X)\cong \frac{\{\alpha^{p,q}\in A^{p,q}(X)|0=\pb\alpha^{p,q}=\p\alpha^{p,q}+\pb\alpha^{p+1,q-1}\}}{\{\p\beta^{p-1,q}+\pb\beta^{p,q-1}|0=\pb\beta^{p-1,q}\}}.
\end{equation*}
Furthermore, the following diagram is commutative:
\begin{equation}\label{eq-drXtodr}
\xymatrix{
  E_r^{p,q}(X) \ar[d] \ar[r]^-{d_r^{p,q} } & E_r^{p+r,q-r+1}(X): \ar[d][\alpha ]\ar@{|->}[r]&
   [d\alpha ]  \\
  \wt{Z}_r^{p,q}/\wt{B}_r^{p,q} \ar[r]^-{ } & \wt{Z}_r^{p+r,q-r+1}/\wt{B}_r^{p+r,q-r+1}: [\alpha^{p,q}]\ar@{|->}[r]& [\p\alpha^{p+r-1,q-r+1}], }
\end{equation}
where the vertical maps are given by \eqref{eq-ErtoA/B} and $\alpha:=\sum_{0\leq i\leq \min\{q,n-p\}}\alpha^{p+i,q-i}\in Z_r^{p,q}$. Forms in $\wt{Z}^{p,q}_r(X)$ and $\wt{B}^{p,q}_r(X)$ are called $E_r$-\textit{closed/exact} $(p,q)$-forms by Popovici-Stelzig-Ugarte \cite{PSU21}.

\subsection{Deformation of closed $r$-filtered $(p,q)$-forms}
Let $\pi: (\mathcal{X}, X)\to (B,0)$ be a small deformation of a compact complex manifold $X$ such that for each $t\in B$ the complex structure on $X_t$ is represented by a Beltrami differential $\phi(t)$. In our situation, we may always assume $B$ is a small polydisc and $0\in B$ is the origin.

First, we briefly recall the definition of (canonical) Dolbeault deformations of $(p,q)$-forms (c.f. \cite{Xia19dDol}).
\begin{definition}
Given $\sigma_0\in \ker\pb\cap A^{p,q}(X)$ and $T\subseteq B$, which is an analytic subset of $B$ containing $0$, a \emph{Dolbeault deformation of $\sigma_0$ on $T$ (w.r.t. $\pi$)} is a family of $(p,q)$-forms $\sigma (t)$ such that
\begin{itemize}
  \item[1.] $\sigma (t)$ is holomorphic in $t$ and $\sigma (0) = \sigma_0$;
  \item[2.] $\pb_{\phi(t)}\sigma(t) = \pb\sigma (t) - \mathcal{L}_{\phi(t)}^{1,0}\sigma (t)=0,~\text{for~any}~t\in T$.
\end{itemize}
The Dolbeault deformation $\sigma (t)=\sum_k\sigma_k$ of $\sigma_0$ is said to be \emph{canonical} if $\sigma(t) = \sigma_0 + \bar{\partial}^*G \mathcal{L}^{1,0}_{\phi(t)}\sigma(t)$ or equivalently, $\sigma_k=\sum_{1\leq j\leq k}\pb^*G_{\pb}\mathcal{L}_{\phi_j}^{1,0}\sigma_{k-j}$ for any $k\geq 1$.
We say the deformation of $\sigma_0$ is \textit{canonically unobstructed} if there is a canonical deformation of $\sigma_0$ on $B$. We say the deformations of $(p,q)$-forms are \textit{canonically unobstructed} if for any class in the Dolbeault cohomology $H_\pb^{p,q}(X)$, it has a representative which is canonically unobstructed.
\end{definition}

\begin{definition}\label{def-filtered-forms}
For any $r\geq1$, forms in the vector space $\bigoplus_{i=0}^{r-1}A^{p+i,q-i}(X)$ will be called \textit{$r$-filtered $(p,q)$-forms on $X$}. In particular, forms in $F^pA^{p+q}(X)=\bigoplus_{i=0}^{\infty}A^{p+i,q-i}(X)$ will be simply called \textit{filtered $(p,q)$-forms on $X$}.
\end{definition}
Now, as a natural extension of the Dolbeault deformations of $(p,q)$-forms, we may consider deformations of $r$-filtered $(p,q)$-forms (this was motivated by the proof of Theorem \ref{mainthm} and Theorem \ref{mainthm-2nd}).
\begin{definition}\label{def-deformation-Dol-Zrpq}
Given
\[
\sigma_0=\sigma^{p,q}_0+\cdots+\sigma^{p+r-1,q-r+1}_0\in \left(\bigoplus_{i=0}^{r-1}A^{p+i,q-i}(X)\right)\cap\ker (d-\p\Pi^{p+r-1,q-r+1})
\]
(such forms may be called\footnote{The $(p,q)$-th component $\sigma^{p,q}_0$ of a closed $r$-filtered $(p,q)$-form $\sigma_0$ is by definition a $E_r$-closed $(p,q)$-form in the sense of Popovici-Stelzig-Ugarte, see \cite{PSU21}.} \textit{closed $r$-filtered $(p,q)$-forms on $X$}) and $T\subseteq B$, which is an analytic subset of $B$ containing $0$, a \textit{deformation} of $\sigma_0$ (w.r.t. $\pi$ ) on $T$ is a family of forms
\[
\sigma (t)=\sigma^{p,q}(t)+\cdots+\sigma^{p+r-1,q-r+1}(t)
\]
such that
\begin{itemize}
  \item[(1)] $\sigma (t)$ is holomorphic in $t$ and $\sigma(0) = \sigma_0$;
  \item[(2)] for any $t\in T$, we have
  \begin{equation}\label{eq-ext-eq}
  d_{\phi(t)}\sigma (t)-\p\Pi^{p+r-1,q-r+1}\sigma (t)=0,
  \end{equation}
  or equivalently, for $r\geq 2$
  \begin{equation}\label{eq-ext-eq-expansion}
      \left\{
      \begin{array}{ll}
       0&=\pb_{\phi(t)}\sigma^{p,q}(t),\\
       0&=\p\sigma^{p,q}(t)+\pb_{\phi(t)}\sigma^{p+1,q-1}(t),\\
        &\cdots\\
       0&=\p\sigma^{p+r-2,q-r+2}(t)+\pb_{\phi(t)}\sigma^{p+r-1,q-r+1}(t).
      \end{array} \right.
\end{equation}
and for $r=1$
\begin{equation}
0=\pb_{\phi(t)}\sigma^{p,q}(t)=\pb\sigma^{p,q}(t)-\mathcal{L}_{\phi(t)}^{1,0}\sigma^{p,q}(t).
\end{equation}
\end{itemize}
\end{definition}
In particular, deformations of $1$-filtered $(p,q)$-forms coincides with Dolbeault deformations of $(p,q)$-forms studied in \cite{Xia19dDol}.

\begin{definition}\label{def-cano-unobst}
Let $y\in \wt{Z}^{p,q}_r(X)$, then by definition there is a $\sigma_0=\sigma^{p,q}_0+\cdots+\sigma^{p+r-1,q-r+1}_0$ such that
\[
d\sigma_0-\p\Pi^{p+r-1,q-r+1}\sigma_0=0\quad\text{and}\quad\sigma^{p,q}_0=y~.
\]
By a \emph{canonical deformation of $\sigma_0$} w.r.t. $\pi$, we mean a family of forms given by
\[
\sigma (t)=\sigma^{p,q}(t)+\cdots+\sigma^{p+r-1,q-r+1}(t),\quad\text{with}\quad \sigma^{p+i,q-i}(t)=\sum_{j=0}^\infty\sigma^{p+i,q-i}_j,~0\leq i\leq r-1,
\]
where each $\sigma^{p+i,q-i}(t)$ is a convergent power series and $\sigma^{p,q}(t)$ is the canonical Dolbeault deformation of $\sigma_0^{p,q}$, i.e.
\[
\sigma^{p,q}_k=\sum_{1\leq j\leq k}\pb^*G_{\pb}\mathcal{L}_{\phi_j}^{1,0}\sigma_{k-j}^{p,q},\quad \forall k\geq1,
\]
where by our convention, each $\sigma^{p,q}_k$ is a $k$-th degree polynomial in $t$ with coefficients in $A^{p,q}(X)$. We say \textit{the deformation of $r$-filtered $(p,q)$-forms on $X$ are canonically unobstructed} if for any
\[
\sigma_0\in\left(\bigoplus_{i=0}^{r-1}A^{p+i,q-i}(X)\right)\cap\ker (d-\p\Pi^{p+r-1,q-r+1}),
\]
and any arbitrary small deformation $\pi: (\mathcal{X}, X)\to (B,0)$ of $X$ with smooth $(B,0)$, there is a canonical deformation of $\sigma_0$ on $B$.
\end{definition}
\subsection{Harmonic representatives for elements in $E_r^{\bullet,\bullet}(X,d_{\phi(t)})$}
If we set
\begin{equation*}
   \wt{Z}^{p,q}_1(X, d_{\phi(t)})=A^{p,q}(X)\cap \ker\pb_{\phi(t)},\quad \wt{B}^{p,q}_1(X)=\pb_{\phi(t)}(A^{p,q-1}(X)),
\end{equation*}
for $r\geq 2$
\begin{align*}
\begin{split}
  \wt{Z}^{p,q}_r(X, d_{\phi(t)})=&\{\alpha^{p,q}\in A^{p,q}(X)|\pb_{\phi(t)}\alpha^{p,q}=0 \text{ and there exist }\\& \alpha^{p+i,q-i}\in A^{p+i,q-i}(X)
 \text{ such that }\\
 & \p\alpha^{p+i-1,q-i+1}+\pb_{\phi(t)}\alpha^{p+i,q-i}=0,\,1\leq i\leq r-1\},
 \end{split}
\end{align*}
and
\begin{align*}
\begin{split}
  \wt{B}^{p,q}_r(X, d_{\phi(t)})=&\{\p\beta^{p-1,q}+\pb_{\phi(t)}\beta^{p,q-1}\in A^{p,q}(X)|\text{ there exist}\\
  & \beta^{p-i,q+i-1}\in A^{p-i,q+i-1}(X),\,2\leq i\leq r-1,\\
  &\text{ satisfying } \p\beta^{p-i,q+i-1}+\pb_{\phi(t)}\beta^{p-i+1,q+i-2}=0,\\
  &\pb_{\phi(t)}\beta^{p-r+1,q+r-2}=0\}.
 \end{split}
\end{align*}
\begin{proposition}\label{prop-pb*}
We have the following orthogonal direct sum decomposition:
\begin{equation}\label{eq-Zrpq-*pb}
\wt{Z}^{p,q}_r(X)=\Big(\ker\pb^*\cap\wt{Z}^{p,q}_r(X)\Big)\oplus \im\pb,
\end{equation}
\begin{equation}\label{eq-Brpq-*pb}
\wt{B}^{p,q}_r(X)=\Big(\ker\pb^*\cap\wt{B}^{p,q}_r(X)\Big)\oplus \im\pb.
\end{equation}
The following natural homomorphisms induced by inclusion are isomorphisms:
\begin{equation}\label{eq-iso-*pb-Z/B}
\frac{\ker\pb^*\cap\wt{Z}^{p,q}_r(X, d_{\phi(t)})}{\ker\pb^*\cap\wt{B}^{p,q}_s(X, d_{\phi(t)})}
\longrightarrow\frac{\wt{Z}^{p,q}_r(X, d_{\phi(t)})}{\wt{B}^{p,q}_s(X, d_{\phi(t)})},\quad r,s\geq 1,
\end{equation}
and
\begin{equation}\label{eq-iso-*pb-Z/Z}
\frac{\ker\pb^*\cap\wt{Z}^{p,q}_r(X, d_{\phi(t)})}{\ker\pb^*\cap\wt{Z}^{p,q}_s(X, d_{\phi(t)})}
\longrightarrow\frac{\wt{Z}^{p,q}_r(X, d_{\phi(t)})}{\wt{Z}^{p,q}_s(X, d_{\phi(t)})},\quad r,s\geq 1.
\end{equation}
\end{proposition}
\begin{proof}First note that: for any $\alpha^{p,q}\in \ker\pb$ we have
\[
\alpha^{p,q}\in \wt{Z}^{p,q}_r(X)\Longleftrightarrow \mathcal{H}\alpha^{p,q}\in \wt{Z}^{p,q}_r(X),
\]
because $\alpha^{p,q}-\mathcal{H}\alpha^{p,q}\in\im\pb\cap A^{p,q}(X)\subseteq\wt{Z}^{p,q}_r(X)$. Then \eqref{eq-Zrpq-*pb} follows from the Hodge decomposition
\[
A^{p,q}(X)=(\ker\pb^*\cap A^{p,q}(X))\oplus(\im\pb\cap A^{p,q}(X)).
\]
\eqref{eq-Brpq-*pb} follows immediately from \eqref{eq-Zrpq-*pb} and the fact that $\im\pb=\wt{B}^{p,q}_1(X)\subseteq \wt{B}^{p,q}_r(X)$. By ~\cite[Prop.\,4.5]{Xia19dDol} we know that \eqref{eq-iso-*pb-Z/B} holds for $r=s=1$:
\[
\frac{\ker\pb^*\cap\wt{Z}^{p,q}_1(X, d_{\phi(t)})}{\ker\pb^*\cap\wt{B}^{p,q}_1(X, d_{\phi(t)})}
\longrightarrow\frac{\wt{Z}^{p,q}_1(X, d_{\phi(t)})}{\wt{B}^{p,q}_1(X, d_{\phi(t)})}=\frac{\ker\pb_{\phi(t)}}{\im\pb_{\phi(t)}}.
\]
Combining this with the fact that
\begin{align*}
\wt{B}^{p,q}_1(X, d_{\phi(t)})&\subseteq\wt{B}^{p,q}_2(X, d_{\phi(t)})\subseteq\cdots\subseteq\wt{B}^{p,q}_\infty(X, d_{\phi(t)})\\
&\subseteq\wt{Z}^{p,q}_\infty(X, d_{\phi(t)})\subseteq\cdots\subseteq\wt{Z}^{p,q}_2(X, d_{\phi(t)})\subseteq\wt{Z}^{p,q}_1(X, d_{\phi(t)}),
\end{align*}
we get for any $r,s\geq 1$,
\[
\frac{\ker\pb^*\cap\wt{Z}^{p,q}_r(X, d_{\phi(t)})}{\ker\pb^*\cap\wt{B}^{p,q}_1(X, d_{\phi(t)})}
\cong\frac{\wt{Z}^{p,q}_r(X, d_{\phi(t)})}{\wt{B}^{p,q}_1(X, d_{\phi(t)})},~
\frac{\ker\pb^*\cap\wt{B}^{p,q}_s(X, d_{\phi(t)})}{\ker\pb^*\cap\wt{B}^{p,q}_1(X, d_{\phi(t)})}
\cong\frac{\wt{B}^{p,q}_s(X, d_{\phi(t)})}{\wt{B}^{p,q}_1(X, d_{\phi(t)})},
\]
which implies \eqref{eq-iso-*pb-Z/B} and \eqref{eq-iso-*pb-Z/Z}.
\end{proof}

\subsection{Characterization of $d_r^{p,q}(X_t)=0$ via canonical deformations}\label{sebsec-characterization-drpq=0}
For any $r=1,2,\cdots,\infty$, we set
\begin{equation}\label{eq-Vrtpq}
\begin{array}{ll}
V_{r,t}^{p,q}:=
&\Big\{ \sigma_0^{p,q}\in \mathcal{H}^{p,q}(X) \mid
\sigma^{p,q}(t)\in \wt{Z}^{p,q}_r(X, d_{\phi(t)})\\
&~\text{where}~\sigma^{p,q}(t)~\text{is the canonical Dolbeault deformation of}~\sigma_0^{p,q} \Big\}.
\end{array}
\end{equation}
Clearly, we have
\[
\mathcal{H}^{p,q}(X)\supseteq V_{1,t}^{p,q}\supseteq V_{2,t}^{p,q}\supseteq V_{3,t}^{p,q}\supseteq\cdots\supseteq V_{\infty,t}^{p,q},\quad\text{for~any}~t\in\Delta,
\]
and
\[
V_{r,0}^{p,q}=\ker\pb^*\cap\wt{Z}^{p,q}_r(X).
\]
In \cite[Prop.\,4.1]{Xia19dDol}, we have shown the following mapping $f_{1,t}^{p,q}$ is an isomorphism,
\begin{equation*}
	\begin{array}{ll}
		f_{1,t}^{p,q}: &V_{1,t}^{p,q} \longrightarrow \ker\pb^*\cap\wt{Z}^{p,q}_1(X, d_{\phi(t)})=\ker\pb^*\cap\ker\pb_{\phi(t)},\\
		&\sigma_0^{p,q}\longmapsto \sigma^{p,q}(t),
	\end{array}
\end{equation*}
where $\sigma^{p,q}(t)$ is the canonical Dolbeault deformation of $\sigma_0^{p,q}$. More generally, let us consider for any $r\geq1$ and $t\in\Delta$,
\begin{equation}\label{eq-frtpq}
\begin{array}{ll}
f_{r,t}^{p,q}: &V_{r,t}^{p,q} \longrightarrow \ker\pb^*\cap\wt{Z}^{p,q}_r(X, d_{\phi(t)}),\\
&\sigma_0^{p,q}\longmapsto \sigma^{p,q}(t),
 \end{array}
\end{equation}
where $\sigma^{p,q}(t)$ is the canonical Dolbeault deformation of $\sigma_0^{p,q}$. This homomorphism $f_{r,t}^{p,q}$ is also an isomorphism because it is just the restriction of $f_{1,t}^{p,q}$.
\begin{proposition}\label{prop-Vrtpq-drt}
Let $\pi: (\mathcal{X}, X)\to (\Delta,0)$ be a complex analytic family over a small disc $\Delta\subset \C$. Then
\[
d_r^{p,q}(X_t)=0 \Longleftrightarrow  V_{r,t}^{p,q}=V_{r+1,t}^{p,q}
\]
for any $t\in\Delta$.
\end{proposition}
\begin{proof}
In view of the commutative diagram \eqref{eq-drt-exp-dXt}, we only need to show
\[
d_r^{p,q}(X, d_{\phi(t)})=0\Longleftrightarrow V_{r,t}^{p,q}=V_{r+1,t}^{p,q}, \quad \text{for~any}~t\in\Delta.
\]
Indeed, it follows from the deformed version of Proposition \ref{prop-kerim-dr} and \eqref{eq-iso-*pb-Z/Z} that for any $(p,q), r\geq 1$ and $t\in\Delta$,
\begin{align*}
d_r^{p,q}(X, d_{\phi(t)})=0&\Longleftrightarrow\wt{Z}^{p,q}_r(X, d_{\phi(t)})= \wt{Z}^{p,q}_{r+1}(X, d_{\phi(t)})\\
&\Longleftrightarrow \ker\pb^*\cap\wt{Z}^{p,q}_r(X, d_{\phi(t)})=\ker\pb^*\cap\wt{Z}^{p,q}_{r+1}(X, d_{\phi(t)})\\
&\Longleftrightarrow V_{r,t}^{p,q}=V_{r+1,t}^{p,q},
\end{align*}
where the last line holds because $V_{r,t}^{p,q}\cong\ker\pb^*\cap\wt{Z}^{p-i,q+i}_r(X, d_{\phi(t)})$.
\end{proof}

\begin{proposition}\label{prop-upp-semicontinuous-V2tpq}
If the Dolbeault deformation of $(p+1,q)$, $(p+1,q-1)$-forms are canonically unobstructed, then $\dim V_{2,t}^{p,q}$ is upper semi-continuous.
\end{proposition}
\begin{proof}
If $\sigma^{p,q}(t)$ is the canonical Dolbeault deformation of a given $\pb$-closed form $\sigma_0^{p,q}$, then by the definition of $\wt{Z}^{p,q}_2(X, d_{\phi(t)})$ and \cite[Prop.\,4.2]{Xia19dDol}, we have
\[
\sigma^{p,q}(t)\in \wt{Z}^{p,q}_2(X, d_{\phi(t)}) \Longleftrightarrow
\mathcal{H}\mathcal{L}_{\phi(t)}^{1,0}\sigma^{p,q}(t)=0~\text{and}~\mathcal{H}_{\phi(t)}\p\sigma^{p,q}(t)=0.
\]
Therefore,
\begin{equation*}
\begin{array}{ll}
V_{2,t}^{p,q}=
&\Big\{ \sigma_0^{p,q}\in \mathcal{H}^{p,q}(X) \mid
\mathcal{H}\mathcal{L}_{\phi(t)}^{1,0}\sigma^{p,q}(t)=0~\text{and}~\mathcal{H}_{\phi(t)}\p\sigma^{p,q}(t)=0\\
&~\text{where}~\sigma^{p,q}(t)~\text{is the canonical Dolbeault deformation of}~\sigma_0^{p,q} \Big\}.
\end{array}
\end{equation*}
Now if the Dolbeault deformation of $(p+1,q)$, $(p+1,q-1)$-forms are canonically unobstructed, then by \cite[Thm.\,5.10]{Xia19dDol} the Hodge number $h^{p+1,q}$ is deformation invariant. But $h^{p+1,q}(X_t)=\dim\ker \square_{\phi(t)}\cap A^{p+1,q}(X)$, where $\square_{\phi(t)}=\pb_{\phi(t)}\pb_{\phi(t)}^*+\pb_{\phi(t)}^*\pb_{\phi(t)}$. So according to \cite[Thm.\,7.4]{Kod86}, $\mathcal{H}_{\phi(t)}:A^{p+1,q}(X)\to \ker \square_{\phi(t)}\cap A^{p+1,q}(X)$ is $C^{\infty}$ differentiable in $t$. Hence, $\dim V_{2,t}^{p,q}$ is upper semi-continuous. Indeed, let $\{e_1,\cdots, e_m\}$ and $\{e_1(t),\cdots, e_l(t)\}$ be a basis of $A^{p,q+1}(X)\cap\ker\square$ and $A^{p+1,q}(X)\cap\ker\square_{\phi(t)}$, respectively. Then
\begin{equation*}
\begin{array}{ll}
V_{2,t}^{p,q}=
&\Big\{ \sigma_0^{p,q}\in \mathcal{H}^{p,q}(X) \mid
\langle \mathcal{L}_{\phi(t)}^{1,0}\sigma^{p,q}(t), e_i\rangle_{L^2}=0,~1\leq i\leq m~\text{and}\\
&\quad\quad\quad\quad\quad\quad\quad\quad \left\langle\p\sigma^{p,q}(t), e_j(t)\right\rangle_{L^2}=0,~1\leq j\leq l,\\
&~\text{where}~\sigma^{p,q}(t)~\text{is the canonical Dolbeault deformation of}~\sigma_0^{p,q} \Big\},
\end{array}
\end{equation*}
from which we see that $\{t\in\Delta\mid \dim V_{2,t}^{p,q}\geq k\}$ is a closed subset of $\Delta$ for each $k$.
\end{proof}

\section{Unobstructed deformations of $(p,q)$-forms and stability of $d_1=0$} \label{sec-}
Let $\pi: (\mathcal{X}, X)\to (\Delta,0)$ be a complex analytic family over a small disc $\Delta\subset \C$ such that for each $t\in \Delta$ the complex structure on $X_t$ is represented by a Beltrami differential $\phi(t)$. In this section, we will prove the Dolbeault deformations of $(p,q)$-forms are unobstructed if
  \[
  \bigoplus_{r\geq 1}d_r^{p,q}(X)=0\quad\text{and}\quad
  \bigoplus_{ r\geq i\geq 1 }d_r^{p-i,q+i}(X)=0.
  \]
\subsection{The obstruction equation}
We want to show the unobstructedness of Dolbeault deformations of $(p,q)$-forms by using
the power series method. This means if we write $\phi(t)=\sum_i\phi_i$ where each $\phi_i$ is a polynomial in $t$ of degree $i$ with coefficients in $A^{0,1}(X, T^{1,0}_X)$, the objective is to find a power series solution $\sigma(t)=\sum_i\sigma_i$ of the \emph{extension equation}
\begin{equation}\label{obs1}
\pb\sigma (t)=\mathcal{L}_{\phi(t)}^{1,0}\sigma (t).
\end{equation}

\subsection{The first main theorem}
In the following, we will use the convenient notation that $\alpha^{p,q}=0$ for $p$ or $q$ not in $\{0,\cdots,n\}$.
For any $\alpha^{p,q}_0\in A^{p,q}(X)\cap \ker\pb$, our strategy is to construct a family of forms
\begin{equation*}
  \alpha(t)=\alpha^{p,q}(t)+\cdots+\alpha^{n,p+q-n}(t)\in F^pA^{p+q}(X)
\end{equation*}
holomorphic on $t$ and
such that
\begin{equation}\label{obs2}
  d \left(e^{i_{\phi(t)}}\alpha(t)\right)=0,\quad \alpha^{p,q}(0)=\alpha_0^{p,q}.
  \end{equation}
In fact, by using the identity $e^{-i_{\phi(t)}}d e^{i_{\phi(t)}}=\pb_{\phi(t)}+\p$, \eqref{obs2} is equivalent to
\begin{equation*}
  (\pb_{\phi(t)}+\p)\alpha(t)=0,
\end{equation*}
or
\begin{align*}
\begin{split}
   \begin{cases}
 \pb_{\phi(t)}\alpha^{p,q}(t)	& =0\\
 \pb_{\phi(t)} \alpha^{p+1,q-1}(t)	+\p\alpha^{p,q}(t)&=0\\
 \cdots&\\
 \pb_{\phi(t)} \alpha^{n,p+q-n}(t)+\p\alpha^{n-1,p+q-n+1}(t)&=0,
 \end{cases}
 \end{split}
\end{align*}
where the first equation is exactly \eqref{obs1}. Assume that $\alpha(t)$ can be expanded as
\begin{equation*}
  \alpha(t)=\sum_{k=0}^\infty\alpha_k,
\end{equation*}
where $\alpha_k$ is the $k$-th order homogeneous part in the expansion of $\alpha(t)$.

\begin{lemma}\label{lemma1}
	$\left( d e^{i_{\phi(t)}}\alpha(t)\right)_{N_1}=0$ for any $N_1\leq N$ if and only if
\begin{align}\label{equ1}
\begin{split}
  \left(d_{\phi(t)}\alpha(t)\right)_{N_1}=0,\quad\text{for~any}~ N_1\leq N.
   \end{split}
\end{align}
\end{lemma}
\begin{proof}
If $\left( d e^{i_{\phi(t)}}\alpha(t)\right)_{N_1}=0$ for any $N_1\leq N$, then
\begin{align*}
\left(d_{\phi(t)}\alpha(t)\right)_{N_1}&=\left(e^{-i_{\phi(t)}}d e^{i_{\phi(t)}}\alpha(t)\right)_{N_1}\\
&=\sum_{0\leq j\leq N_1}\left(e^{-i_{\phi(t)}}\right)_j\left(d e^{i_{\phi(t)}}\alpha(t)\right)_{N_1-j}\\
&=0
\end{align*}
for any $N_1\leq N$.
	
Conversely, if \eqref{equ1} holds,  then
\begin{equation*}
\left(d e^{i_{\phi(t)}}\alpha(t)\right)_{N_1}=\left(e^{i_{\phi(t)}}d_{\phi(t)}\alpha(t)\right)_{N_1}
=\sum_{0\leq j\leq N_1}\left(e^{i_{\phi(t)}}\right)_j\left(d_{\phi(t)}\alpha(t)\right)_{N_1-j}=0.
\end{equation*}
\end{proof}
Notice that
\[
(d_{\phi(t)}\alpha)_k=\left((d-\mathcal{L}_{\phi(t)}^{1,0})\alpha\right)_k=d\alpha_k-\sum_{1\leq j\leq k}\mathcal{L}_{\phi_j}^{1,0}\alpha_{k-j},\quad \forall k\geq 1.
\]
Our aim is then reduced to solve the following linear system of equations:
\[
d\alpha_k=\sum_{1\leq j\leq k}\mathcal{L}_{\phi_j}^{1,0}\alpha_{k-j},~\alpha_k\in F^{p}A^{p+q}(X),\quad \forall k\geq 1.
\]

The following observation will be essential for us:
\begin{lemma}\label{lem-induction}
If $(d_{\phi(t)}\alpha)_{N_1}=0$ for any $N_1\leq N$, then
\begin{equation}\label{eq-obstruction-to-d}
\left(d(e^{i_{\phi(t)}}-1)\alpha\right)_{N_1}=-\left(\mathcal{L}_{\phi(t)}^{1,0}\alpha\right)_{N_1},\quad \forall N_1\leq N+1.
\end{equation}
\end{lemma}
\begin{proof}We have
\begin{align*}
\begin{split}
d(e^{i_{\phi(t)}}-1)\alpha
&=d e^{i_{\phi(t)}}\alpha-d\alpha\\
&=e^{i_{\phi(t)}}d_{\phi(t)}\alpha-d\alpha\\
&=(e^{i_{\phi(t)}}-1)d_{\phi(t)}\alpha+d_{\phi(t)}\alpha-d\alpha\\
&=(e^{i_{\phi(t)}}-1)d_{\phi(t)}\alpha-\mathcal{L}_{\phi(t)}^{1,0}\alpha.
\end{split}
\end{align*}
Notice that if $(d_{\phi(t)}\alpha)_{N_1}=0$ for any $N_1\leq N$, then
\[
((e^{i_{\phi(t)}}-1)d_{\phi(t)}\alpha)_{N_1}=\sum_{k>0}( e^{i_{\phi(t)}}-1)_k (d_{\phi(t)}\alpha)_{N_1-k}=0,\quad \forall N_1\leq N+1,
\]
which implies \eqref{eq-obstruction-to-d}.
\end{proof}

\begin{lemma}\label{lem-induction-p,q-fixed}Let $\pi: (\mathcal{X}, X)\to (\Delta,0)$ be a complex analytic family over a small disc $\Delta\subset \C$. Assume
  \[
  \bigoplus_{r\geq 1}d_r^{p,q}(X)=0\quad\text{and}\quad
  F^{p}A^{p+q+1}(X)\cap dA^{p+q}(X) =dF^{p}A^{p+q}(X).
  \]
Then for any given $\alpha_0\in\ker d\cap F^{p}A^{p+q}(X)$, the following system of equations can be solved inductively:
\begin{equation}\label{eq-lem-p,q-fixed}
\left\{
\begin{array}{ll}
d\alpha_k&=\sum_{1\leq j\leq k}\mathcal{L}_{\phi_j}^{1,0}\alpha_{k-j},~\alpha_k\in F^{p}A^{p+q}(X),\quad  k\geq 0,\\
\alpha_k^{p,q}&=\sum_{1\leq j\leq k}\pb^*G_{\pb}\mathcal{L}_{\phi_j}^{1,0}\alpha^{p,q}_{k-j},~ k\geq 1,
\end{array} \right.
\end{equation}
where the summation $\sum_{1\leq j\leq k}\mathcal{L}_{\phi_j}^{1,0}\alpha_{k-j}$ is regarded as $0$ when $k=0$ and $\phi(t)=\sum_j\phi_j$ is the Beltrami differential determined by $\pi$.
\end{lemma}
\begin{proof}
{\bf Step 0:} In this case, \eqref{eq-lem-p,q-fixed} is reduced to $d\alpha_0=0$ which is clearly solvable by our assumption.

{\bf Step $\leq N$:} Assume that \eqref{eq-lem-p,q-fixed} can be solved for $k\leq N$ ($N\geq0$), namely, there exists $\alpha_0,\alpha_1,\cdots,\alpha_N\in F^{p}A^{p+q}(X)$ such that
\begin{equation}\label{eq-lem-p,q-fixed-kleqN}
\left\{
\begin{array}{ll}
d\alpha_k&=\sum_{1\leq j\leq k}\mathcal{L}_{\phi_j}^{1,0}\alpha_{k-j},\quad  0\leq k\leq N,\\
\alpha_k^{p,q}&=\sum_{1\leq j\leq k}\pb^*G_{\pb}\mathcal{L}_{\phi_j}^{1,0}\alpha^{p,q}_{k-j},\quad  1\leq k\leq N,
\end{array} \right.
\end{equation}
{\bf Step $N+1$}: We need to solve $\alpha_{N+1}\in F^{p}A^{p+q}(X)$ for
\begin{equation}\label{eq-lem-p,q-fixed-k=N+1}
\left\{
\begin{array}{ll}
d\alpha_{N+1}&=\sum_{1\leq j\leq N+1}\mathcal{L}_{\phi_j}^{1,0}\alpha_{N+1-j},\\
\alpha_{N+1}^{p,q}&=\sum_{1\leq j\leq N+1}\pb^*G_{\pb}\mathcal{L}_{\phi_j}^{1,0}\alpha^{p,q}_{N+1-j}.
\end{array} \right.
\end{equation}
First, it follows from \eqref{eq-obstruction-to-d} that
\begin{equation}\label{eq-obstruction-to-d-N+1-lem}
d\left( (e^{i_{\phi(t)}}-1)\alpha \right)_{N+1}=-\left( \mathcal{L}_{\phi(t)}^{1,0}\alpha\right)_{N+1}.
\end{equation}
Therefore, by using our assumption again, we get
\[
-\left( \mathcal{L}_{\phi(t)}^{1,0}\alpha\right)_{N+1}\in dA^{p+q}(X)\cap F^{p}A^{p+q+1}(X)=dF^{p}A^{p+q}(X).
\]
So there exists $\tilde{\alpha}_{N+1}\in F^{p}A^{p+q}(X)$ such that
$d\tilde{\alpha}_{N+1}=\left( \mathcal{L}_{\phi(t)}^{1,0}\alpha\right)_{N+1}$, in particular,
\[
\pb\tilde{\alpha}_{N+1}^{p,q}=\left( \mathcal{L}_{\phi(t)}^{1,0}\alpha\right)_{N+1}^{p,q+1}=\sum_{1\leq j\leq N+1}\mathcal{L}_{\phi_j}^{1,0}\alpha_{N+1-j}^{p,q}.
\]
It follows from the Hodge decomposition of $\sum_{1\leq j\leq N+1}\mathcal{L}_{\phi_j}^{1,0}\alpha_{N+1-j}^{p,q}$ that
\[
\pb(\pb^*G_{\pb}\sum_{1\leq j\leq N+1}\mathcal{L}_{\phi_j}^{1,0}\alpha_{N+1-j}^{p,q})=\sum_{1\leq j\leq N+1}\mathcal{L}_{\phi_j}^{1,0}\alpha_{N+1-j}^{p,q}
\]
which implies
\[
\tilde{\alpha}_{N+1}^{p,q}-\pb^*G_{\pb}\sum_{1\leq j\leq N+1}\mathcal{L}_{\phi_j}^{1,0}\alpha_{N+1-j}^{p,q}\in\ker\pb\cap A^{p,q}(X).
\]
By using $(2)$ of Proposition \ref{prop-dr=0}, there exists $\gamma_{N+1}\in\ker d\cap F^{p}A^{p+q}(X)$ such that
\[
\gamma^{p,q}_{N+1}=\tilde{\alpha}_{N+1}^{p,q}-\pb^*G_{\pb}\sum_{1\leq j\leq N+1}\mathcal{L}_{\phi_j}^{1,0}\alpha_{N+1-j}^{p,q}.
\]
Then $\alpha_{N+1}:=\tilde{\alpha}_{N+1}-\gamma_{N+1}$ is a solution of \eqref{eq-lem-p,q-fixed-k=N+1}. The proof is complete.
\end{proof}

\begin{thm}\label{mainthm}
Let $\pi: (\mathcal{X}, X)\to (\Delta,0)$ be a complex analytic family over a small disc $\Delta\subset \C$. Assume
  \[
  \bigoplus_{r\geq 1}d_r^{p,q}(X)=0\quad\text{and}\quad
  \bigoplus_{ r\geq i\geq 1 }d_r^{p-i,q+i}(X)=0.
  \]
Then
\begin{enumerate}
  \item The Dolbeault deformations of $(p,q)$-forms on $X$ are canonically unobstructed.
  \item For any $t\in\Delta$ we have
  \[
     \bigoplus_{r\geq 1}d_r^{p,q}(X_t)=0.
  \]
\end{enumerate}
\end{thm}
\begin{proof}
First of all, notice that the condition $\bigoplus_{ r\geq i\geq 1} d_r^{p-i,q+i}(X)=0$ is equivalent to $F^{p}A^{p+q+1}(X)\cap dA^{p+q}(X) =dF^{p}A^{p+q}(X)$ in view of Corollary \ref{coro-Fp+1capdA=dFp+1}. So our assumption is identical to that of Lemma \ref{lem-induction-p,q-fixed}.

$(1)$ For any given $\hat{\alpha}_0^{p,q}\in \ker\pb\cap A^{p,q}(X)$, by using Proposition \ref{prop-dr=0}, there exists $\alpha_0\in F^{p}A^{p+q}(X)\cap\ker d$ such that $\alpha_0^{p,q}=\hat{\alpha}_0^{p,q}$. By Lemma \ref{lem-induction-p,q-fixed}, there exists $\alpha_1,\alpha_2,\cdots\in F^{p}A^{p+q}(X)$ such that
\begin{equation*}
\left\{
\begin{array}{ll}
d\alpha_k&=\sum_{1\leq j\leq k}\mathcal{L}_{\phi_j}^{1,0}\alpha_{k-j},\quad   k\geq 1,\\
\alpha_k^{p,q}&=\sum_{1\leq j\leq k}\pb^*G_{\pb}\mathcal{L}_{\phi_j}^{1,0}\alpha^{p,q}_{k-j},\quad  k\geq 1.
\end{array} \right.
\end{equation*}
In particular, the Dolbeault deformation of $\hat{\alpha}_0^{p,q}$ is canonically unobstructed. Because $\hat{\alpha}_0^{p,q}$ is arbitrary, this shows the Dolbeault deformations of $(p,q)$-forms are canonically unobstructed.

$(2)$ We want to prove this by applying Proposition \ref{prop-Vrtpq-drt}. For this purpose, we try to find (for any given $\pb$-harmonic $\alpha_0^{p,q}$) a convergent power series solution
\[
\alpha(t)=\sum_{k\geq0}\alpha_k,\quad \text{with}~\alpha_k\in F^{p}A^{p+q}(X)~\text{for~any}~k,
\]
to the following
\begin{equation}
\left\{
\begin{array}{ll}
d_{\phi(t)}\alpha(t)&=0,\quad\text{for~any}~t\in\Delta,\\
\alpha_k^{p,q}&=\sum_{1\leq j\leq k}\pb^*G_{\pb}\mathcal{L}_{\phi_j}^{1,0}\alpha^{p,q}_{k-j},\quad k\geq 1.
\end{array} \right.
\end{equation}
To accomplish this, we first show for any given $\hat{\alpha}_0^{p,q}\in \ker\pb\cap A^{p,q}(X)$ the following system of equations has a convergent power series solution $\alpha(t)=\sum_{k\geq0}\alpha_k$:
\begin{equation}\label{eq-d=lie-expansion-p,q-fixed}
\left\{
\begin{array}{ll}
d(\alpha_k-\hat{\alpha}_k^{p,q})&=\sum_{1\leq j\leq k}\mathcal{L}_{\phi_j}^{1,0}(\alpha_{k-j}-\hat{\alpha}_{k-j}^{p,q})-\p\hat{\alpha}_k^{p,q},\quad \forall k\geq 1,\\
\alpha_k-\hat{\alpha}_k^{p,q}&\in F^{p+1}A^{p+q}(X),~\forall k\geq 0,\quad\text{and}~d\alpha_0=0,\\
\hat{\alpha}_k^{p,q}&=\sum_{1\leq j\leq k}\pb^*G_{\pb}\mathcal{L}_{\phi_j}^{1,0}\hat{\alpha}^{p,q}_{k-j},~\forall k\geq 1.
\end{array} \right.
\end{equation}

For $k=0$, it follows from Proposition \ref{prop-dr=0} that there exists $\alpha_0\in F^{p}A^{p+q}(X)\cap\ker d$ such that $\alpha_0^{p,q}=\hat{\alpha}_0^{p,q}\Rightarrow \alpha_0-\hat{\alpha}_0^{p,q}\in F^{p+1}A^{p+q}(X)$. According to Proposition \ref{prop-canonicalsolution-d-eq}, a canonical solution of \eqref{eq-d=lie-expansion-p,q-fixed} for $k=0$ is given by
\[
\alpha_0:=(1-d_{p+1}^*G_{p+1}\p)\hat{\alpha}_0^{p,q}.
\]

Assume \eqref{eq-d=lie-expansion-p,q-fixed} can be solved for $k\leq N$ with $N\geq0$, we need to show it can also be solved for $k=N+1$.

Indeed, since \eqref{eq-d=lie-expansion-p,q-fixed} can be solved for $k\leq N$, there exists $\alpha_0,\alpha_1,\cdots,\alpha_N$ such that
\begin{equation}\label{eq-d=lie-expansion-p,q-fixed-leqN}
\left\{
\begin{array}{ll}
d(\alpha_k-\hat{\alpha}_k^{p,q})&=\sum_{1\leq j\leq k}\mathcal{L}_{\phi_j}^{1,0}(\alpha_{k-j}-\hat{\alpha}_{k-j}^{p,q})-\p\hat{\alpha}_k^{p,q},\quad  1\leq k\leq N,\\
\alpha_k-\hat{\alpha}_k^{p,q}&\in F^{p+1}A^{p+q}(X),~1\leq k\leq N,\quad\text{and}~d\alpha_0=0,\\
\hat{\alpha}_k^{p,q}&=\sum_{1\leq j\leq k}\pb^*G_{\pb}\mathcal{L}_{\phi_j}^{1,0}\hat{\alpha}^{p,q}_{k-j},~ 1\leq k\leq N.
\end{array} \right.
\end{equation}
Combining this with the fact that
\[
\pb\hat{\alpha}_k^{p,q}=\sum_{1\leq j\leq k}\mathcal{L}_{\phi_j}^{1,0}\hat{\alpha}^{p,q}_{k-j},~ k\geq 1,
\]
which holds since we already showed in $(1)$ that the Dolbeault deformations of $(p,q)$-forms are canonically unobstructed, we get
\begin{equation*}
\left\{
\begin{array}{ll}
d\alpha_k&=\sum_{1\leq j\leq k}\mathcal{L}_{\phi_j}^{1,0}\alpha_{k-j},\quad  1\leq k\leq N,\\
\alpha_k^{p,q}&=\sum_{1\leq j\leq k}\pb^*G_{\pb}\mathcal{L}_{\phi_j}^{1,0}\alpha^{p,q}_{k-j},\quad  1\leq k\leq N.
\end{array} \right.
\end{equation*}
So by Lemma \ref{lem-induction-p,q-fixed} there exists $\tilde{\alpha}_{N+1}\in F^{p}A^{p+q}(X)$ such that
\[
d\tilde{\alpha}_{N+1}=\sum_{1\leq j\leq N+1}\mathcal{L}_{\phi_j}^{1,0}\alpha_{N+1-j},
\]
and (it follows from the second line of \eqref{eq-d=lie-expansion-p,q-fixed-leqN} that $\alpha_k^{p,q}=\hat{\alpha}_k^{p,q}$ for $1\leq k\leq N$)
\[
\tilde{\alpha}_{N+1}^{p,q}=\sum_{1\leq j\leq N+1}\pb^*G_{\pb}\mathcal{L}_{\phi_j}^{1,0}\alpha^{p,q}_{N+1-j}
=\sum_{1\leq j\leq N+1}\pb^*G_{\pb}\mathcal{L}_{\phi_j}^{1,0}\hat{\alpha}^{p,q}_{N+1-j}=\hat{\alpha}^{p,q}_{N+1}.
\]
It follows that
\begin{align*}
&\sum_{1\leq j\leq N+1}\mathcal{L}_{\phi_j}^{1,0}(\alpha_{N+1-j}-\hat{\alpha}_{N+1-j}^{p,q})-\p\hat{\alpha}_{N+1}^{p,q}\\
=&d(\tilde{\alpha}_{N+1}-\hat{\alpha}_{N+1}^{p,q})\in d\left(F^{p+1}A^{p+q}(X)\right).
\end{align*}
which by Proposition \ref{prop-canonicalsolution-d-eq} implies
\[
\alpha_{N+1}:=\hat{\alpha}_{N+1}^{p,q}+d_{p+1}^*G_{p+1}\left(\sum_{1\leq j\leq N+1}\mathcal{L}_{\phi_j}^{1,0}(\alpha_{N+1-j}-\hat{\alpha}_{N+1-j}^{p,q})-\p\hat{\alpha}_{N+1}^{p,q} \right)
\]
is a solution of \eqref{eq-d=lie-expansion-p,q-fixed} for $k=N+1$.

It follows from slight modifications of standard arguments (see \cite[pp.162]{MK71} and \cite[pp.2950]{Xia19dDol}) that the power series $\alpha(t)=\sum_{k\geq0}\alpha_k$ defined by
\[
\alpha_{k}:=\hat{\alpha}_{k}^{p,q}+d_{p+1}^*G_{p+1}\left(\sum_{1\leq j\leq k}\mathcal{L}_{\phi_j}^{1,0}(\alpha_{k-j}-\hat{\alpha}_{k-j}^{p,q})-\p\hat{\alpha}_{k}^{p,q} \right),\quad k\geq 1,
\]
converges and gives a smooth form in $F^{p}A^{p+q}(X)$ for any fixed small $t\in\Delta$. This power series $\alpha(t)=\sum_{k\geq0}\alpha_k$ is thus a convergent solution of \eqref{eq-d=lie-expansion-p,q-fixed}.

Therefore, for any $\hat{\alpha}_0^{p,q}\in \mathcal{H}^{p,q}(X)$ and $t\in\Delta$ we can find a solution $\alpha(t)=\sum_{k\geq0}\alpha_k\in F^{p}A^{p+q}(X)$ to the following equation
\begin{equation*}
\left\{
\begin{array}{ll}
d_{\phi(t)}\alpha(t)&=0,\\
\alpha_k^{p,q}&=\sum_{1\leq j\leq k}\pb^*G_{\pb}\mathcal{L}_{\phi_j}^{1,0}\alpha^{p,q}_{k-j},\quad k\geq 1,~\text{and}~\alpha_0^{p,q}=\hat{\alpha}_0^{p,q},
\end{array} \right.
\end{equation*}
which implies $\hat{\alpha}_0^{p,q}\in V_{\infty,t}^{p,q}$ (see Subsection \ref{sebsec-characterization-drpq=0} for the definition of $V_{r,t}^{p,q}$). Hence, $V_{\infty,t}^{p,q}=\mathcal{H}^{p,q}(X)$. But we already know in $(1)$ that the Dolbeault deformations of $(p,q)$-forms are canonically unobstructed which is equivalent to $V_{1,t}^{p,q}=\mathcal{H}^{p,q}(X)$ for any $t\in\Delta$. As a result, we get $V_{\infty,t}^{p,q}=V_{1,t}^{p,q}$ for any $t\in\Delta$. Hence $\bigoplus_{r\geq 1}d_r^{p,q}(X_t)=0$ for any $t\in\Delta$ by Proposition \ref{prop-Vrtpq-drt}.
\end{proof}
It follows immediately from Theorem \ref{mainthm} that
\begin{corollary}\label{coro-h0,q}
The deformations of $(0,q)$-forms on $X$ are canonically unobstructed if $\bigoplus_{r\geq 1}d_r^{0,q}(X)=0$. In particular, $h^{0,1}$ is deformation invariant if $\bigoplus_{r\geq 1}d_r^{0,1}(X)=0$.
\end{corollary}
As observed in \cite{HX24A}, canonical unobstructedness of $(p,q)$-forms is equivalent to the deformation invariance of $\sum_{i=0}^{q}(-1)^{q-i}h^{p,i}$, we thus have the following
\begin{corollary}\label{coro-altsum-hp,q}
The alternating sum of Hodge numbers: $\sum_{i=0}^{q}(-1)^{q-i}h^{p,i}$, is deformation invariant if $\bigoplus_{r\geq 1}d_r^{p,q}(X)=0=\bigoplus_{ r\geq i\geq 1} d_r^{p-i,q+i}(X)$.
\end{corollary}
\begin{proof}We will sketch the proof for the reader's convenience. In fact, it is proved in \cite[Thm.\,1.2]{Xia19dDol} that for any given complex analytic family $\pi: (\mathcal{X}, X)\to (\Delta,0)$, we have
\begin{equation}\label{eq-vtpq}
h^{p,q}(X)=h^{p,q}(X_t)+v^{p,q}_t+v^{p,q-1}_t,
\end{equation}
where $v^{p,q}_t:=\dim H^{p,q}_{\pb}(X)-\dim \ker\pb_{\phi(t)}\cap\ker\pb^*\cap A^{p,q}(X) \geq 0$. It follows from \eqref{eq-vtpq} that
\begin{equation}\label{eq-vtq=altsumhi}
v^{p,q}_t=\sum_{i=0}^{q}(-1)^{q-i}h^{p,i}(X)-\sum_{i=0}^{q}(-1)^{q-i}h^{p,i}(X_t).
\end{equation}
Since the Dolbeault deformations of $(p,q)$-forms are canonically unobstructed iff $v^{p,q}_t=0$ for every $t\in\Delta$ (c.f. the proof of \cite[Prop.\,4.10]{Xia19dDol}), we conclude that canonical unobstructedness of $(p,q)$-forms is equivalent to the deformation invariance of $\sum_{i=0}^{q}(-1)^{q-i}h^{p,i}$.
\end{proof}

In general, for the deformation invariance of a single Hodge number $h^{p,q}$, Combining Theorem \ref{mainthm} with \cite[Thm\,1.2]{Xia19dDol} (in the case of $E=\Omega^p$) we have the following
\begin{corollary}\label{coro-def-invariance-hpq}
Assume
  \[
  \bigoplus_{\substack{j=0,1\\ r\geq 1}}d_r^{p,q-j}(X)=0=
  \bigoplus_{\substack{j=0,1\\ r\geq i\geq 1}}d_r^{p-i,q+i-j}(X).
  \]
Then the Hodge number $h^{p,q}$ is deformation invariant.
\end{corollary}

\begin{remark}
It is proved in \cite[Thm.\,1.3]{RZ18} that if $H_{BC}^{p+1,q}(X)\to H_{\partial}^{p+1,q}(X)$ is injective, $\p_{A,\bar{\partial}}^{p-1,q+1}=0$ and $\dim H_{\bar{\partial}_t}^{p,q-1}(X_t)$ is deformation invariant then $\dim H_{\bar{\partial}_t}^{p,q}(X_t)$ is also deformation invariant. Here, $\p_{A,\bar{\partial}}^{p-1,q+1}$ is the natural homomorphism
\[
\p_{A,\pb}^{p-1,q+1}:H_{BC}^{p-1,q+1}(X)\longrightarrow H_{\pb}^{p,q+1}(X)
\]
induced by $\p$. When $\clubsuit,\diamondsuit\in\{A,BC,\pb,\pb(\ker\p)\}$, homomorphisms
\[
\p_{\clubsuit,\diamondsuit}^{\bullet,\bullet}:H_{\clubsuit}^{\bullet,\bullet}(X)\longrightarrow H_{\diamondsuit}^{\bullet+1,\bullet}(X)
\]
can be defined in a similar way, c.f. \cite{Xia19dDol}. It is easy to see that
\begin{itemize}
  \item $\p_{A,\pb}^{p-1,q+1}=0\Longleftrightarrow H_{\pb}^{p,q+1}(X)\to H_{A}^{p,q+1}(X)$ is injective $\Longleftrightarrow X\in \mathbb{S}^{p,q+1}$;
  \item $\p_{A,BC}^{p,q}=0\Longleftrightarrow H_{BC}^{p+1,q}(X)\to H_{\p}^{p+1,q}(X)$ is injective $\Longleftrightarrow X\in \mathbb{B}^{p+1,q}$;
  \item $\p_{\pb,BC}^{p,q}=0\Longleftrightarrow H_{BC}^{p,q}(X)\to H_{\pb}^{p,q}(X)$ is surjective $\Longleftrightarrow X\in \mathcal{B}^{p+1,q}$,
\end{itemize}
where $\mathbb{S}^{p,q+1}$, $\mathbb{B}^{p+1,q}$ and $\mathcal{B}^{p+1,q}$ represent some conditions related to $\p\pb$-lemma, see \cite{RZ15,RZ18} for more information. In summary, we get the following
\renewcommand\arraystretch{1.5}
\begin{table}[!htbp]
\caption{Conditions \textrm{I}+\textrm{II} $\Longrightarrow$ deformation invariance of $h^{p,q}$}
\centering
\begin{center}
\begin{tabular}{|c|c|c|}
\hline
 & Condition \textrm{I} & Condition \textrm{II}  \\
\hline
$(1)$\cite{RZ18} & $\p_{A,\pb}^{p-1,q+1}=0=\p_{A,BC}^{p,q}$ & deformation invariance of $\dim H_{\pb}^{p,q-1}(X)$  \\
\hline
$(2)$\cite{Xia19dDol} & $\p_{A,\pb}^{p-1,q+1}=0=\p_{A,\pb(\ker\p)}^{p,q}$ & $\p_{A,\pb}^{p-1,q}=0=\p_{A,\pb(\ker\p)}^{p,q-1}$  \\
\hline
$(3)$\cite{RZ22} & $\p_{A,BC}^{p-1,q+1}=0=\p_{\pb,BC}^{p,q}$ &  deformation invariance of $\dim H_{\pb}^{p,q-1}(X)$ \\
\hline
$(4)$Corollary \ref{coro-def-invariance-hpq} &
$\bigoplus_{r\geq 1}d_r^{p,q}(X)=0$ & $\bigoplus_{r\geq 1}d_r^{p,q-1}(X)=0$  \\

 & $\bigoplus_{ r\geq i\geq 1} d_r^{p-i,q+i}(X)=0$
 &  $\bigoplus_{ r\geq i\geq 1} d_r^{p-i,q-1+i}(X)=0$\\
\hline
\end{tabular}
\end{center}
\end{table}
\end{remark}
\newpage
\section{Deformation stability of $d_r=0$ for general $r\geq1$}
In this section, we will study the deformation stability of the condition $d_r^{p,q}=0$ for general $r\geq1$.

\subsection{Formally solve the deformation equation}
\begin{lemma}\label{lem-induction-p,q-fixed-sec4}Let $\pi: (\mathcal{X}, X)\to (\Delta,0)$ be a complex analytic family over a small disc $\Delta\subset \C$ and $r\in\Z^+$.
For any given $\hat{\alpha}_0^{p,q}+\cdots+\hat{\alpha}_0^{p+r-1,q-r+1}\in\ker d_{\Pi_{r}}$, let
	\[
	\hat{\alpha}(t)=\hat{\alpha}^{p,q}(t)+\cdots+\hat{\alpha}^{p+r-1,q-r+1}(t)=\sum_k\hat{\alpha}^{p,q}_k+\cdots+\sum_k\hat{\alpha}^{p+r-1,q-r+1}_k,
	\]
be a canonical deformation of $\sum_{i=0}^{r-1}\hat{\alpha}_0^{p+i,q-i}\in\ker d_{\Pi_{r}}$ on $\Delta$. Assume
\begin{equation}\label{eq-asumption-lem2-dr=0}
F^{p+r}A^{p+q+1}(X)\cap dA^{p+q}(X)= dF^{p+r}A^{p+q}(X).
\end{equation}
Then the following system of equations can be solved inductively:
	\begin{equation}\label{eq-Z_r+1-p,q-r-canonical}
		\left\{
		\begin{array}{ll}
			d\alpha_k&=\sum_{1\leq j\leq k}\mathcal{L}_{\phi_j}^{1,0}\alpha_{k-j},~\alpha_k\in F^{p}A^{p+q}(X),\quad k\geq 0, \\
			\alpha_k^{p+i,q-i}&=\hat{\alpha}_k^{p+i,q-i},~0\leq i\leq r-1,~ k\geq 0.
		\end{array} \right.
	\end{equation}
where $\phi(t)=\sum_j\phi_j$ is the Beltrami differential determined by $\pi$.	
\end{lemma}
\begin{proof}It is enough to solve inductively the following equation of formal power series
\begin{equation}\label{eq-formally-ab}
d_{\phi(t)}\Big( \hat{\alpha}(t)+\beta(t)\Big)=0.
\end{equation}
for $\beta(t)=\sum_k\beta_k\in F^{p+r}A^{p+q}(X)$. Since $\hat{\alpha}(t)$ is a deformation of the $r$-filtered form $\hat{\alpha}_0$, by Definition \ref{def-deformation-Dol-Zrpq} we have
\begin{equation}\label{eq-alphat}
d_{\phi(t)}\hat{\alpha}(t)=d\hat{\alpha}(t)-\mathcal{L}_{\phi(t)}^{1,0}\hat{\alpha}(t)=\p\hat{\alpha}^{p+r-1,q-r+1}(t),~\text{for~any}~ t\in\Delta.
\end{equation}
This implies \eqref{eq-formally-ab} is equivalent to the following
\begin{equation}\label{eq-formally-ab-vary}
d_{\phi(t)}\beta(t)+\p\hat{\alpha}^{p+r-1,q-r+1}(t)=0.
\end{equation}
In what follows, we will show \eqref{eq-formally-ab-vary} can be solved inductively for $\beta(t)=\sum_k\beta_k\in F^{p+r}A^{p+q}(X)$. In fact, for $k=0$, we need to solve $\beta_0\in F^{p+r}A^{p+q}(X)$ for
\[
d\beta_0+\p\hat{\alpha}^{p+r-1,q-r+1}_0=0,\quad\text{or}\quad d\beta_0=-\p\hat{\alpha}^{p+r-1,q-r+1}_0
\]
which is possible because $\hat{\alpha}_0\in\ker d_{\Pi_{r}}$ implies
\[
 \p\hat{\alpha}^{p+r-1,q-r+1}_0=d\hat{\alpha}_0\in F^{p+r}A^{p+q+1}(X)\cap dA^{p+q}(X)= dF^{p+r}A^{p+q}(X),
\]
where the last equality is \eqref{eq-asumption-lem2-dr=0}.

Assume \eqref{eq-formally-ab-vary} can be solved for $k\leq N$, we need to show it can also be solved for $k=N+1$, that is, there is a $\beta_{N+1}\in F^{p+r}A^{p+q}(X)$ such that
\begin{equation}\label{eq-lem2-N+1}
d\beta_{N+1}-\left( \mathcal{L}_{\phi(t)}^{1,0}\beta\right)_{N+1}+\p\hat{\alpha}^{p+r-1,q-r+1}_{N+1}=0.
\end{equation}
The induction hypothesis implies \eqref{eq-formally-ab} can also be solved for $k\leq N$ and applying Lemma \ref{lem-induction}, we get
\begin{align*}
	d\Big( (e^{i_{\phi(t)}}-1)(\hat{\alpha}+\beta) \Big)_{N+1}=&-\left( \mathcal{L}_{\phi(t)}^{1,0}(\hat{\alpha}+\beta)\right)_{N+1}\\
=&-\left( d\hat{\alpha}-\p\hat{\alpha}^{p+r-1,q-r+1}+\mathcal{L}_{\phi(t)}^{1,0}\beta\right)_{N+1},
\end{align*}
which implies,
\[
\p\hat{\alpha}^{p+r-1,q-r+1}_{N+1}-\left(\mathcal{L}_{\phi(t)}^{1,0}\beta\right)_{N+1}=d\Big( (e^{i_{\phi(t)}}-1)(\hat{\alpha}+\beta)+\hat{\alpha}\Big)_{N+1}\in F^{p+r}A^{p+q+1}(X)\cap dA^{p+q}(X).
\]
By using \eqref{eq-asumption-lem2-dr=0} again, we have
\[
\left(\p\hat{\alpha}^{p+r-1,q-r+1}-\mathcal{L}_{\phi(t)}^{1,0}\beta\right)_{N+1}\in dF^{p+r}A^{p+q}(X),
\]
which means \eqref{eq-lem2-N+1} can be solved.
\end{proof}

\subsection{The second main theorem}
\begin{thm}\label{mainthm-2nd}
Let $\pi: (\mathcal{X}, X)\to (\Delta,0)$ be a complex analytic family over a small disc $\Delta\subset \C$. Assume for some $r\geq1$,
\begin{itemize}
  \item the deformation of $r$-filtered $(p,q)$-forms on $X$ are canonically unobstructed;
  \item there holds
\begin{equation}\label{eq-asumption-thm2-FdA=dF}
F^{p+r}A^{p+q+1}(X)\cap dA^{p+q}(X)= dF^{p+r}A^{p+q}(X).
\end{equation}
or equivalently,
\begin{equation}\label{eq-asumption-thm2-dr=0}
\bigoplus_{\lambda\geq i\geq 1}d_\lambda^{p+r-i,q-r+i}(X)=0.
\end{equation}
\end{itemize}
Let $\mu\geq r$. If $\dim V_{\mu,t}^{p,q}$ as defined in \eqref{eq-Vrtpq} is upper semi-continuous: $\dim V_{\mu,0}^{p,q}\geq\dim V_{\mu,t}^{p,q}$ for any $t\in\Delta$, then
\[
\bigoplus_{\lambda\geq \mu}d_\lambda^{p,q}(X_t)=0,\quad \text{for~any}~t\in\Delta.
\]
\end{thm}
\begin{proof}The equivalence of \eqref{eq-asumption-thm2-FdA=dF} with \eqref{eq-asumption-thm2-dr=0} follows from Corollary \ref{coro-Fp+1capdA=dFp+1} by replacing $(p,q)$ with $(p+r-1,q-r+1)$.

First, for any given $\hat{\alpha}_0^{p,q}\in \wt{Z}_{r}^{p,q}(X)$ there exists $\hat{\alpha}_0^{p+i,q-i}$ for $1\leq i\leq r-1$ such that $\sum_{i=0}^{r-1}\hat{\alpha}_0^{p+i,q-i}\in\ker d_{\Pi_r}$, i.e.
\begin{equation}\label{eq-d_Pir-closed}
d(\hat{\alpha}_0^{p,q}+\cdots+\hat{\alpha}_0^{p+r-1,q-r+1})-\p\hat{\alpha}_0^{p+r-1,q-r+1}=0.
\end{equation}
Since the deformation of $r$-filtered $(p,q)$-forms on $X$ are canonically unobstructed on $\Delta$, there exists convergent power series $\hat{\alpha}^{p+i,q-i}=\hat{\alpha}^{p+i,q-i}(t)$ for $0\leq i\leq r-1$ such that
\begin{equation}\label{eq-r-filtered-unobstruction}
	\left\{
	\begin{array}{ll}
d_{\phi(t)}(\hat{\alpha}^{p,q}+\cdots+\hat{\alpha}^{p+r-1,q-r+1})-\p\hat{\alpha}^{p+r-1,q-r+1}=0,~\text{for~any}~t\in\Delta,\\
\hat{\alpha}_k^{p,q}=\sum_{j=1}^k\pb^*G_{\pb}\mathcal{L}_{\phi_j}^{1,0}\hat{\alpha}^{p,q}_{k-j},~k\geq 1.
\end{array} \right.
\end{equation}
We will show the following system of equations has a convergent power series solution $\alpha(t)=\sum_{k\geq0}\alpha_k$ with each $\alpha_k\in F^pA^{p+q}(X)$,
\begin{equation}\label{eq-d=lie-expansion-p,q-fixed-r}
\left\{
\begin{split}
&d(\alpha_k-\sum_{i=0}^{r-1}\hat{\alpha}_k^{p+i,q-i})\\
=&\sum_{j=1}^k\mathcal{L}_{\phi_j}^{1,0}(\alpha_{k-j}-\sum_{i=0}^{r-1}\hat{\alpha}_{k-j}^{p+i,q-i})-\p\hat{\alpha}_k^{p+r-1,q-r+1},\quad \forall k\geq 1,\\
&\alpha_k-\sum_{i=0}^{r-1}\hat{\alpha}_k^{p+i,q-i}\in F^{p+r}A^{p+q}(X),~\forall k\geq 0,\quad\text{and}~d\alpha_0=0.
\end{split} \right.
\end{equation}
For $k=0$, we need to solve $\alpha_0$ for
\[
\alpha_0-\sum_{i=0}^{r-1}\hat{\alpha}_0^{p+i,q-i}\in F^{p+r}A^{p+q}(X),\quad\text{and}~d\alpha_0=0.
\]
This can be solved because $\sum_{i=0}^{r-1}\hat{\alpha}_0^{p+i,q-i}\in\ker d_{\Pi_r}$ and it follows from \eqref{eq-asumption-thm2-FdA=dF} that
\[
d\sum_{i=0}^{r-1}\hat{\alpha}_0^{p+i,q-i}=\p\hat{\alpha}_0^{p+r-1,q-r+1}=d\beta_0~\text{for~some}~\beta_0\in F^{p+r}A^{p+q}(X).
\]
So we may set $\alpha_0=\sum_{i=0}^{r-1}\hat{\alpha}_0^{p+i,q-i}-\beta_0$.

Assume \eqref{eq-d=lie-expansion-p,q-fixed-r} can be solved for $k\leq N$, we need to show it can also be solved for $k=N+1$. Indeed, from \eqref{eq-r-filtered-unobstruction} we get
\[
d(\sum_{i=0}^{r-1}\hat{\alpha}_k^{p+i,q-i})
=\sum_{j=1}^k\mathcal{L}_{\phi_j}^{1,0}(\sum_{i=0}^{r-1}\hat{\alpha}_{k-j}^{p+i,q-i})-\p\hat{\alpha}_k^{p+r-1,q-r+1},\quad \forall k\geq 1,
\]
which combining with \eqref{eq-d=lie-expansion-p,q-fixed-r} for $k\leq N$ gives
\begin{equation}\label{eq-4}
\left\{
\begin{split}
&d\alpha_k=\sum_{j=1}^k\mathcal{L}_{\phi_j}^{1,0}\alpha_{k-j},\quad \text{for~any}~1\leq k\leq N,\\
&\alpha_k-\sum_{i=0}^{r-1}\hat{\alpha}_k^{p+i,q-i}\in F^{p+r}A^{p+q}(X),~\text{for~any}~1\leq k\leq N,\quad\text{and}~d\alpha_0=0.
\end{split} \right.
\end{equation}
Notice that \eqref{eq-4} is exactly \eqref{eq-Z_r+1-p,q-r-canonical} for $k\leq N$. So we may apply Lemma \ref{lem-induction-p,q-fixed-sec4} and conclude that there is a
$\tilde{\alpha}_{N+1}\in F^{p}A^{p+q}(X)$ such that
\begin{equation*}
\left\{
\begin{split}
&d\tilde{\alpha}_{N+1}=\sum_{j=1}^{N+1}\mathcal{L}_{\phi_j}^{1,0}\alpha_{N+1-j},\\
&\tilde{\alpha}_{N+1}-\sum_{i=0}^{r-1}\hat{\alpha}_{N+1}^{p+i,q-i}\in F^{p+r}A^{p+q}(X).
\end{split} \right.
\end{equation*}
Therefore,
\begin{align*}
&\sum_{j=1}^{N+1}\mathcal{L}_{\phi_j}^{1,0}(\alpha_{N+1-j}-\sum_{i=0}^{r-1}\hat{\alpha}_{N+1-j}^{p+i,q-i})-\p\hat{\alpha}_{N+1}^{p+r-1,q-r+1}\\
=&d(\tilde{\alpha}_{N+1}-\sum_{i=0}^{r-1}\hat{\alpha}_{N+1}^{p+i,q-i})\in d(F^{p+r}A^{p+q}(X)).
\end{align*}
which by Proposition \ref{prop-canonicalsolution-d-eq} implies
\[
\alpha_{N+1}:=\sum_{i=0}^{r-1}\hat{\alpha}_{N+1}^{p+i,q-i}+d_{p+r}^*G_{p+r}
\left(\sum_{j=1}^{N+1}\mathcal{L}_{\phi_j}^{1,0}(\alpha_{N+1-j}-\sum_{i=0}^{r-1}\hat{\alpha}_{N+1-j}^{p+i,q-i})-\p\hat{\alpha}_{N+1}^{p+r-1,q-r+1} \right)
\]
is a solution of \eqref{eq-d=lie-expansion-p,q-fixed-r} for $k=N+1$.

Combining \eqref{eq-r-filtered-unobstruction} and \eqref{eq-d=lie-expansion-p,q-fixed-r} we see that for any $\hat{\alpha}_0^{p,q}\in  \wt{Z}_{r}^{p,q}(X)$ and $t\in\Delta$, the power series $\alpha(t)=\sum_{k\geq0}\alpha_k\in F^{p}A^{p+q}(X)$ defined by: for each $k>0$,
\[
\alpha_{k}:=\sum_{i=0}^{r-1}\hat{\alpha}_{k}^{p+i,q-i}+d_{p+r}^*G_{p+r}
\left(\sum_{j=1}^{k}\mathcal{L}_{\phi_j}^{1,0}(\alpha_{k-j}-\sum_{i=0}^{r-1}\hat{\alpha}_{k-j}^{p+i,q-i})-\p\hat{\alpha}_{k}^{p+r-1,q-r+1} \right),
\]
is a solution of
\begin{equation}\label{eq-canonical-alphat}
	\left\{
	\begin{array}{ll}
		d_{\phi(t)}\alpha(t)&=0,\quad~\text{with}~\alpha_0-\sum_{i=0}^{r-1}\hat{\alpha}_0^{p+i,q-i}\in F^{p+r}A^{p+q}(X),\\
		\alpha_k^{p,q}&=\sum_{1\leq j\leq k}\pb^*G_{\pb}\mathcal{L}_{\phi_j}^{1,0}\alpha^{p,q}_{k-j},\quad k\geq 1,
	\end{array} \right.
\end{equation}
where $\sum_{i=0}^{r-1}\hat{\alpha}_0^{p+i,q-i}\in\ker d_{\Pi_r}$ is a closed $r$-filtered $(p,q)$-form whose $(p,q)$-component is $\hat{\alpha}_0^{p,q}$ and $\sum_{i=0}^{r-1}\hat{\alpha}^{p+i,q-i}(t)=\sum_{i=0}^{r-1}\sum_{k\geq 0}\hat{\alpha}^{p+i,q-i}_k$ is a canonical deformation of $\sum_{i=0}^{r-1}\hat{\alpha}_0^{p+i,q-i}$.

From these we get that for any $\mu\geq r$,
\[
\ker\pb^*\cap\wt{Z}_{\mu}^{p,q}(X)\subseteq\ker\pb^*\cap\wt{Z}_{r}^{p,q}(X)\subseteq V_{\infty,t}^{p,q}(\subseteq V_{\mu,t}^{p,q}),\quad \text{for~any}~t\in\Delta.
\]
Since $V_{\mu,0}^{p,q}=\ker\pb^*\cap\wt{Z}_{\mu}^{p,q}(X)$ and $\dim V_{\mu,t}^{p,q}$ is upper semi-continuous, we get
\[
V_{\mu,t}^{p,q}=V_{\infty,t}^{p,q},\quad \text{for~any}~t\in\Delta.
\]
Then it follows from Proposition \ref{prop-Vrtpq-drt} that
\[
\bigoplus_{\lambda\geq \mu}d_\lambda^{p,q}(X_t)=0,\quad \text{for~any}~t\in\Delta.
\]
\end{proof}
\begin{remark}
We note that if the deformation of $r$-filtered $(p,q)$-forms on $X$ are canonically unobstructed, it follows immediately that
\[
(V_{r,0}^{p,q}=)\ker\pb^*\cap\wt{Z}_{r}^{p,q}(X)\subseteq V_{r,t}^{p,q},\quad \text{for~any}~t\in\Delta.
\]
\end{remark}

\begin{corollary}\label{coro-to-mainthm2}
Let $\pi: (\mathcal{X}, X)\to (\Delta,0)$ be a complex analytic family over a small disc $\Delta\subset \C$.
\begin{enumerate}
\item Assume
\begin{itemize}
  \item the deformations of $(p,q)$-forms on $X$ are canonically unobstructed or equivalently, the alternating sum of Hodge numbers $\sum_{i=0}^{q}(-1)^{q-i}h^{p,i}$ is deformation invariant;
  \item there holds
\begin{equation}
\bigoplus_{\lambda\geq i\geq 1}d_\lambda^{p+1-i,q-1+i}(X)=0.
\end{equation}
\end{itemize}
Then
\[
\bigoplus_{\lambda\geq 1}d_\lambda^{p,q}(X_t)=0,\quad \text{for~any}~t\in\Delta.
\]
\item Assume
\begin{itemize}
  \item the deformations of $(p+1,q)$-forms, $(p+1,q-1)$-forms and $2$-filtered $(p,q)$-forms on $X$ are canonically unobstructed;
  \item there holds
\begin{equation}
\bigoplus_{ \lambda\geq i\geq 1 }d_\lambda^{p+2-i,q-2+i}(X)=0.
\end{equation}
\end{itemize}
Then
\[
\bigoplus_{\lambda\geq 2}d_\lambda^{p,q}(X_t)=0,\quad \text{for~any}~t\in\Delta.
\]
\end{enumerate}
\end{corollary}
\begin{proof}
$(1)$ First note that $\dim V_{1,t}^{p,q}$ is upper semi-continuous by the arguments in the paragraph just after \cite[Definition\,4.8]{Xia19dDol} (where $V_{1,t}^{p,q}$ is denoted by $V_t$). The conclusion then follows from Theorem \ref{mainthm-2nd} by letting $r=1$.
$(2)$ By assumption, the deformations of $(p+1,q)$-forms and $(p+1,q-1)$-forms on $X$ are canonically unobstructed, so the function $\dim V_{2,t}^{p,q}$ is upper semi-continuous by Proposition \ref{prop-upp-semicontinuous-V2tpq}. The conclusion then follows from Theorem \ref{mainthm-2nd} by letting $r=2$.
\end{proof}

\subsection{Unobstructedness for deformation of $r$-filtered $(p,q)$-forms}
\begin{thm}\label{thm-deformation of r-filtered}
Let $\pi: (\mathcal{X}, X)\to (\Delta,0)$ be a complex analytic family over a small disc $\Delta\subset \C$. Assume for some $r\geq1$,
\begin{itemize}
  \item the deformation of $r$-filtered $(p,q)$-forms on $X$ are canonically unobstructed;
  \item the deformation of $(p+r,q-r)$-forms on $X$ are canonically unobstructed;
  \item there holds
\begin{equation}\label{eq-deformation of r-filtered}
\bigoplus_{\lambda\geq i\geq 1}d_\lambda^{p+r-i,q-r+i}(X)=0.
\end{equation}
\end{itemize}
Then the deformation of $(r+1)$-filtered $(p,q)$-forms on $X$ are canonically unobstructed.
\end{thm}
\begin{proof}For any given $\sum_{i=0}^{r}\hat{\alpha}_0^{p+i,q-i}\in\ker d_{\Pi_{r+1}}$, we need to construct a canonical deformation of $\sum_{i=0}^{r}\hat{\alpha}_0^{p+i,q-i}$ on $\Delta$.

Indeed, since $\sum_{i=0}^{r}\hat{\alpha}_0^{p+i,q-i}\in\ker d_{\Pi_{r+1}}$ implies $\sum_{i=0}^{r-1}\hat{\alpha}_0^{p+i,q-i}\in\ker d_{\Pi_{r}}$, there is a canonical deformation $\sum_{i=0}^{r-1}\hat{\alpha}^{p+i,q-i}(t)$ of $\sum_{i=0}^{r-1}\hat{\alpha}_0^{p+i,q-i}$ on $\Delta$ where $\hat{\alpha}^{p+i,q-i}(t)=\sum_{k\geq0}\hat{\alpha}_k^{p+i,q-i}$ for $0\leq i\leq r-1$. The proof of Theorem \ref{mainthm-2nd} shows there is a convergent power series solution $\alpha(t)=\sum_{k\geq0}\alpha_k$ of \eqref{eq-canonical-alphat}. Now since $\hat{\alpha}_0^{p+r-1,q-r+1}=\alpha_0^{p+r-1,q-r+1}$ and
\[
\p\hat{\alpha}_0^{p+r-1,q-r+1}+\pb\hat{\alpha}_0^{p+r,q-r}=\p\alpha_0^{p+r-1,q-r+1}+\pb\alpha_0^{p+r,q-r}=0,
\]
we have $\hat{\alpha}_0^{p+r,q-r}-\alpha_0^{p+r,q-r}\in\ker\pb$. Let $\theta(t)=\sum_{k\geq0}\theta_k$ be the canonical Dolbeault deformation of $\hat{\alpha}_0^{p+r,q-r}-\alpha_0^{p+r,q-r}$, then by definition
\begin{equation}\label{eq-5}
\theta_0=\hat{\alpha}_0^{p+r,q-r}-\alpha_0^{p+r,q-r}\quad\text{and}\quad\pb_{\phi(t)}\theta(t)=0,~\text{for~any}~t\in\Delta.
\end{equation}
Because $\alpha(t)=\sum_{k\geq0}\alpha_k\in\ker d_{\phi(t)}$ for any $t\in\Delta$, we have
\[
0=\pb_{\phi(t)}\alpha^{p,q}(t)=\p\alpha^{p,q}(t)+\pb_{\phi(t)}\alpha^{p+1,q-1}(t)=\cdots=\p\alpha^{p+r-1,q-r+1}(t)+\pb_{\phi(t)}\alpha^{p+r,q-r}(t).
\]
Combining this with \eqref{eq-5}, we get
\begin{align*}
0=\pb_{\phi(t)}\alpha^{p,q}(t)=&\p\alpha^{p,q}(t)+\pb_{\phi(t)}\alpha^{p+1,q-1}(t)
\\=\cdots=&\p\alpha^{p+r-1,q-r+1}(t)+\pb_{\phi(t)}(\alpha^{p+r,q-r}(t)+\theta(t)).
\end{align*}
Note that $\alpha^{p+r,q-r}_0+\theta_0=\hat{\alpha}_0^{p+r,q-r}$, we see that $\sum_{i=0}^{r}\alpha^{p+i,q-i}(t)+\theta(t)$ is a canonical deformation of $\sum_{i=0}^{r}\hat{\alpha}_0^{p+i,q-i}\in\ker d_{\Pi_{r+1}}$.
\end{proof}

\begin{corollary}\label{coro-deformation of r-filtered-1}
Let $\pi: (\mathcal{X}, X)\to (\Delta,0)$ be a complex analytic family over a small disc $\Delta\subset \C$. Assume for some $r\geq2$,
\begin{itemize}
  \item the deformation of $(p+i,q-i)$-forms on $X$ are canonically unobstructed for $0\leq i\leq r-1$;
  \item there holds
\begin{equation}\label{eq-r-filtered-1}
\bigoplus_{ \substack{\lambda\geq 1,\\ 2\leq i\leq r} }d_\lambda^{p+r-i,q-r+i}(X)=0
=\bigoplus_{ \lambda\geq i\geq 2 }d_\lambda^{p+1-i,q-1+i}(X).
\end{equation}
\end{itemize}
Then the deformation of $r$-filtered $(p,q)$-forms on $X$ are canonically unobstructed.
\end{corollary}
\begin{proof}We notice that \eqref{eq-r-filtered-1} implies $\bigoplus_{\lambda\geq i\geq 1}d_\lambda^{p+k-i,q-k+i}(X)=0$ holds for $k=1,2,\cdots, r-1$. By applying Theorem \ref{thm-deformation of r-filtered} for $r=1$, we have the deformation of $2$-filtered $(p,q)$-forms on $X$ are canonically unobstructed. Then by applying Theorem \ref{thm-deformation of r-filtered} for $r=2$, we have the deformation of $3$-filtered $(p,q)$-forms on $X$ are canonically unobstructed. The conclusion follows by repeating this process.
\end{proof}

\begin{corollary}\label{coro-deformation of r-filtered-2}
Let $\pi: (\mathcal{X}, X)\to (\Delta,0)$ be a complex analytic family over a small disc $\Delta\subset \C$. Assume for some $r\geq1$, there holds
\begin{equation}
\bigoplus_{ \substack{\lambda\geq 1\\ 1\leq i\leq r+1}} d_\lambda^{p+r-i,q-r+i}(X)=0=\bigoplus_{\lambda\geq i\geq 2}d_\lambda^{p-i,q+i}(X).
\end{equation}
Then the deformation of $r$-filtered $(p,q)$-forms on $X$ are canonically unobstructed.
\end{corollary}
\begin{proof}Notice that for $r=1$, this is exactly the first conclusion of Theorem \ref{mainthm}. For $r\geq2$, this follows from Theorem \ref{mainthm} and Corollary \ref{coro-deformation of r-filtered-1}.
\end{proof}

\section{Examples}\label{sec-example}
In this section, we will provide some concrete examples of the degeneration of Fr\"olicher spectral sequence to examine the conditions in our results.

\subsection{The case of degeneration at $E_2$}
\begin{example}\label{example Case III-(2)}
Let $G$ be the matrix Lie group defined by
\[
G := \left\{
\left(
\begin{array}{ccc}
 1 & z^1 & z^3 \\
 0 &  1  & z^2 \\
 0 &  0  &  1
\end{array}
\right) \in \mathrm{GL}(3;\mathbb{C}) \mid z^1,\,z^2,\,z^3 \in\mathbb{C} \right\}\cong \mathbb{C}^3,
\]
where the product is the one induced by matrix multiplication. This is usually called the \emph{Heisenberg group}. Consider the discrete subgroup $\Gamma$ defined by
\[
\Gamma := \left\{
\left(
\begin{array}{ccc}
 1 & \omega^1 & \omega^3 \\
 0 &  1  & \omega^2 \\
 0 &  0  &  1
\end{array}
\right) \in G \mid \omega^1,\,\omega^2,\,\omega^3 \in\mathbb{Z}[\sqrt{-1}] \right\},
\]
The quotient $X=G/\Gamma$ is called the \emph{Iwasawa manifold}.
\end{example}
\subsubsection{The Kuranishi family of $X$}
A basis of $H^0(X,\Omega^1)$ is given by\footnote{$\{\varphi^1, \varphi^2, \varphi^3\}$ is also a basis of left invariant $(1,0)$-forms on $G$.}
\[
\varphi^1 = d z^1,~ \varphi^2 = d z^2,~ \varphi^3 = d z^3-z^1\,d z^2,
\]
and a dual basis $\theta^1, \theta^2, \theta^3\in H^0(X,T_X^{1,0})$ is given by
\[
\theta^1=\frac{\partial}{\partial z^1},~\theta^2=\frac{\partial}{\partial z^2} + z^1\frac{\partial}{\partial z^3},~\theta^3=\frac{\partial}{\partial z^3}.
\]
The Iwasawa manifold $X$ is equipped with the ($G$-left invariant) Hermitian metric $\sum_{i=1}^3\varphi^i\otimes\bar{\varphi}^i$ which induces a Hermitian metric on $\mathfrak{g}$ (the Lie algebra of $G$). The Beltrami differential of the Kuranishi family of $X$ is given by (c.f. \cite{Nak75})
\[
\phi(t) = \sum_{i=1}^3\sum_{\lambda=1}^2t_{i\lambda}\theta^i\bar{\varphi}^{\lambda} - D(t)\theta^3\bar{\varphi}^{3},~\text{with}~D(t)=t_{11}t_{22}-t_{21}t_{12},
\]
and the Kuranishi space of $X$ is
\[
\mathcal{B}=\{t=(t_{11}, t_{12}, t_{21}, t_{22}, t_{31}, t_{32})\in \mathbb{C}^6\mid |t_{i\lambda}|<\epsilon, i=1, 2, 3, \lambda=1,2 \},
\]
where $\epsilon>0$ is sufficiently small. Set
\[
\phi_1=\sum_{i=1}^3\sum_{\lambda=1}^2t_{i\lambda}\theta^i\bar{\varphi}^{\lambda},~ \phi_2 = -D(t)\theta^3\bar{\varphi}^{3}.
\]

\subsubsection{The Fr\"olicher spectral sequence $(E_r^{p,q},d_r^{p,q})$ of $X$}
It was shown \cite[Coro.\,4.3]{CFUG99} that for compact nilmanifolds with nilpotent complex structure, in particular the Iwasawa manifold and the manifold considered in Example \ref{example-E3}, the Fr\"olicher spectral sequence $(E_r^{p,q},d_r^{p,q})$ may be computed by using the left invariant forms only. In fact, for the Iwasawa manifold, $\dim d_r^{p,q}(X)$ has been computed in \cite{HX24D} (building on \cite{Ang13,Ste21,KQ20}) where we may find that
\begin{itemize}
  \item $d_{r}^{p,q}=0$ for any $(p,q)$ and $r\geq 2$;
  \item $d_{1}^{p,q}\neq 0$ if and only if $(p,q)=(1,0), (1,1), (1,2), (1,3)$.
\end{itemize}

\subsubsection{Computing the canonical deformations of $(p,q)$-forms}
\begin{description}
  \item[Obstructed cases] the deformations of $(p,q)$-forms are obstructed for
  \[(p,q)=(1,0), (2,0), (1,2), (2,2).\]

  \item[Unobstructed cases] the deformations of $(p,q)$-forms are unobstructed for
  \[(p,q)\neq (1,0), (2,0), (1,2), (2,2).\]
\end{description}
In order to see this, we compute the (canonical) deformations as in \cite[Sec.\,7]{Xia19dDol} where the cases $(p,q)=(1,0), (1,1), (2,0)$ has already been done. Here we only present the computations in the case $(p,q)=(1,2), (2,2), (2,1)$. In fact,
\[
H^{1,2}_{\pb}(X)=\mathbb{C}\{\varphi^{1\overline{23}}, \varphi^{2\overline{23}}, \varphi^{3\overline{23}}, \varphi^{1\overline{13}}, \varphi^{2\overline{13}}, \varphi^{3\overline{13}} \},
\]
we set $\sigma_0 = a_{123}\varphi^{1\overline{23}}+a_{223}\varphi^{2\overline{23}}+a_{323}\varphi^{3\overline{23}}+a_{113}\varphi^{1\overline{13}}+a_{213}\varphi^{2\overline{13}}+a_{313} \varphi^{3\overline{13}}$, then
\begin{align*}
\mathcal{L}_{\phi_1}^{1,0}\sigma_0
&= (a_{313}t_{22} - a_{323}t_{21})\varphi^{1\overline{123}} - (a_{313}t_{12} - a_{323}t_{11})\varphi^{2\overline{123}}
\end{align*}
is exact if and only if $\mathcal{L}_{\phi_1}^{1,0}\sigma_0=0$, i.e.
\begin{equation*}
\left\{
\begin{array}{rcl}
a_{313}t_{22} - a_{323}t_{21} &=& 0 \\[5pt]
a_{313}t_{12} - a_{323}t_{11} &=& 0.
\end{array}
\right.
\end{equation*}
This shows the deformation of $(1,2)$-forms are obstructed. Similarly, for
\[
H^{2,2}_{\pb}(X)=\mathbb{C}\{\varphi^{12\overline{23}}, \varphi^{23\overline{23}}, \varphi^{13\overline{23}}, \varphi^{12\overline{13}}, \varphi^{23\overline{13}}, \varphi^{13\overline{13}} \},
\]
we set $\sigma_0 = a_{1223}\varphi^{12\overline{23}}+a_{2323}\varphi^{23\overline{23}}+a_{1323}\varphi^{13\overline{23}}+a_{1213}\varphi^{12\overline{13}}+a_{2313} \varphi^{23\overline{13}}+a_{1313}\varphi^{13\overline{13}}$, then
\begin{align*}
\mathcal{L}_{\phi_1}^{1,0}\sigma_0
&= (a_{1313}t_{12}-a_{1323}t_{11}-a_{2323}t_{21} -a_{2313}t_{22})\varphi^{12\overline{123}}
\end{align*}
is exact if and only if $a_{1313}t_{12}-a_{1323}t_{11}-a_{2323}t_{21} -a_{2313}t_{22}=0$ which shows the deformation of $(2,2)$-forms are obstructed.
Finally, for
\[
H^{2,1}_{\pb}(X)=\mathbb{C}\{\varphi^{12\bar{1}}, \varphi^{23\bar{1}}, \varphi^{13\bar{1}}, \varphi^{12\bar{2}}, \varphi^{23\bar{2}}, \varphi^{13\bar{2}} \},
\]
we set $\sigma_0 = a_{121}\varphi^{12\bar{1}}+a_{231}\varphi^{23\bar{1}}+a_{131}\varphi^{13\bar{1}}+a_{122}\varphi^{12\bar{2}}+a_{232}\varphi^{23\bar{2}}+a_{132}\varphi^{13\bar{2}}$, then
\begin{align*}
\mathcal{L}_{\phi_1}^{1,0}\sigma_0
&= (a_{231}t_{22}+a_{131}t_{11}-a_{232}t_{21} -a_{131}t_{11})\varphi^{12\overline{12}}\\
&=-\pb(a_{231}t_{22}+a_{131}t_{11}-a_{232}t_{21} -a_{131}t_{11})\varphi^{12\overline{3}},
\end{align*}
and
\begin{align*}
\sigma_1&=\pb^*G_{\pb} \mathcal{L}_{\phi_1}^{1,0}\sigma_0\\
&=\pb^*G_{\pb}\pb(a_{232}t_{21}+a_{131}t_{11}-a_{231}t_{22}-a_{131}t_{11})\varphi^{12\overline{3}} \\
&=P_{(\ker\pb)^\perp}(a_{232}t_{21}+a_{131}t_{11}-a_{231}t_{22}-a_{131}t_{11})\varphi^{12\overline{3}}\\
&=(a_{232}t_{21}+a_{131}t_{11}-a_{231}t_{22}-a_{131}t_{11})\varphi^{12\overline{3}},
\end{align*}
where $P_{(\ker\pb)^\perp}$ is the projection operator onto $(\ker\pb)^\perp$. The last equality holds since $\varphi^{12\overline{3}}$ is orthogonal to $\ker\pb$. Furthermore,
\[
\mathcal{L}_{\phi_1}^{1,0}\sigma_1=\mathcal{L}_{\phi_2}^{1,0}\sigma_0=\mathcal{L}_{\phi_2}^{1,0}\sigma_1 = 0\Rightarrow \sigma_{k}=\pb^*G_{\pb}\sum_{j=1}^{k}\mathcal{L}_{\phi_j}^{1,0}\sigma_{k-j}=0,~ k>1.
\]
We see that the deformation of $(2,1)$-forms are unobstructed.
\subsubsection{Examinations of results}
Now we examine the sufficiency and necessity of our conditions
\begin{equation}\label{eq-condition}
\bigoplus_{r\geq 1}d_r^{p,q}(X)=0\quad\text{and}\quad
  \bigoplus_{ r\geq i\geq 1 }d_r^{p-i,q+i}(X)=0,
\end{equation}
in Theorem \ref{mainthm}.
\begin{description}
\item[Sufficiency] \eqref{eq-condition} is satisfied for $(p,q)\neq (1,0), (1,1), (1,2), (1,3), (2,0), (2,1), (2,2)$, i.e.
     \[
    (p,q)=(0,1),(0,2),(0,3),(3,0),(3,1),(3,2),(2,3),(3,3).
     \]
     By Theorem \ref{mainthm} the deformations of $(p,q)$-forms are unobstructed in these cases. In other (unobstructed) situations, \eqref{eq-condition} is not satisfied. In fact, the deformations of $(1,1)$-forms are unobstructed but $d_1^{1,1}\neq0$; $(2,1)$-forms: $d_1^{1,2}\neq 0$; $(1,3)$-forms: $d_1^{1,3}\neq 0$;
\item[Necessity] The obstructed cases shows $d_1^{p,q}=0$ and $d_1^{p-1,q+1}=0$ in \eqref{eq-condition} are necessary conditions for the unobstructed deformations of $(p,q)$-forms. In fact, for $(1,0)$-forms: $d_1^{1,0}\neq 0$; $(2,0)$-forms: $d_1^{1,1}\neq 0$; $(1,2)$-forms: $d_1^{1,2}\neq 0$; $(2,2)$-forms: $d_1^{1,3}\neq 0$.
\end{description}

On the other hand, by applying Theorem \ref{mainthm}, we get $\bigoplus_{\lambda\geq 1}d_\lambda^{p,q}(X_t)=0$ for any small $t$ when
\[
(p,q)=(0,1),(0,2),(0,3),(3,0),(3,1),(3,2),(2,3),(3,3).
\]
Furthermore, by applying Corollary \ref{coro-to-mainthm2}, $(1)$, we get $\bigoplus_{\lambda\geq 1}d_\lambda^{p,q}(X_t)=0$ for any small $t$ when (in these cases Theorem \ref{mainthm} is not applicable)
\[
(p,q)=(2,1).
\]
It is known that for complex parallelizable nilmanifold (including the Iwasawa manifold $X$) there holds $E_2\cong E_\infty$ (see \cite[Thm.\,9]{CFG91}). When $(t_{11}, t_{12}, t_{21}, t_{22})\neq0$, the small deformation $X_t$ is not complex parallelizable, our results give some useful information about the degenerating behavior of $d_r^{p,q}(X_t)$ (c.f. \cite{Fla20,HX24D}).

\subsection{The case of degeneration at $E_3$}
\begin{example}\label{example-E3}
Let $X=(\Gamma\setminus G, J)$ be a complex nilmanifold endowed with a left invariant complex structure $J$ such that the underlying Lie algebra $\g\cong \mathfrak{h}_{15}$ and the structure equation is given by
\begin{equation}\label{eq-streq-h15}
\left\{
\begin{array}{rcl}
d\omega^1 &=& 0,  \\
d\omega^2 &=& \omega^{1\bar{1}},  \\
d\omega^3 &=& \omega^{12}+3\omega^{1\bar{2}}+\omega^{2\bar{1}},
\end{array}
\right.
\end{equation}
where $\C\{\omega^1, \omega^2, \omega^3\}=\g^{1,0*}$ is the dual basis of $\C\{e_1, e_2, e_3\}=\g^{1,0}$. See \cite[pp.\,267]{COUV16} for more information.
\end{example}

\subsubsection{The Kuranishi family of $X$}
Because
\begin{align*}
[x,y]
&=\sum_{j=1}^3\omega^j[x,y]e_j+\sum_{j=1}^3\omega^{\bar{j}}[x,y]\overline{e_j}\\
&=-\sum_{j=1}^3d\omega^j(x,y)e_j-\sum_{j=1}^3d\omega^{\bar{j}}(x,y)\overline{e_j},\quad\forall x,y\in\g_\C,
\end{align*}
so the non-zero brackets are
\[
[e_1,\overline{e_1}]=-e_2+\overline{e_2},\quad [e_1,e_2]=-e_3,\quad [e_1,\overline{e_2}]=-3e_3+\overline{e_3},\quad [e_2,\overline{e_1}]=-e_3+3\overline{e_3}.
\]
We have
\[
\pb e_1=-\sum_{j=1}^3\omega^{\bar{j}}\ot[\overline{e_j},e_1]^{1,0}=-\omega^{\bar{1}}\ot e_2-3\omega^{\bar{2}}\ot e_3,~\pb e_2=-\omega^{\bar{1}}\ot e_3,~\pb e_3=0
\]
which implies that
\begin{align*}
H^{0,1}_{\pb}(X,T_X^{1,0})=&\frac{\C\{\omega^{\bar{1}}\ot e_2, \omega^{\bar{1}}\ot e_3, \omega^{\bar{2}}\ot e_3, \omega^{\bar{3}}\ot e_3+\omega^{\bar{2}}\ot e_2,
3\omega^{\bar{3}}\ot e_3-\omega^{\bar{1}}\ot e_1\}}{\C\{\omega^{\bar{1}}\ot e_2+3\omega^{\bar{2}}\ot e_3, \omega^{\bar{1}}\ot e_3 \} }\\
\cong &\C\{\omega^{\bar{1}}\ot e_2, \omega^{\bar{3}}\ot e_3+\omega^{\bar{2}}\ot e_2,
3\omega^{\bar{3}}\ot e_3-\omega^{\bar{1}}\ot e_1\}.
\end{align*}
We compute the following
\begin{align*}
[\omega^{\bar{1}}\ot e_2, \omega^{\bar{1}}\ot e_2]&=0,\quad [\omega^{\bar{3}}\ot e_3, \omega^{\bar{3}}\ot e_3]=0, \\
[\omega^{\bar{1}}\ot e_2, \omega^{\bar{3}}\ot e_3+\omega^{\bar{2}}\ot e_2]&=0,\quad [\omega^{\bar{2}}\ot e_2, \omega^{\bar{1}}\ot e_1]=-\omega^{\overline{12}}\ot e_3,\\
[\omega^{\bar{1}}\ot e_2, 3\omega^{\bar{3}}\ot e_3-\omega^{\bar{1}}\ot e_1]&=0,\quad [\omega^{\bar{2}}\ot e_2, \omega^{\bar{2}}\ot e_2]=0,\\
[\omega^{\bar{3}}\ot e_3+\omega^{\bar{2}}\ot e_2, \omega^{\bar{3}}\ot e_3+\omega^{\bar{2}}\ot e_2]&=6\omega^{\overline{12}}\ot e_3,\quad [\omega^{\bar{3}}\ot e_3, \omega^{\bar{1}}\ot e_1]=-\omega^{\overline{12}}\ot e_3,\\
[\omega^{\bar{3}}\ot e_3+\omega^{\bar{2}}\ot e_2, 3\omega^{\bar{3}}\ot e_3-\omega^{\bar{1}}\ot e_1]&=9\omega^{\overline{12}}\ot e_3,\\
[3\omega^{\bar{3}}\ot e_3-\omega^{\bar{1}}\ot e_1, 3\omega^{\bar{3}}\ot e_3-\omega^{\bar{1}}\ot e_1]&=6\omega^{\overline{12}}\ot e_3.
\end{align*}
Let $\phi_1=t_1\omega^{\bar{1}}\ot e_2+t_2(\omega^{\bar{3}}\ot e_3+\omega^{\bar{2}}\ot e_2)+t_3(3\omega^{\bar{3}}\ot e_3-\omega^{\bar{1}}\ot e_1)$, then
\[
[\phi_1,\phi_1]=(6t_2^2+18t_2t_3+6t_3^2)\omega^{\overline{12}}\ot e_3.
\]
Note that $\omega^{\overline{12}}\ot e_3=\pb(\omega^{\overline{3}}\ot e_3)$, we get
\[
\phi_2:=(3t_2^2+9t_2t_3+3t_3^2)\omega^{\overline{3}}\ot e_3\Longrightarrow \pb\phi_2=\frac{1}{2}[\phi_1,\phi_1].
\]
From
\[
[\phi_1,\phi_2]=(3t_2^2+9t_2t_3+3t_3^2)(3t_2+t_3)\omega^{\overline{12}}\ot e_3,
\]
we get
\[
\phi_3:=(3t_2^2+9t_2t_3+3t_3^2)(3t_2+t_3)\omega^{\overline{3}}\ot e_3\Longrightarrow \pb\phi_3=\frac{1}{2}([\phi_1,\phi_2]+[\phi_2,\phi_1]),
\]
and similarly,
\[
\phi_4:=(3t_2^2+9t_2t_3+3t_3^2)(3t_2+t_3)^2\omega^{\overline{3}}\ot e_3\Longrightarrow \pb\phi_4=\frac{1}{2}\sum_{j=1}^{3}[\phi_j,\phi_{4-j}].
\]
Continuing in this way and notice that
\[
[\phi_j,\phi_k]=0,~\forall j,k\geq2\Longrightarrow \frac{1}{2}\sum_{j=1}^{k-1}[\phi_j,\phi_{k-j}]=[\phi_1,\phi_{k-1}],~\forall k\geq3,
\]
we get a complex analytic family (isomorphic to the Kuranishi family of $X$) whose Beltrami differentials are given by
\begin{equation}\label{eq-a family}
\left\{
\begin{array}{rcl}
\phi(t)&=&\phi_1+\phi_2+\phi_3+\cdots,\\
 \phi_1&=&t_1\omega^{\bar{1}}\ot e_2+t_2(\omega^{\bar{3}}\ot e_3+\omega^{\bar{2}}\ot e_2)+t_3(3\omega^{\bar{3}}\ot e_3-\omega^{\bar{1}}\ot e_1),\\
 \phi_{j}&=& (3t_2^2+9t_2t_3+3t_3^2)(3t_2+t_3)^{j-2}\omega^{\overline{3}}\ot e_3,~\forall j\geq 2.
\end{array}
\right.
\end{equation}

\subsubsection{The Fr\"olicher spectral sequence $(E_r^{p,q},d_r^{p,q})$ of $X$}
Based on the following computations,
\begin{itemize}
  \item $\pb\omega^1=\pb\omega^{\bar{1}}=\pb\omega^{\bar{2}}=0$; $\pb\omega^2=\omega^{1\bar{1}}$; $\pb\omega^3=3\omega^{1\bar{2}}+\omega^{2\bar{1}}$; $\pb\omega^{\bar{3}}=\omega^{\overline{12}}$; $\pb\omega^{12}=0$; $\pb\omega^{23}=-\omega^{13\bar{1}}+3\omega^{12\bar{2}}$; $\pb\omega^{13}=-\omega^{12\bar{1}}$;
  \item $\pb\omega^{1\bar{1}}=\pb\omega^{1\bar{2}}=\pb\omega^{2\bar{1}}=0$; $\pb\omega^{1\bar{3}}=-\omega^{1\overline{12}}=-\pb\omega^{2\bar{2}}$; $\pb\omega^{2\bar{3}}=\omega^{1\overline{13}}-\omega^{2\overline{12}}$; $\pb\omega^{3\bar{1}}=-3\omega^{1\overline{12}}$; $\pb\omega^{3\bar{2}}=\omega^{2\overline{12}}$; $\pb\omega^{3\bar{3}}=3\omega^{1\overline{23}}+\omega^{2\overline{13}}-\omega^{3\overline{12}}$; $\pb\omega^{\overline{12}}=\pb\omega^{\overline{13}}=\pb\omega^{\overline{23}}=0$;
  \item $\pb\omega^{12\bar{1}}=\pb\omega^{12\bar{2}}=\pb\omega^{13\bar{1}}=0$; $\pb\omega^{12\bar{3}}=\omega^{12\overline{12}}=-\frac{1}{3}\pb\omega^{23\bar{1}}$; $\pb\omega^{13\bar{2}}=-\omega^{12\overline{12}}$; $\pb\omega^{13\overline{3}}=\omega^{13\overline{12}}-\omega^{12\overline{13}}$; $\pb\omega^{23\overline{2}}=-\omega^{13\overline{12}}$; $\pb\omega^{23\overline{3}}=-\omega^{13\overline{13}}+3\omega^{12\overline{23}}+\omega^{23\overline{12}}$;
  \item $\pb\omega^{1\overline{12}}=\pb\omega^{1\overline{13}}=\pb\omega^{1\overline{23}}=\pb\omega^{2\overline{12}}=\pb\omega^{2\overline{13}}=\pb\omega^{3\overline{12}}=0$; $\pb\omega^{2\overline{23}}=\omega^{1\overline{123}}=-\frac{1}{3}\pb\omega^{3\overline{13}}$; $\pb\omega^{3\overline{23}}=\omega^{2\overline{123}}$;
  \item $\pb\omega^{123\overline{1}}=\pb\omega^{123\overline{2}}=0$; $\pb\omega^{123\overline{3}}=-\omega^{123\overline{12}}$;
  \item $\pb\omega^{12\overline{12}}=\pb\omega^{12\overline{13}}=\pb\omega^{12\overline{23}}=\pb\omega^{13\overline{12}}=\pb\omega^{13\overline{13}}=\pb\omega^{23\overline{12}}=0$; $\pb\omega^{13\overline{23}}=-\omega^{12\overline{123}}$; $\pb\omega^{23\overline{13}}=-3\omega^{12\overline{123}}, \pb\omega^{23\overline{23}}=-\omega^{13\overline{123}}$;
  \item $\pb\omega^{1\overline{123}}=\pb\omega^{2\overline{123}}=\pb\omega^{1\overline{123}}=0$;
  \item $\pb\omega^{123\overline{12}}=\pb\omega^{123\overline{13}}=\pb\omega^{123\overline{23}}=0$;
\end{itemize}
A basis of $\wt{Z}_1^{p,q}$ may be listed as follows

\renewcommand\arraystretch{1.5}
\begin{table}[!htbp]
\centering
\begin{center}
\begin{tabular}{|c|c|c|c|}
\hline
$\wt{Z}_1^{1,0}$ & $\omega^1;$ & $\wt{Z}_1^{0,1}$ & $\omega^{\bar{1}}, \omega^{\bar{2}}; $  \\
\hline
$\wt{Z}_1^{2,0}$ & $\omega^{12};$&$\wt{Z}_1^{1,1}$ & $\omega^{1\bar{2}}, 3\omega^{1\bar{3}}-\omega^{3\bar{1}}, 3\omega^{2\bar{2}}+\omega^{3\bar{1}};$  \\
 &  &  & $\omega^{1\bar{1}}=\pb\omega^2, 3\omega^{1\bar{2}}+\omega^{2\bar{1}}=\pb\omega^3$ \\
\hline
$\wt{Z}_1^{0,2}$ & $\omega^{\overline{13}}, \omega^{\overline{23}};$  &   & \\
 & $\omega^{\overline{12}}=\pb\omega^{\bar{3}}$ &  &  \\
\hline
$\wt{Z}_1^{3,0}$ & $\omega^{123};$ &  $\wt{Z}_1^{2,1}$ & $\omega^{12\bar{2}}, 3\omega^{12\bar{3}}+\omega^{23\bar{1}}, \omega^{13\bar{2}}+\omega^{12\bar{3}};$   \\
 &  &  & $3\omega^{12\bar{2}}-\omega^{13\bar{1}}=\pb\omega^{23}, -\omega^{12\bar{1}}=\pb\omega^{13}$; \\
\hline
$\wt{Z}_1^{0,3}$ & $\omega^{\overline{123}};$ & $\wt{Z}_1^{1,2}$  & $\omega^{1\overline{23}}, \omega^{2\overline{13}},  3\omega^{2\overline{23}}+\omega^{3\overline{13}}; \omega^{1\overline{13}}=\pb(\omega^{2\bar{3}}+\omega^{3\bar{2}})$, \\
 &  &  & $ \omega^{1\overline{12}}=-\pb\omega^{1\bar{3}}=\pb\omega^{2\bar{2}}=-\frac{1}{3}\pb\omega^{3\bar{1}}, $ \\
 &  &   & $\omega^{2\overline{12}}=\pb\omega^{3\bar{2}}, 3\omega^{1\overline{23}}+\omega^{2\overline{13}}-\omega^{3\overline{12}}=\pb\omega^{3\overline{3}}$ \\
\hline
$\wt{Z}_1^{3,1}$ & $\omega^{123\overline{1}}, \omega^{123\overline{2}}; $ & $\wt{Z}_1^{2,2}$  & $\omega^{12\overline{23}}, \omega^{13\overline{13}}, 3\omega^{13\overline{23}}-\omega^{23\overline{13}}; $ \\
&  &  & $\omega^{12\overline{12}}=\pb\omega^{12\overline{3}}=-\frac{1}{3}\pb\omega^{23\overline{1}}=-\pb\omega^{13\overline{2}}$, \\
&  &  & $\omega^{13\overline{12}}=-\pb\omega^{23\overline{2}}, \omega^{12\overline{13}}=-\pb(\omega^{23\overline{2}}+\omega^{13\overline{3}}),$\\
 &  &   & $\omega^{23\overline{12}}-\omega^{13\overline{13}}+3\omega^{12\overline{23}}=\pb\omega^{23\bar{3}}$ \\
\hline
$\wt{Z}_1^{1,3}$ & $\omega^{3\overline{123}}; \omega^{2\overline{123}}=\pb\omega^{3\overline{23}},$ &   &  \\
 &$\omega^{1\overline{123}}=\pb\omega^{2\overline{23}}=-\frac{1}{3}\pb\omega^{3\overline{13}}$   &  &  \\
\hline
$\wt{Z}_1^{3,2}$ & $\omega^{123\overline{23}}, \omega^{123\overline{13}};$ & $\wt{Z}_1^{2,3}$  & $\omega^{23\overline{123}};$ \\
 & $\omega^{123\overline{12}}=-\pb\omega^{123\overline{3}}$ &  & $\omega^{12\overline{123}}=\frac{1}{2}\pb\omega^{23\overline{13}}, \omega^{13\overline{123}}=-\pb\omega^{23\overline{23}}$\\
\hline
$\wt{Z}_1^{3,3}$ & $\omega^{123\overline{123}};$ &   &  \\
\hline
\end{tabular}
\end{center}
\end{table}

\begin{itemize}
  \item $d\omega^1=0\Rightarrow\omega^1\in\wt{Z}_\infty^{1,0}$; $d\omega^{\bar{1}}=0\Rightarrow\omega^{\bar{1}}\in\wt{Z}_\infty^{0,1}$; $0=\p\omega^{\bar{2}}+\pb\omega^2=\p\omega^2\Rightarrow\omega^{\bar{2}}\in\wt{Z}_\infty^{0,1}$.
  \item $d\omega^{12}=0\Rightarrow\omega^{12}\in\wt{Z}_\infty^{2,0}$.
  \item $d\omega^{1\bar{2}}=0\Rightarrow\omega^{1\bar{2}}\in\wt{Z}_\infty^{1,1}$; $\p(3\omega^{1\bar{3}}-\omega^{3\bar{1}})=8\omega^{12\bar{1}}=-8\pb\omega^{13}$ and $\p\omega^{13}=0\Rightarrow(3\omega^{1\bar{3}}-\omega^{3\bar{1}})\in \wt{Z}_\infty^{1,1}$; $\p(3\omega^{2\bar{2}}+\omega^{3\bar{1}})=-2\omega^{12\bar{1}}=2\pb\omega^{13}\Rightarrow(3\omega^{2\bar{2}}+\omega^{3\bar{1}})\in\wt{Z}_\infty^{1,1}$;
  \item $\p\omega^{\overline{13}}=\pb\omega^{1\bar{3}}, \p\omega^{1\bar{3}}=-3\pb\omega^{13}$ and $\p\omega^{13}=0\Rightarrow\omega^{\overline{13}}\in\wt{Z}_\infty^{0,2}$; $\p\omega^{\overline{23}}=-\omega^{1\overline{13}}+3\omega^{2\overline{12}}=\pb(2\omega^{3\bar{2}}-\omega^{2\bar{3}})\Rightarrow\omega^{\overline{23}}\in\wt{Z}_2^{0,2}$, but $\p(2\omega^{3\bar{2}}-\omega^{2\bar{3}}+\ker\pb)=\pb\omega^{23}-\omega^{13\bar{1}}+\p(\ker\pb)\notin\im\pb\Rightarrow\omega^{\overline{23}}\notin\wt{Z}_3^{0,2}$.
  \item $d\omega^{123}=0\Rightarrow\omega^{123}\in\wt{Z}_\infty^{3,0}$; $d\omega^{\overline{123}}=0\Rightarrow\omega^{\overline{123}}\in\wt{Z}_\infty^{0,3}$.
  \item $d\omega^{12\bar{2}}=d\omega^{13\bar{1}}=d(3\omega^{12\bar{3}}+\omega^{23\bar{1}})=d(\omega^{13\bar{2}}+\omega^{12\bar{3}})=0\Rightarrow
      \omega^{12\bar{2}}, \omega^{13\bar{1}}, 3\omega^{12\bar{3}}+\omega^{23\bar{1}}, \omega^{13\bar{2}}+\omega^{12\bar{3}}\in\wt{Z}_\infty^{0,3}$.
  \item $\p\omega^{1\overline{23}}=-3\pb\omega^{12\overline{3}}$ and $\p\omega^{12\overline{3}}=0\Rightarrow\omega^{1\overline{23}}\in\wt{Z}_\infty^{1,2}$;
       $\p\omega^{2\overline{13}}=-\pb\omega^{12\overline{3}}$ and $\p\omega^{12\overline{3}}=0\Rightarrow\omega^{2\overline{13}}\in\wt{Z}_\infty^{1,2}$; $\p(3\omega^{2\overline{23}}+\omega^{3\overline{13}})=\pb(2\omega^{13\overline{3}}+3\omega^{23\overline{2}})\Rightarrow\omega^{2\overline{13}}\in\wt{Z}_2^{1,2}$,
       but $\p(2\omega^{13\overline{3}}+3\omega^{23\overline{2}}+\ker\pb)=3\omega^{123\overline{1}}+\p(\ker\pb)\notin\im\pb\Rightarrow 3\omega^{2\overline{23}}+\omega^{3\overline{13}}\notin\wt{Z}_3^{1,2}$;
  \item $d\omega^{123\overline{1}}=d\omega^{123\overline{2}}=0\Rightarrow
      \omega^{123\overline{1}}, \omega^{123\overline{2}}\in\wt{Z}_\infty^{3,1}$;
  \item $d\omega^{12\overline{23}}=d\omega^{13\overline{13}}=0$; $\p(3\omega^{13\overline{23}}-\omega^{23\overline{13}})=-8\omega^{123\overline{12}}=8\pb\omega^{123\overline{3}}$ and $\p\omega^{123\overline{3}}=0\Rightarrow 3\omega^{13\overline{23}}-\omega^{23\overline{13}}\in\wt{Z}_\infty^{2,2}$;
  \item $d\omega^{2\overline{123}}=0\Rightarrow\omega^{2\overline{123}}\in\wt{Z}_\infty^{1,3}$;
  \item $d\omega^{123\overline{23}}=d\omega^{123\overline{13}}=0\Rightarrow\omega^{123\overline{23}},\omega^{123\overline{13}}\in\wt{Z}_\infty^{3,2}$;
  \item $d\omega^{23\overline{123}}=0\Rightarrow\omega^{23\overline{123}}\in\wt{Z}_\infty^{2,3}$;
\end{itemize}
We see that
\begin{itemize}
\item $d_1^{p,q}=0$ for any $(p,q)$;
\item $0=d_2^{p,q}=d_3^{p,q}=\cdots$ for $(p,q)\neq(0,2), (1,2)$;
\item $d_2^{p,q}\neq0$ and $0=d_3^{p,q}=d_4^{p,q}=\cdots$ for $(p,q)=(0,2), (1,2)$.
\end{itemize}
In summary, we have
\begin{equation}\label{eq-dr=0examp-2}
d_r^{p,q}\neq0\Longleftrightarrow r=2~\text{and}~(p,q)=(0,2), (1,2).
\end{equation}

\subsubsection{Computing the canonical deformations of $(p,q)$-forms}
As a result of \eqref{eq-dr=0examp-2}, we see that the conditions
$\bigoplus_{r\geq 1}d_r^{p,q}(X)=0=
  \bigoplus_{ r\geq i\geq 1 }d_r^{p-i,q+i}(X)$
of Theorem \ref{mainthm} is satisfied for all $(p,q)\neq (0,2), (1,2), (2,0), (3,0)$. So the deformations of $(p,q)$-forms are canonically unobstructed for all $(p,q)\neq (0,2), (1,2), (2,0), (3,0)$.

In order to compute deformations of $(p,q)$-forms, we first compute the Lie derivatives:
\begin{itemize}
\item $\mathcal{L}_{\phi_1}^{1,0}\omega^{1}=\mathcal{L}_{\phi_1}^{1,0}\omega^{\bar{1}}=\mathcal{L}_{\phi_1}^{1,0}\omega^{\bar{2}}=0$;
\item $\mathcal{L}_{\phi_1}^{1,0}\omega^{2}=t_2\omega^{1\bar{1}}=t_2\pb\omega^{2}$;\quad $\mathcal{L}_{\phi_1}^{1,0}\omega^{\bar{3}}=(3t_2+t_3)\omega^{\overline{12}}$; \\ $\mathcal{L}_{\phi_1}^{1,0}\omega^{3}=t_1\omega^{1\bar{1}}+(3t_3+2t_2)\omega^{1\bar{2}}+(10t_3+3t_2)\omega^{2\bar{1}}$;
\item $\mathcal{L}_{\phi_1}^{1,0}\omega^{12}=\mathcal{L}_{\phi_1}^{1,0}\omega^{1\bar{2}}=0$;\quad $\mathcal{L}_{\phi_1}^{1,0}\omega^{1\bar{3}}=-(3t_2+t_3)\omega^{1\overline{12}}$;\quad $\mathcal{L}_{\phi_1}^{1,0}\omega^{2\bar{2}}=t_2\omega^{1\overline{12}};$ \\ $\mathcal{L}_{\phi_1}^{1,0}\omega^{3\bar{1}}=-(3t_3+2t_2)\omega^{1\overline{12}};$
\item $\mathcal{L}_{\phi_1}^{1,0}\omega^{12\bar{2}}=\mathcal{L}_{\phi_1}^{1,0}\omega^{13\bar{1}}=0$;\quad $\mathcal{L}_{\phi_1}^{1,0}\omega^{12\bar{3}}=(3t_2+t_3)\omega^{12\overline{12}}$; \\ $\mathcal{L}_{\phi_1}^{1,0}\omega^{23\bar{1}}=-(3t_3+2t_2)\omega^{12\overline{12}}$;\quad $\mathcal{L}_{\phi_1}^{1,0}\omega^{13\bar{2}}=-(10t_3+3t_2)\omega^{12\overline{12}}$;
\item $\mathcal{L}_{\phi_1}^{1,0}\omega^{1\overline{23}}=\mathcal{L}_{\phi_1}^{1,0}\omega^{2\overline{13}}=0$;\quad $\mathcal{L}_{\phi_1}^{1,0}\omega^{3\overline{13}}=-(3t_3+2t_2)\omega^{1\overline{123}}$;
\quad $\mathcal{L}_{\phi_1}^{1,0}\omega^{2\overline{23}}=t_2\omega^{1\overline{123}}$;
\item $\mathcal{L}_{\phi_2}^{1,0}\omega^{1}=\mathcal{L}_{\phi_2}^{1,0}\omega^{2}=0$,\quad $\mathcal{L}_{\phi_2}^{1,0}\omega^{3}=A(3\omega^{2\bar{1}}+\omega^{1\bar{2}})$ where $A:=3t_2^2+9t_2t_3+3t_3^2$;
\item  $\mathcal{L}_{\phi_2}^{1,0}\omega^{\bar{1}}=\mathcal{L}_{\phi_2}^{1,0}\omega^{\bar{2}}=\mathcal{L}_{\phi_2}^{1,0}\omega^{\bar{3}}=0$;
\item  $\mathcal{L}_{\phi_2}^{1,0}\omega^{1\bar{2}}=\mathcal{L}_{\phi_2}^{1,0}\omega^{1\bar{3}}=\mathcal{L}_{\phi_2}^{1,0}\omega^{2\bar{2}}=0$; $\mathcal{L}_{\phi_2}^{1,0}\omega^{3\bar{1}}=-A\omega^{1\overline{12}}$;
\item For $j\geq 2$, note that $\phi_j=(3t_2+t_3)^{j-2}\phi_2$.
\end{itemize}

Equip $X$ with the Hermitian metric $\sum_{i=1}^3\omega^{i\bar{i}}$, we are now ready to compute the canonical deformations of $(p,q)$-forms for the complex analytic family defined by \eqref{eq-a family}. We will only do this for $(p,q)=(0,2), (1,2), (2,0), (3,0), (1,1), (2,1)$.

For $(p,q)=(2,0)$. Let $\sigma_0=\omega^{12}$, then $0=\mathcal{L}_{\phi_1}^{1,0}\sigma_0=\mathcal{L}_{\phi_2}^{1,0}\sigma_0=\cdots$ implies the deformation of $(2,0)$-forms are unobstructed.

For $(p,q)=(1,1)$. Let $\sigma_0=a_1\omega^{1\bar{2}}+a_2(3\omega^{1\bar{3}}-\omega^{3\bar{1}})+a_3(3\omega^{2\bar{2}}+\omega^{3\bar{1}})$, then
\[
\mathcal{L}_{\phi_1}^{1,0}\sigma_0=-7a_2t_2\omega^{1\overline{12}}+a_3(t_2-3t_3)\omega^{1\overline{12}}=[a_3(t_2-3t_3)-7a_2t_2]\pb\omega^{1\bar{3}}.
\]
So
\begin{align*}
\sigma_1&=\pb^*G_{\pb}\mathcal{L}_{\phi_1}^{1,0}\sigma_0=[a_3(t_2-3t_3)-7a_2t_2]\pb^*G_{\pb}\pb\omega^{1\bar{3}}\\
&=[a_3(t_2-3t_3)-7a_2t_2]P_{(\ker\pb)^\bot}\omega^{1\bar{3}}\\
&=\frac{a_3(t_2-3t_3)-7a_2t_2}{\sqrt{11}}(\omega^{1\bar{3}}-\omega^{2\bar{2}}+3\omega^{3\bar{1}}),
\end{align*}
where we used the fact that $(\ker\pb)^\bot\cap A^{1,1}(X)=\C\{\omega^{2\bar{3}}, \omega^{3\bar{2}}, \omega^{3\bar{3}}, \omega^{1\bar{3}}-\omega^{2\bar{2}}+3\omega^{3\bar{1}}\}$. Similarly, we have
\[
\mathcal{L}_{\phi_1}^{1,0}\sigma_1=-\frac{a_3(t_2-3t_3)-7a_2t_2}{\sqrt{11}}\times 10(t_2+t_3)\omega^{1\overline{12}},\quad \mathcal{L}_{\phi_2}^{1,0}\sigma_0=A(a_2-a_3)\omega^{1\overline{12}},
\]
which implies $\sigma_2=\pb^*G_{\pb}(\mathcal{L}_{\phi_2}^{1,0}\sigma_0+\mathcal{L}_{\phi_1}^{1,0}\sigma_1)=p_2(\omega^{1\bar{3}}-\omega^{2\bar{2}}+3\omega^{3\bar{1}})$, where $p_2$ is a (homogeneous) polynomial of degree $2$ in $(t_2,t_3)$. Now because $\mathcal{L}_{\phi_2}^{1,0}(\omega^{1\bar{3}}-\omega^{2\bar{2}}+3\omega^{3\bar{1}})=-3A\omega^{1\overline{12}}$, we have $\sigma_3=\pb^*G_{\pb}(\mathcal{L}_{\phi_3}^{1,0}\sigma_0+\mathcal{L}_{\phi_2}^{1,0}\sigma_1+\mathcal{L}_{\phi_1}^{1,0}\sigma_2)=p_3(\omega^{1\bar{3}}-\omega^{2\bar{2}}+3\omega^{3\bar{1}})$ where $p_3$ is a polynomial of degree $3$ in $(t_2,t_3)$. Since $\phi_j=(3t_2+t_3)^{j-2}\phi_2$ for $j\geq 2$, this process can always be continued, we conclude that the deformation of $(1,1)$-forms are unobstructed.

For $(p,q)=(0,2)$. Let $\sigma_0=a_1\omega^{\overline{13}}+a_2\omega^{\overline{23}}$, then
\[
\mathcal{L}_{\phi_1}^{1,0}\sigma_0=0\Rightarrow
\sigma_1=\pb^*G_{\pb}\mathcal{L}_{\phi_1}^{1,0}\sigma_0=0.
\]
It can be checked easily that $\sigma_k=\sum_{1\leq j\leq k}\pb^*G_{\pb}\mathcal{L}_{\phi_j}^{1,0}\sigma_{k-j}=0$ for any $k\geq 2$. It follows that the deformation of $(0,2)$-forms are unobstructed.

For $(p,q)=(3,0)$. Let $\sigma_0=\omega^{123}$, then
\[
\mathcal{L}_{\phi_1}^{1,0}\sigma_0=0\Rightarrow
\sigma_1=\pb^*G_{\pb}\mathcal{L}_{\phi_1}^{1,0}\sigma_0=0.
\]
It can be checked easily that $\sigma_k=\sum_{1\leq j\leq k}\pb^*G_{\pb}\mathcal{L}_{\phi_j}^{1,0}\sigma_{k-j}=0$ for any $k\geq 2$. It follows that the deformation of $(3,0)$-forms are unobstructed.

For $(p,q)=(2,1)$. Let $\sigma_0=a_1\omega^{12\bar{2}}+a_2(3\omega^{12\bar{3}}+\omega^{23\bar{1}})+a_3(\omega^{13\bar{2}}+\omega^{12\bar{3}})$, then
\begin{align*}
\mathcal{L}_{\phi_1}^{1,0}\sigma_0=&(7a_2t_2-9a_3t_3)\omega^{12\overline{12}}=(7a_2t_2-9a_3t_3)\pb\omega^{12\overline{3}}\\
\Longrightarrow\sigma_1=&\pb^*G_{\pb}\mathcal{L}_{\phi_1}^{1,0}\sigma_0=\frac{7a_2t_2-9a_3t_3}{\sqrt{11}}(\omega^{12\overline{3}}-\omega^{13\overline{2}}-3\omega^{23\overline{1}}),
\end{align*}
where we used the fact that $(\ker\pb)^\bot\cap A^{2,1}(X)=\C\{\omega^{13\bar{3}}, \omega^{23\bar{2}}, \omega^{23\bar{3}}, \omega^{12\overline{3}}-\omega^{13\overline{2}}-3\omega^{23\overline{1}}\}$. Similarly, we have
\[
\mathcal{L}_{\phi_1}^{1,0}\sigma_1=\frac{(7a_2t_2-9a_3t_3)(12t_2+20t_3)}{\sqrt{11}}\omega^{12\overline{12}},\quad \mathcal{L}_{\phi_2}^{1,0}\sigma_0=-A(a_2+3a_3)\omega^{12\overline{12}},
\]
which implies $\sigma_2=\pb^*G_{\pb}(\mathcal{L}_{\phi_2}^{1,0}\sigma_0+\mathcal{L}_{\phi_1}^{1,0}\sigma_1)=p_2(\omega^{12\overline{3}}-\omega^{13\overline{2}}-3\omega^{23\overline{1}})$, where $p_2$ is a polynomial of degree $2$ in $(t_2,t_3)$. Now because $\mathcal{L}_{\phi_2}^{1,0}(\omega^{12\overline{3}}-\omega^{13\overline{2}}-3\omega^{23\overline{1}})=6A\omega^{12\overline{12}}$, we have $\sigma_2=\pb^*G_{\pb}(\mathcal{L}_{\phi_3}^{1,0}\sigma_0+\mathcal{L}_{\phi_2}^{1,0}\sigma_1+\mathcal{L}_{\phi_1}^{1,0}\sigma_2)
=p_3(\omega^{12\overline{3}}-\omega^{13\overline{2}}-3\omega^{23\overline{1}})$ where $p_3$ is a polynomial of degree $3$ in $(t_2,t_3)$. Since $\phi_j=(3t_2+t_3)^{j-2}\phi_2$ for $j\geq 2$, this process can always be continued, we conclude that the deformation of $(2,1)$-forms are unobstructed.

For $(p,q)=(1,2)$. Let $\sigma_0=a_1\omega^{1\overline{23}}+a_2\omega^{2\overline{13}}+a_3(3\omega^{2\overline{23}}+\omega^{3\overline{13}})$, then
\[
\mathcal{L}_{\phi_1}^{1,0}\sigma_0=-a_3(3t_3+2t_2)\omega^{1\overline{123}}=-a_3(3t_3+2t_2)\pb\omega^{2\overline{23}}.
\]
So
\[
\sigma_1=\pb^*G_{\pb}\mathcal{L}_{\phi_1}^{1,0}\sigma_0=\frac{-a_3(3t_3+2t_2)}{\sqrt{10}}(\omega^{2\overline{23}}-3\omega^{3\overline{13}}),
\]
where we used the fact that $(\ker\pb)^\bot\cap A^{1,2}(X)=\C\{\omega^{3\overline{23}}, \omega^{2\overline{23}}-3\omega^{3\overline{13}}\}$. Similarly, we have
\[
\mathcal{L}_{\phi_1}^{1,0}\sigma_1=\frac{-a_3(3t_3+2t_2)(9t_3+7t_2)}{\sqrt{10}}\omega^{1\overline{123}},\quad \mathcal{L}_{\phi_2}^{1,0}\sigma_0=-a_3A\omega^{1\overline{123}},
\]
which implies $\sigma_2=\pb^*G_{\pb}(\mathcal{L}_{\phi_2}^{1,0}\sigma_0+\mathcal{L}_{\phi_1}^{1,0}\sigma_1)=p_2(\omega^{2\overline{23}}-3\omega^{3\overline{13}})$, where $p_2$ is a polynomial of degree $2$ in $(t_2,t_3)$. Now because $\mathcal{L}_{\phi_2}^{1,0}(\omega^{2\overline{23}}-3\omega^{3\overline{13}})=3A\omega^{1\overline{123}}$, we have $\sigma_3=\pb^*G_{\pb}(\mathcal{L}_{\phi_3}^{1,0}\sigma_0+\mathcal{L}_{\phi_2}^{1,0}\sigma_1+\mathcal{L}_{\phi_1}^{1,0}\sigma_2)
=p_3(\omega^{2\overline{23}}-3\omega^{3\overline{13}})$ where $p_3$ is a polynomial of degree $3$ in $(t_2,t_3)$. Since $\phi_j=(3t_2+t_3)^{j-2}\phi_2$ for $j\geq 2$, this process can always be continued, we see that the deformation of $(1,2)$-forms are unobstructed.

\subsubsection{Examinations of results}
Combining these computations with Theorem \ref{mainthm}, we conclude that in this example the deformations of $(p,q)$-forms are unobstructed for any $(p,q)$. On the other hand, by applying Theorem \ref{mainthm}, we have $\bigoplus_{\lambda\geq 1}d_\lambda^{p,q}(X_t)=0$ for any small $t$ when
\[
(p,q)\neq (0,2), (1,2), (2,0), (3,0).
\]
Furthermore, by applying Corollary \ref{coro-to-mainthm2}, $(1)$, we get $\bigoplus_{\lambda\geq 1}d_\lambda^{p,q}(X_t)=0$ for any small $t$ whenever $(p,q)\neq (0,2), (1,2), (1,1), (2,1)$. So we conclude that
\begin{center}
$\bigoplus_{\lambda\geq 1}d_\lambda^{p,q}(X_t)=0$ for any small $t$ when $(p,q)\neq (0,2), (1,2)$.
\end{center}

\appendix
\section{Hodge decomposition for filtered forms}\label{appendix-A}
Let $X$ be a compact complex manifold of dimension $n$ equipped with a fixed Hermitian metric. Let us consider the following complex vector bundle on $X$:
\[
E_{\geq p}^{p+q}=\bigoplus_{\substack{i\geq p,\\i+j=p+q} }\bigoplus^n_{j=0}\left( (\wedge^iT^{*1,0}_X)\wedge(\wedge^jT^{*0,1}_X) \right)\subset \wedge^{p+q}T^*_{X,\C},
\]
where $T^{*1,0}_X$ is the holomorphic cotangent bundle of $X$ and $T^*_{X,\C}$ is the complexified cotangent bundle of $X$. Then the space of smooth sections of $E_{\geq p}^{p+q}$ is exactly the space of filtered forms (see Definition \ref{def-filtered-forms}), i.e.
\[
\Gamma(X,E_{\geq p}^{p+q})=F^pA^{p+q}(X).
\]
As is well known, we have
\[
d:\Gamma(X,E_{\geq p}^{p+q})\longrightarrow \Gamma(X,E_{\geq p}^{p+q+1}).
\]
In order to construct a convergent solution of the deformation equation, we need to find canonical solutions to the $d$-equation in the spaces $F^pA^{\bullet}(X)$ (see the proof of Theorem \ref{mainthm} and Theorem \ref{mainthm-2nd}). The usual $d$-Laplacian $\triangle_d=dd^*+d^*d$ is not suitable here because the image of $F^pA^{\bullet}(X)$ under the operator $d^*$ may lie outside $F^pA^{\bullet}(X)$. indeed, in general we only know $d^*:F^pA^{p+q}(X)\to F^{p-1}A^{p+q-1}(X)$. To remedy this, we set
\[
d^*_p=\Pi^{\geq p}d^*,
\]
where $\Pi^{\geq p}:A^{\bullet}(X)\to F^pA^{\bullet}(X)$ is the linear projection onto $F^pA^{\bullet}(X)$. Then
\begin{equation}
d^*_p:\Gamma(X,E_{\geq p}^{p+q+1})\longrightarrow \Gamma(X,E_{\geq p}^{p+q}).
\end{equation}
Set $\triangle_p=dd^*_p+d^*_pd$, then we have
\begin{equation}
\triangle_p:\Gamma(X,E_{\geq p}^{p+q})\longrightarrow \Gamma(X,E_{\geq p}^{p+q}).
\end{equation}
It can be checked (in the same way as in the usual case for $\triangle_d$) that $d^*_p$ is the formal adjoint of $d:F^pA^{p+q}(X)\to F^pA^{p+q+1}(X)$ and $\triangle_p$ is a formally self-adjoint, elliptic differential operator. Then by classic theory about such operators (see e.g. \cite[pp.289]{Dem12}), we have the following
\begin{proposition}\label{prop-hodge-for-filtered-forms}
For each pair $(p,q)$, there is a linear operator
\[
G_p:\Gamma(X,E_{\geq p}^{p+q})\longrightarrow \Gamma(X,E_{\geq p}^{p+q}),
\]
called the Green operator of $\triangle_p$, such that
\[
1=\mathcal{H}_p+\triangle_pG_p=\mathcal{H}_p+G_p\triangle_p,\quad\text{on}~\Gamma(X,E_{\geq p}^{p+q})=F^pA^{p+q}(X),
\]
where $\mathcal{H}_p:F^pA^{p+q}(X)\to \ker\triangle_p\cap F^pA^{p+q}(X)$ is the orthogonal projection onto the finite dimensional vector space $\ker\triangle_p\cap F^pA^{p+q}(X)$. Equivalently, there is an orthogonal decomposition
\begin{equation}\label{eq-Hodge-orthogonal decomposition}
F^pA^{p+q}(X)=\ker\triangle_p\cap F^pA^{p+q}(X)\oplus d\left(F^pA^{p+q-1}(X)\right)\oplus d^*_p\left(F^pA^{p+q+1}(X)\right).
\end{equation}
\end{proposition}

\begin{proposition}\label{prop-canonicalsolution-d-eq}
Let $r\geq0$. For any $y\in d\left(F^{p+r}A^{p+q}(X)\right)\subset F^{p+r}A^{p+q+1}(X)$, there is a canonical solution $x\in F^{p+r}A^{p+q}(X)$ to the equation $dx=y$ given by
\[
x=d^*_{p+r}G_{p+r}y.
\]
In particular, for any $\hat{\alpha}_0^{p,q}\in\ker\pb$, if there is a $\alpha\in F^pA^{p+q}(X)\cap\ker d$ such that $\alpha^{p,q}=\hat{\alpha}_0^{p,q}$, then $\alpha$ can be taken as the following canonical form:
\[
\alpha=(1-d^*_{p+1}G_{p+1}\p)\hat{\alpha}_0^{p,q}.
\]
\end{proposition}
\begin{proof}By Proposition \ref{prop-hodge-for-filtered-forms}, $\mathcal{H}_{p+r}y=0$ and
\[
y=\mathcal{H}_{p+r}y+\triangle_{p+r}G_{p+r}y=dd_{p+r}^*G_{p+r}y+d_{p+r}^*dG_{p+r}y
\]
which implies
\[
y-dd_{p+r}^*G_{p+r}y=d_{p+r}^*dG_{p+r}y\in d\left(F^{p+r}A^{p+q}(X)\right)\cap d_{p+r}^*\left(F^{p+r}A^{p+q+2}(X)\right)=0,
\]
where we have used the Hodge decomposition \eqref{eq-Hodge-orthogonal decomposition}. Therefore, we have
\[
d(d_{p+r}^*G_{p+r}y)=y.
\]

The second statement follows immediately by noting that in this particular case we have $-\p\hat{\alpha}_0^{p,q}=d(\alpha-\hat{\alpha}_0^{p,q})$ and $\alpha-\hat{\alpha}_0^{p,q}\in F^{p+1}A^{p+q}(X)$.
\end{proof}

\section{Various conditions about degenerations of $(E_r^{\bullet,\bullet}(X),d_r^{\bullet,\bullet})$}\label{appendix-B}
Let $(E_r^{\bullet,\bullet}(X),d_r^{\bullet,\bullet})$ be the Fr\"olicher spectral sequence of a complex manifold $X$, we say the Fr\"olicher spectral sequence \textit{degenerates at the $r$-th page} if
\[
E_r^{\bullet,\bullet}(X)\cong E_\infty^{\bullet,\bullet}(X).
\]
\subsection{Degenerations of $(E_r^{\bullet,\bullet}(X),d_r^{\bullet,\bullet})$ and $d^{p,q}_r=0$}
It is well-known that $E_1^{p,q}(X)\cong E_2^{p,q}(X)$ if and only if
\[
\forall \alpha^{p,q}, \alpha^{p-1,q}\in\ker\pb,\quad \p\alpha^{p,q}, \p\alpha^{p-1,q}\in\im\pb.
\]
From the commutative diagram \eqref{eq-drXtodr} and the fact $E_{r+1}^{p,q}(X)\cong \ker d_r^{p,q}/\im d_r^{p-r,q+r-1}$, we see that $E_2^{p,q}(X)\cong E_3^{p,q}(X)$ if and only if
\begin{itemize}
\item for any $\alpha^{p,q}$ and $\alpha^{p+1,q-1}$ with $0=\pb\alpha^{p,q}=\p\alpha^{p,q}+\pb\alpha^{p+1,q-1}$ there exists $\beta^{p+1,q-1}$ and $\beta^{p+2,q-2}$ such that
\[
\p\alpha^{p+1,q-1}=\p\beta^{p+1,q-1}+\pb\beta^{p+2,q-2};
\]
\item for any $\alpha^{p-2,q+1}$ and $\alpha^{p-1,q}$ with $0=\pb\alpha^{p-2,q+1}=\p\alpha^{p-2,q+1}+\pb\alpha^{p-1,q}$ there exists $\beta^{p-2,q+1}$ and $\beta^{p-1,q}$
 such that
\[
\p\alpha^{p-1,q}=\p\beta^{p-2,q+1}+\pb\beta^{p-1,q}.
\]
\end{itemize}
More generally, by using \eqref{eq-drXtodr} we have the following
\begin{proposition}\label{prop-Er=Er+1}For $r\geq 2$, $E_r^{p,q}(X)\cong E_{r+1}^{p,q}(X)$ if and only if
\begin{itemize}
\item the homomorphism
\[
d_r^{p,q}:\wt{Z}_r^{p,q}(X)/\wt{B}_r^{p,q}(X)\to \wt{Z}_r^{p+r,q-r+1}(X)/\wt{B}_r^{p+r,q-r+1}(X)
\]
is zero i.e. for any $\alpha^{p,q}$ with
\[
0=\pb\alpha^{p,q}=\p\alpha^{p,q}+\pb\alpha^{p+1,q-1}=\cdots=\p\alpha^{p+r-2,q-r+2}+\pb\alpha^{p+r-1,q-r+1}
\]
for some $\alpha^{p+1,q-1}, \cdots, \alpha^{p+r-1,q-r+1}$, there exists $\beta^{p+r,q-r}, \cdots, \beta^{p+1,q-1}$ such that
\begin{equation}\label{eq-dralpha=0}
      \left\{
      \begin{array}{ll}
       \p\alpha^{p+r-1,q-r+1}&=\p\beta^{p+r-1,q-r+1}+\pb\beta^{p+r,q-r},\\
       0&=\p\beta^{p+r-2,q-r+2}+\pb\beta^{p+r-1,q-r+1},\\
        &\cdots\\
       0&=\p\beta^{p+1,q-1}+\pb\beta^{p+2,q-2},\\
       0&=\pb\beta^{p+1,q-1}.
      \end{array} \right.
\end{equation}
\item the homomorphism
\[
d_r^{p-r,q+r-1}:\wt{Z}_r^{p-r,q+r-1}(X)/\wt{B}_r^{p-r,q+r-1}(X)\to \wt{Z}_r^{p,q}(X)/\wt{B}_r^{p,q}(X)
\]
is zero i.e. for any $\alpha^{p-r,q+r-1}$ with
\[
0=\pb\alpha^{p-r,q+r-1}=\p\alpha^{p-r,q+r-1}+\pb\alpha^{p-r+1,q+r-2}=\cdots=\p\alpha^{p-2,q+1}+\pb\alpha^{p-1,q}
\]
for some $\alpha^{p-r+1,q+r-2}, \cdots, \alpha^{p-1,q}$, there exists $\beta^{p,q-1}, \cdots, \beta^{p-r+1,q+r-2}$ such that
\begin{equation}
      \left\{
      \begin{array}{ll}
       \p\alpha^{p-1,q}&=\p\beta^{p-1,q}+\pb\beta^{p,q-1},\\
       0&=\p\beta^{p-2,q+1}+\pb\beta^{p-1,q},\\
        &\cdots\\
       0&=\p\beta^{p-r+1,q+r-2}+\pb\beta^{p-r+2,q+r-3},\\
       0&=\pb\beta^{p-r+1,q+r-2}.
      \end{array} \right.
\end{equation}
\end{itemize}
\end{proposition}

\begin{proposition}\label{prop-kerim-dr}
The following statements holds:
\begin{enumerate}
\item For any $\alpha^{p,q}\in \wt{Z}^{p,q}_r(X)$, we have
\[
d_r^{p,q}\alpha^{p,q}=0\in \wt{Z}_r^{p+r,q-r+1}(X)/\wt{B}_r^{p+r,q-r+1}(X)
\]
if and only if $\alpha^{p,q}\in \wt{Z}_{r+1}^{p,q}(X)$, that is, there exists $\alpha'^{p+1,q-1},\cdots, \alpha'^{p+r,q-r}$ such that
\[
0=\pb\alpha^{p,q}=\p\alpha^{p,q}+\pb\alpha'^{p+1,q-1}=\cdots=\p\alpha'^{p+r-1,q-r+1}+\pb\alpha'^{p+r,q-r}.
\]
In particular, we have
\[
\ker d_r^{p,q}\cong \frac{\wt{Z}_{r+1}^{p,q}(X)}{\wt{B}_{r}^{p,q}(X)}\quad\text{and}\quad\im d_r^{p,q}\cong \frac{\wt{Z}_{r}^{p,q}(X)}{\wt{Z}_{r+1}^{p,q}(X)}.
\]
\item We have
\[
\im d_r^{p-r,q+r-1}\cong \frac{\wt{B}_{r+1}^{p,q}(X)}{\wt{B}_{r}^{p,q}(X)}~.
\]
\end{enumerate}
\end{proposition}
\begin{proof}$(1)$ In fact, assume $d_r^{p,q}\alpha^{p,q}\in \wt{B}_r^{p+r,q-r+1}(X)$. Then it follows from \eqref{eq-dralpha=0} that
\begin{align*}
0&=\pb\alpha^{p,q}=\p\alpha^{p,q}+\pb(\alpha^{p+1,q-1}-\beta^{p+1,q-1}),\\
&\cdots\\
0&=\p(\alpha^{p+r-2,q-r+2}-\beta^{p+r-2,q-r+2})+\pb(\alpha^{p+r-1,q-r+1}-\beta^{p+r-1,q-r+1}),
\end{align*}
and
\[
\p(\alpha^{p+r-1,q-r+1}-\beta^{p+r-1,q-r+1})=\pb\beta^{p+r,q-r}.
\]

Conversely, assume there exists $\alpha'^{p+1,q-1},\cdots, \alpha'^{p+r,q-r}$ such that
\[
0=\pb\alpha^{p,q}=\p\alpha^{p,q}+\pb\alpha'^{p+1,q-1}=\cdots=\p\alpha'^{p+r-1,q-r+1}+\pb\alpha'^{p+r,q-r}.
\]
We may set $\beta^{p+r-1,q-r+1}=-\alpha'^{p+r,q-r}$ and $\beta^{p+r-i,q-r+i}=0$ for all $i\geq 1$, then \eqref{eq-dralpha=0} is satisfied
which implies $d_r^{p,q}\alpha^{p,q}\in \wt{B}_r^{p+r,q-r+1}$.

Now, $\im d_r^{p,q}\cong \frac{\wt{Z}_{r}^{p,q}(X)}{\wt{Z}_{r+1}^{p,q}(X)}$ follows immediately from $\ker d_r^{p,q}\cong \frac{\wt{Z}_{r+1}^{p,q}(X)}{\wt{B}_{r}^{p,q}(X)}$ since
\[
\im d_r^{p,q}\cong \frac{E_r^{p,q}(X)}{\ker d_r^{p,q}}\cong\frac{\wt{Z}_{r}^{p,q}(X)/\wt{B}_{r}^{p,q}(X)}{\wt{Z}_{r+1}^{p,q}(X)/\wt{B}_{r}^{p,q}(X)}.
\]

$(2)$ Set
\begin{align*}
\begin{split}
  C^{p,q}_r:=\{&\p\alpha^{p-1,q}\mid 0=\p\alpha^{p-2,q+1}+\pb\alpha^{p-1,q}=\\
  &\cdots=\p\alpha^{p-r,q+r-1}+\pb\alpha^{p-r+1,q+r-2}=\pb\alpha^{p-r,q+r-1} \}.
 \end{split}
\end{align*}
It follows from the commutative diagram \eqref{eq-drXtodr} that (we write $\wt{B}_{r}^{p,q}=\wt{B}_{r}^{p,q}(X)$ for convenience)
\[
\im d_r^{p-r,q+r-1}=\frac{C^{p,q}_r }{C^{p,q}_r\cap \wt{B}_{r}^{p,q}}\cong \frac{C^{p,q}_r+\wt{B}_{r}^{p,q} }{ \wt{B}_{r}^{p,q} }=\frac{\wt{B}_{r+1}^{p,q}}{\wt{B}_{r}^{p,q}}~,
\]
where we have used the fact $C^{p,q}_r+\wt{B}_{r}^{p,q}=\wt{B}_{r+1}^{p,q}$: indeed, assume $\p\alpha^{p-1,q}\in C^{p,q}_r$ and $\p\beta^{p-1,q}+\pb\beta^{p,q-1}\in \wt{B}_{r}^{p,q}$, then $\p(\alpha^{p-1,q}+\beta^{p-1,q})+\pb\beta^{p,q-1}\in\wt{B}_{r+1}^{p,q}$
\begin{equation*}
      \left\{
      \begin{array}{ll}
       0&=\p(\alpha^{p-2,q+1}+\beta^{p-2,q+1} )+\pb(\alpha^{p-1,q}+\beta^{p-1,q} ),\\
        &\cdots\\
       0&=\p(\alpha^{p-r+1,q+r-2}+\beta^{p-r+1,q+r-2} )+\pb(\alpha^{p-r+2,q+r-3}+\beta^{p-r+2,q+r-3} ),\\
       0&=\p\alpha^{p-r,q+r-1} +\pb(\alpha^{p-r+1,q+r-2}+\beta^{p-r+1,q+r-2} ),\\
       0&=\pb\alpha^{p-r,q+r-1},
      \end{array} \right.
\end{equation*}
which implies $C^{p,q}_r+\wt{B}_{r}^{p,q}\subseteq \wt{B}_{r+1}^{p,q}$. The inverse inclusion $C^{p,q}_r+\wt{B}_{r}^{p,q}\supseteq \wt{B}_{r+1}^{p,q}$ follows from $\wt{B}_{r+1}^{p,q}\subseteq C^{p,q}_r+\wt{B}_{1}^{p,q}$ and $\wt{B}_{1}^{p,q}\subseteq \wt{B}_{r}^{p,q}$.
\end{proof}

\subsection{The condition $\bigoplus_{\lambda\geq1}d^{p,q}_\lambda=0$}
\begin{proposition}\label{prop-dr=0}
For any $(p,q)\in\mathbb{N}\times\mathbb{N}$, the following statements are equivalent:
\begin{enumerate}
  \item $d_\lambda^{p,q}(X)=0$ for all $\lambda\geq 1$;
  \item For any given $\alpha^{p,q}\in \ker\pb\cap A^{p,q}(X)$ there exists $\alpha\in \ker d\cap F^pA^{p+q}(X)$ such that $\Pi^{p,q}\alpha=\alpha^{p,q}$.
\end{enumerate}
\end{proposition}
\begin{proof}$(1)\Rightarrow(2)$: Assume $d_\lambda^{p,q}=0$ for all $\lambda\geq 1$ and $\alpha^{p,q}\in \ker\pb\cap A^{p,q}(X)$. By Proposition \ref{prop-kerim-dr} we have
\[
\alpha^{p,q}\in\wt{Z}^{p,q}_1(X)=\wt{Z}^{p,q}_2(X)=\cdots=\wt{Z}^{p,q}_\infty(X).
\]
In particular, $\alpha^{p,q}\in\wt{Z}^{p,q}_\infty(X)$ implies there exist $\alpha'^{p+1,q-1},\cdots, \alpha'^{p+q,0}$ such that
\[
0=\pb\alpha^{p,q}=\p\alpha^{p,q}+\pb\alpha'^{p+1,q-1}=\cdots
=\p\alpha'^{p+q-1,1}+\pb\alpha'^{p+q,0}=\p\alpha'^{p+q,0}.
\]
Set $\alpha=\alpha^{p,q}+\alpha'^{p+1,q-1}+\cdots+\alpha'^{p+q,0}$, then $\alpha\in \ker d\cap F^pA^{p+q}(X)$ such that $\Pi^{p,q}\alpha=\alpha^{p,q}$.

$(2)\Rightarrow(1)$: Conversely, assume for any given $\alpha^{p,q}\in \ker\pb\cap A^{p,q}(X)$ there exist $\alpha\in \ker d\cap F^pA^{p+q}(X)$ such that $\Pi^{p,q}\alpha=\alpha^{p,q}$. This holds if and only if
\[
\wt{Z}_1^{p,q}(X)=\wt{Z}_{\infty}^{p,q}(X).
\]
Because $\wt{Z}_{\lambda+1}^{p,q}(X)\subseteq \wt{Z}_\lambda^{p,q}(X)$ for any $\lambda\geq 1$, we get
\[
\wt{Z}_1^{p,q}(X)=\wt{Z}_2^{p,q}(X)=\cdots=\wt{Z}_{\infty}^{p,q}(X).
\]
But according to Proposition \ref{prop-kerim-dr}, $d_\lambda^{p,q}=0$ if and only if $\wt{Z}_{\lambda+1}^{p,q}(X)=\wt{Z}_\lambda^{p,q}(X)$. Hence $d_\lambda^{p,q}=0$ for all $\lambda\geq 1$.
\end{proof}

\subsection{Degenerations of $(E_r^{\bullet,\bullet}(X),d_r^{\bullet,\bullet})$ and surjectivity of $\Pi_{j,k}^{p,\ker}$}
For any fixed natural numbers $k,j\in\mathbb{Z}^+$ with $k\geq j$, let us consider the following natural projection operators,
\begin{equation}\label{eq-Pi-k,j-p,ker}
\Pi_{k,j}^{p,\ker}:\ker d_{\Pi_k}\cap F^p_{k}A^{p+q}(X)\longrightarrow\ker\big(d_{\Pi_j}:F^p_{j}A^{p+q}(X)\to F^p_{j}A^{p+q+1}(X)\big),
\end{equation}
where
\begin{align*}
d_{\Pi_j}:&=d-\p\Pi^{p+j-1,q-j+1},\\
F^p_{j}A^{p+q}(X):&=A^{p,q}(X)\oplus\cdots\oplus A^{p+j-1,q-j+1}(X).
\end{align*}
For $k=\infty$, there are also natural projections
\begin{equation}\label{eq-Pi-infty,j-p,ker}
\Pi_{\infty,j}^{p,\ker}:\ker d\cap F^pA^{p+q}(X)\longrightarrow\ker\big(d_{\Pi_j}:F^p_{j}A^{p+q}(X)\to F^p_{j}A^{p+q+1}(X)\big).
\end{equation}
Note that
\begin{align*}
&\ker d\cap F^pA^{p+q}(X)\supseteq\cdots\supseteq\ker d_{\Pi_k}\cap F^p_{k}A^{p+q}(X)\supseteq\cdots\\
&\ker d_{\Pi_j}\cap F^p_{j}A^{p+q}(X)\supseteq\cdots \supseteq\ker d_{\Pi_1}\cap F^p_{1}A^{p+q}(X)=\ker\pb\cap A^{p,q}(X).
\end{align*}

\begin{proposition}\label{prop-Pi-k,j-p,ker-surj} Let $k,j\in\mathbb{Z}^+\cup\infty$ with $k\geq j$. Then the projection operator $\Pi_{k,j}^{p,\ker}$ defined as above is surjective
if and only if $\wt{Z}^{p+i,q-i}_{j-i}(X)=\wt{Z}^{p+i,q-i}_{k-i}(X)$ holds for any $0\leq i\leq j-1$.
\end{proposition}
\begin{proof}We will only give the proof when $k\in\mathbb{N}$ because the case $k=\infty$ can be proved in the same way. Let $k,j\in\mathbb{N}$ with $k\geq j$.

Assume $\wt{Z}^{p+i,q-i}_{j-i}(X)=\wt{Z}^{p+i,q-i}_{k-i}(X)$ holds for any $0\leq i\leq j-1$. Let $\alpha^{p,q}+\cdots+\alpha^{p+j-1,q-j+1}\in\ker d_{\Pi_j}\cap F^p_{j}A^{p+q}(X)$, then
\[
0=\pb\alpha^{p,q}=\p\alpha^{p,q}+\pb\alpha^{p+1,q-1}=\cdots=\p\alpha^{p+j-2,q-j+2}+\pb\alpha^{p+j-1,q-j+1},
\]
in particular, $\alpha^{p,q}\in\wt{Z}^{p,q}_{j}(X)=\wt{Z}^{p,q}_{k}(X)$. There exist $\alpha_1^{p+1,q-1},\cdots,\alpha_1^{p+k-1,q-k+1}$ such that
\begin{equation}\label{eq-alpha1}
0=\pb\alpha^{p,q}=\p\alpha^{p,q}+\pb\alpha_1^{p+1,q-1}=\cdots=\p\alpha_1^{p+k-2,q-k+2}+\pb\alpha_1^{p+k-1,q-k+1}.
\end{equation}
Combining these equalities, we get
\begin{equation}\label{eq-alpha-alpha1}
\begin{split}
0=&\pb(\alpha^{p+1,q-1}-\alpha_1^{p+1,q-1}),\\
0=&\p(\alpha^{p+1,q-1}-\alpha_1^{p+1,q-1})+\pb(\alpha^{p+2,q-2}-\alpha_1^{p+2,q-2}),\\
&\vdots\\
0=&\p(\alpha^{p+j-2,q-j+2}-\alpha_1^{p+j-2,q-j+2})+\pb(\alpha^{p+j-1,q-j+1}-\alpha_1^{p+j-1,q-j+1}),
\end{split}
\end{equation}
which implies $\alpha^{p+1,q-1}-\alpha_1^{p+1,q-1}\in\wt{Z}^{p+1,q-1}_{j-1}(X)=\wt{Z}^{p+1,q-1}_{k-1}(X)$. There exist
\[
\alpha_2^{p+2,q-2},\cdots,\alpha_2^{p+k-1,q-k+1}
\]
such that
\begin{equation}\label{eq-alpha2}
\begin{split}
0&=\pb(\alpha^{p+1,q-1}-\alpha_1^{p+1,q-1}),\\
0&=\p(\alpha^{p+1,q-1}-\alpha_1^{p+1,q-1})+\pb\alpha_2^{p+2,q-2},\\
&\vdots\\
0&=\p\alpha_2^{p+k-2,q-k+2}+\pb\alpha_2^{p+k-1,q-k+1}.
\end{split}
\end{equation}
By using \eqref{eq-alpha1} and \eqref{eq-alpha2}, we get
\begin{equation}\label{eq-alpha1+alpha2}
\begin{split}
0&=\pb\alpha^{p,q}=\p\alpha^{p,q}+\pb\alpha^{p+1,q-1},\\
0&=\p\alpha^{p+1,q-1}+\pb(\alpha_1+\alpha_2)^{p+2,q-2},\\
&\vdots\\
0&=\p(\alpha_1+\alpha_2)^{p+k-2,q-k+2}+\pb(\alpha_1+\alpha_2)^{p+k-1,q-k+1}.
\end{split}
\end{equation}
Moreover, it follows from \eqref{eq-alpha-alpha1} and \eqref{eq-alpha2} that
\begin{align*}
0&=\pb(\alpha-\alpha_1-\alpha_2)^{p+2,q-2},\\
0&=\p(\alpha-\alpha_1-\alpha_2)^{p+2,q-2}+\pb(\alpha-\alpha_1-\alpha_2)^{p+3,q-3},\\
&\vdots\\
0&=\p(\alpha-\alpha_1-\alpha_2)^{p+j-2,q-j+2}+\pb(\alpha-\alpha_1-\alpha_2)^{p+j-1,q-j+1},
\end{align*}
which implies $(\alpha-\alpha_1-\alpha_2)^{p+2,q-2}\in\wt{Z}^{p+2,q-2}_{j-2}(X)=\wt{Z}^{p+2,q-2}_{k-2}(X)$. Continuing this process, by using $\wt{Z}^{p+i,q-i}_{j-i}(X)=\wt{Z}^{p+i,q-i}_{k-i}(X)$ for $i=3,4,\cdots,j-1$, there exist
\[
\alpha_i^{p+i,q-i},\cdots,\alpha_i^{p+k-1,q-k+1},\quad 1\leq i\leq j,
\]
such that
\begin{equation}\label{eq-alpha1++j}
\begin{split}
0&=\pb\alpha^{p,q}=\p\alpha^{p,q}+\pb\alpha^{p+1,q-1}=\cdots=\p\alpha^{p+j-2,q-j+2}+\pb\alpha^{p+j-1,q-j+1},\\
0&=\p\alpha^{p+j-1,q-j+1}+\pb(\sum_{i=1}^{j}\alpha_i)^{p+2,q-2},\\
&\vdots\\
0&=\p(\sum_{i=1}^{j}\alpha_i)^{p+k-2,q-k+2}+\pb(\sum_{i=1}^{j}\alpha_i)^{p+k-1,q-k+1}.
\end{split}
\end{equation}
and for any $1\leq l\leq j-1$ we have
\begin{equation}\label{eq-alphal}
\begin{split}
0&=\pb(\alpha-\sum_{i=1}^{l}\alpha_i)^{p+l,q-l},\\
0&=\p(\alpha-\sum_{i=1}^{l}\alpha_i)^{p+l,q-l}+\pb\alpha_{l+1}^{p+l+1,q-l-1},\\
&\vdots\\
0&=\p\alpha_{l+1}^{p+k-2,q-k+2}+\pb\alpha_{l+1}^{p+k-1,q-k+1},
\end{split}
\end{equation}
and for any $1\leq l\leq j-2$,
\begin{equation}\label{eq-lem-l2}
\begin{split}
0&=\pb(\alpha-\sum_{i=1}^{l+1}\alpha_i)^{p+l+1,q-l-1},\\
0&=\p(\alpha-\sum_{i=1}^{l+1}\alpha_i)^{p+l+1,q-l-1}+\pb(\alpha-\sum_{i=1}^{l+1}\alpha_i)^{p+l+2,q-l-2},\\
&\vdots\\
0&=\p(\alpha-\sum_{i=1}^{l+1}\alpha_i)^{p+j-2,q-j+2}+\pb(\alpha-\sum_{i=1}^{l+1}\alpha_i)^{p+j-1,q-j+1}.
\end{split}
\end{equation}
Therefore, we see from \eqref{eq-alpha1++j} that
\[
\sum_{i=0}^{j-1}\alpha^{p+i,q-i}+\sum_{m=j}^{k-1}\sum_{i=1}^{j}\alpha_i^{p+m,q-m}\in\ker d_{\Pi_k}\cap F^p_{k}A^{p+q}(X)
\]
and
\[
\Pi_{k,j}^{p,\ker}(\sum_{i=0}^{j-1}\alpha^{p+i,q-i}+\sum_{m=j}^{k-1}\sum_{i=1}^{j}\alpha_i^{p+m,q-m})=\sum_{i=0}^{j-1}\alpha^{p+i,q-i}\in\ker d_{\Pi_j}\cap F^p_{j}A^{p+q}(X).
\]
This shows $\Pi_{k,j}^{p,\ker}$ is surjective.

Conversely, assume $\Pi_{k,j}^{p,\ker}$ is surjective. For any $0\leq i\leq j-1$, let $\alpha^{p+i,q-i}\in\wt{Z}^{p+i,q-i}_{j-i}(X)$, we need to show $\alpha^{p+i,q-i}\in\wt{Z}^{p+i,q-i}_{k-i}(X)$. Indeed, by the definition of $\wt{Z}^{p+i,q-i}_{j-i}(X)$ there exist
\[
\alpha_0^{p+i+1,q-i-1},\cdots,\alpha_0^{p+j-1,q-j+1}
\]
such that (we let $\alpha^{p,q}_0=\cdots=\alpha^{p+i-1,q-i+1}_0=0$)
\begin{align*}
0&=\pb\alpha^{p,q}_0=\p\alpha^{p,q}_0+\pb\alpha^{p+1,q-1}_0=\cdots=\p\alpha^{p+i-1,q-i+1}_0+\pb\alpha^{p+i,q-i}\\
&=\p\alpha^{p+i,q-i}+\pb\alpha^{p+i+1,q-i-1}_0=\cdots=\p\alpha^{p+j-2,q-j+2}_0+\pb\alpha^{p+j-1,q-j+1}_0,
\end{align*}
which implies
\[
\alpha^{p+i,q-i}+\alpha_0^{p+i+1,q-i-1}+\cdots+\alpha_0^{p+j-1,q-j+1}\in\ker d_{\Pi_j}\cap F^p_{j}A^{p+q}(X).
\]
By using the surjectivity of $\Pi_{k,j}^{p,\ker}$ there exist $\alpha_0^{p+j,q-j},\cdots,\alpha_0^{p+k,q-k}$ such that
\[
\alpha^{p+i,q-i}+\alpha_0^{p+i+1,q-i-1}+\cdots+\alpha_0^{p+k,q-k}\in\ker d_{\Pi_k}\cap F^p_{k}A^{p+q}(X).
\]
In particular, $\alpha^{p+i,q-i}\in\wt{Z}^{p+i,q-i}_{k-i}(X)$.
\end{proof}
The following Corollary is a generalization of Proposition \ref{prop-dr=0}.
\begin{corollary}\label{coro-Pi-p+r,r,ker-surj}
Let $r\geq1$ be an integer. The following conditions are equivalent:
\begin{enumerate}
  \item The projection $\Pi_{\infty,r}^{p,\ker}$ defined by $\eqref{eq-Pi-infty,j-p,ker}$ is surjective;
  \item We have
  \[
\bigoplus_{ \substack{\lambda\geq r-i,\\ 0\leq i\leq r-1} }d_\lambda^{p+i,q-i}(X)=0;
\]
\end{enumerate}
\end{corollary}
\begin{proof}This follows immediately from Proposition \ref{prop-Pi-k,j-p,ker-surj} and Proposition \ref{prop-kerim-dr}.
\end{proof}
Now, we examine the following condition
\[
   F^{p+1}A^{p+1+q}(X)\cap dA^{p+q}(X)= dF^{p+1}A^{p+q}(X),
\]
in other words, the de Rham differential $d$ is strict with respect to the filtration $F^{\blacktriangledown}A^\bullet(X)$ in degrees $\blacktriangledown=p+1$ and $\bullet=p+1+q$ (see \cite{Fri,CS22} for related discussions).
\begin{corollary}\label{coro-Fp+1capdA=dFp+1}
Let $r\geq1$ be an integer. The following conditions are equivalent:
\begin{enumerate}
  \item The projection $\Pi^{0,\ker}_{\infty,p+1}$
\[
\ker d\cap A^{p+q}(X)\longrightarrow\ker\big(d-\p\Pi^{p,q}:\bar{F}^qA^{p+q}(X)\to \bar{F}^qA^{p+q+1}(X)\big),
\]
is surjective, where $\bar{F}^qA^{p+q}(X)=F^0_{p+1}A^{p+q}(X)$ and $\bar{F}^qA^{p+q}(X)=F^0_{p+2}A^{p+q+1}(X)$;
  \item We have
  \[
\bigoplus_{ \substack{\lambda\geq i+1,\\ i\geq 0} }d_\lambda^{p-i,q+i}(X)=0;
\]
  \item We have
\begin{equation}\label{Fp+1dA=dFp+1}
  F^{p+1}A^{p+q+1}(X)\cap dA^{p+q}(X)= dF^{p+1}A^{p+q}(X).
\end{equation}
\end{enumerate}
\end{corollary}
\begin{proof}
For any $d\sum_{i\in\Z}\alpha^{p+i,q-i}\in dA^{p+q}(X)$, it is clear that
\[
d\sum_{i\in\Z}\alpha^{p+i,q-i}\in F^{p+1}A^{p+q+1}(X)\Longleftrightarrow d(\sum_{i\leq0}\alpha^{p+i,q-i})-\p\alpha^{p,q}=0,
\]
i.e. $\sum_{i\leq0}\alpha^{p+i,q-i}\in \ker\big(d-\p\Pi^{p,q}\big)\cap\bar{F}^qA^{p+q}(X)$.
It follows that for any $d\sum_{i\in\Z}\alpha^{p+i,q-i}\in F^{p+1}A^{p+q+1}(X)\cap dA^{p+q}(X)$, we have $d\sum_{i\in\Z}\alpha^{p+i,q-i}\in dF^{p+1}A^{p+q}(X)$ iff there exists $\beta\in F^{p+1}A^{p+q}(X)$ such that
\[
d\Big(\sum_{i\leq0}\alpha^{p+i,q-i}+\sum_{i\geq1}(\alpha-\beta)^{p+i,q-i}\Big)=0.
\]
This shows $(3)$ implies $(1)$. Conversely, for any $d\sum_{i\in\Z}\alpha^{p+i,q-i}\in F^{p+1}A^{p+q+1}(X)\cap dA^{p+q}(X)$, it follows from $(1)$ that there is a $\gamma\in F^{p+1}A^{p+q}(X)$ such that
\[
d\Big(\sum_{i\leq0}\alpha^{p+i,q-i}-\gamma\Big)=0~\text{or}~
d\Big(\sum_{i\leq0}\alpha^{p+i,q-i}\Big)=d\gamma.
\]
The equivalence of $(1)$ and $(2)$ follows directly from Corollary \ref{coro-Pi-p+r,r,ker-surj}.
\end{proof}

\bibliographystyle{alpha}
\bibliography{reference}
\end{document}